\documentclass{article}

\usepackage{microtype}
\usepackage{graphicx}
\usepackage{subfigure}
\usepackage{booktabs}

\usepackage{amsthm}
\usepackage{amsmath, amssymb}

\usepackage{mathtools}
\DeclarePairedDelimiter{\ceil}{\lceil}{\rceil}

\usepackage{subfigure}

\usepackage{hyperref}

\usepackage[accepted]{icml2021}

\usepackage{enumitem}
\newcommand{\less}{\leqslant}
\newcommand{\gre}{\geqslant}

\newcommand{\argmin}{\mathop{\rm argmin}}

\newcommand{\defn}{\ensuremath{: \, =}}
\newcommand{\real}{\mathbb{R}}

\newcommand{\wtilde}{\widetilde}
\newcommand{\vne}{v_\mathrm{ne}}
\newcommand{\vnsk}{v_\mathrm{nsk}}
\newcommand{\xne}{x_\mathrm{ne}}
\newcommand{\xnsk}{x_\mathrm{nsk}}
\newcommand{\gx}{\nabla f(x)}
\newcommand{\newtondec}{\lambda_f(x)}
\newcommand{\appnewtondec}{\wtilde \lambda_f(x)}
\newcommand{\mbR}{\mathbb{R}}

\newcommand{\lra}[1]{\left\langle #1\right\rangle}

\newtheorem{theo}{Theorem}[section]

\newtheorem{lem}{Lemma}[section]
\newtheorem{prop}{Proposition}[section]
\newtheorem{cor}{Corollary}[section]
\newtheorem{nota}{Notation}[section]
\newtheorem{de}{Definition}[section]
\newtheorem{exa}{Example}[section]
\newtheorem{as}{Assumption}[section]
\newtheorem{alg}{Algorithm}[section]
\newcommand{\btheo}{\begin{theo}}
\newcommand{\bde}{\begin{de}}
\newcommand{\ble}{\begin{lem}}
\newcommand{\bpr}{\begin{prop}}
\newcommand{\bno}{\begin{nota}}
\newcommand{\bex}{\begin{exa}}
\newcommand{\bcor}{\begin{cor}}
\newcommand{\spro}{\begin{proof}}
\newcommand{\bas}{\begin{as}}
\newcommand{\balg}{\begin{alg}}
\newcommand{\etheo}{\end{theo}}
\newcommand{\ede}{\end{de}}
\newcommand{\ele}{\end{lem}}
\newcommand{\epr}{\end{prop}}
\newcommand{\eno}{\end{nota}}
\newcommand{\eex}{\end{exa}}
\newcommand{\ecor}{\end{cor}}
\newcommand{\fpro}{\end{proof}}
\newcommand{\eas}{\end{as}}
\newcommand{\ealg}{\end{alg}}
\newtheorem{theos}{Theorem}
\newtheorem{props}{Proposition}
\newtheorem{lems}{Lemma}
\newtheorem{cors}{Corollary}
\newtheorem{exas}{Example}
\newtheorem{algs}{Algorithm}
\newtheorem{asss}{Assumption}
\newtheorem{defns}{Definition}
\newcommand{\btheos}{\begin{theos}}
\newcommand{\etheos}{\end{theos}}
\newcommand{\bprops}{\begin{props}}
\newcommand{\eprops}{\end{props}}
\newcommand{\bdes}{\begin{defns}}
\newcommand{\edes}{\end{defns}}
\newcommand{\blems}{\begin{lems}}
\newcommand{\elems}{\end{lems}}
\newcommand{\bcors}{\begin{cors}}
\newcommand{\ecors}{\end{cors}}
\newcommand{\bexs}{\begin{exas}}
\newcommand{\eexs}{\end{exas}}
\newcommand{\balgs}{\begin{algs}}
\newcommand{\ealgs}{\end{algs}}
\newcommand{\bass}{\begin{asss}}
\newcommand{\eass}{\end{asss}}

\icmltitlerunning{Adaptive Newton Sketch}

\begin{document}

\twocolumn[
\icmltitle{Adaptive Newton Sketch: Linear-time Optimization with Quadratic Convergence and Effective Hessian Dimensionality}

\icmlsetsymbol{equal}{*}

\begin{icmlauthorlist}
\icmlauthor{Jonathan Lacotte}{stanford}
\icmlauthor{Yifei Wang}{stanford}
\icmlauthor{Mert Pilanci}{stanford}
\end{icmlauthorlist}

\icmlaffiliation{stanford}{Department of Electrical Engineering, Stanford University}

\icmlcorrespondingauthor{Jonathan Lacotte}{lacotte@stanford.edu}
\icmlcorrespondingauthor{Yifei Wang}{wangyf18@stanford.edu}

\icmlkeywords{Newton method, Sketching, Quadratic convergence}

\vskip 0.3in
]

\printAffiliationsAndNotice{}

\begin{abstract}
We propose a randomized algorithm with quadratic convergence rate for convex optimization problems with a self-concordant, composite, strongly convex objective function. Our method is based on performing an approximate Newton step using a random projection of the Hessian. Our first contribution is to show that, at each iteration, the embedding dimension (or sketch size) can be as small as the effective dimension of the Hessian matrix. Leveraging this novel fundamental result, we design an algorithm with a sketch size proportional to the effective dimension and which exhibits a quadratic rate of convergence. This result dramatically improves on the classical linear-quadratic convergence rates of state-of-the-art sub-sampled Newton methods. However, in most practical cases, the effective dimension is not known beforehand, and this raises the question of how to pick a sketch size as small as the effective dimension while preserving a quadratic convergence rate. Our second and main contribution is thus to propose an adaptive sketch size algorithm with quadratic convergence rate and which does not require prior knowledge or estimation of the effective dimension: at each iteration, it starts with a small sketch size, and increases it until quadratic progress is achieved. Importantly, we show that the embedding dimension remains proportional to the effective dimension throughout the entire path and that our method achieves state-of-the-art computational complexity for solving convex optimization programs with a strongly convex component. We discuss and illustrate applications to linear and quadratic programming, as well as logistic regression and other generalized linear models.
\end{abstract}

\section{Introduction}
\label{sectionintroduction}

We consider a composite optimization problem of the form
\begin{align}
\label{eqnmainoptimizationprobleminitial}
    x^* \defn \argmin_{x \in \real^d} \left\{f(x) \defn f_0(x) + g(x)\right\}\,,
\end{align}
where $f_0, g : \real^d \to \overline{\real}$ are both closed, twice differentiable convex functions. Here, we denote $\overline{\real} \defn \real \cup \{+\infty\}$ and by $\mathbf{dom}\,f$ the domain of $f$. We are interested in the structured setting where forming the Hessian matrix $\nabla^2 f_0(x)$ is prohibitively expensive, but we have available at small computational cost a Hessian matrix square-root $\nabla^2 f_0(x)^{1/2}$, that is, a matrix $\nabla^2 f_0(x)^{1/2}$ of dimensions $n \times d$ such that $(\nabla^2 f_0(x)^{1/2})^\top \nabla^2 f_0(x)^{1/2} = \nabla^2 f_0(x)$ for some integer $n \gre d$, and $n$ eventually very large. Moreover, we assume the function $g$ to be $\mu$-strongly convex, i.e., $\nabla^2 g(x) \succeq \mu I_d$.

Large-scale optimization problems of this form are becoming ever more common in applications, due to the increasing dimensionality of data (e.g., genomics, medicine, high-dimensional models). Typically, the function $f_0$ may represent an objective value we aim to minimize over a convex set $\mathcal{C} \subseteq \real^d$, that is, we aim to solve $\min_{x \in \mathcal{C}} f_0(x)$. A common practice to turn this constrained optimization problem into an unconstrained one is to add to the objective function a \emph{penalty} or {barrier} function $g(x)$ which encodes $\mathcal{C}$ (e.g., logarithmic barrier functions for polyhedral constraints or $\ell_p$-norm regularization for $\ell_p$-ball constraints). In many cases of practical interest, a matrix square-root $\nabla^2 f_0(x)^{1/2}$ can be computed efficiently. For instance, in the broad context of empirical risk minimization, the function $f_0$ has the separable form $f_0(x) = \sum_{i=1}^m \ell_i(a_i^\top x)$ where the functions $\ell_i$ are twice-differentiable and convex. In this case, a suitable Hessian matrix square root is given by the $n\times d$ matrix $\nabla^2 f_0(x)^{1/2} \defn \mathbf{diag}(\ell_i^{\prime \prime}(a_i^\top x)^{1/2}) \, A$. On the other hand, we assume that the Hessian of the function $g$ is well-structured, so that its computation is relatively cheap in comparison to that of $f_0$. For instance, if the constraint set is the unit simplex (i.e., $x \gre 0$ and $\mathbf{1}^\top x \less 1$), then the Hessian of the associated logarithmic barrier function is a diagonal matrix plus a rank one matrix. Other examples include problems for which $g$ has a separable structure such as typical regularizers for ill-posed inverse problems (e.g., graph regularization $g(x) = \frac{1}{2} \sum_{i,j \in E} (x_i - x_j)^2$, $\ell_p$-norms with $p > 1$ or  approximations of $\ell_1$-norm).

Second-order methods such as the Newton's method enjoy superior convergence in both theory and practice compared to first-order methods, that is, quadratic convergence rate versus $1/T^2$ for accelerated gradient descent. A common issue in first-order methods is the tuning of step size~\cite{asi2019importance}, whose optimal choice depends on the strong convexity and smoothness of the underlying problem. In contrast, whenever the objective function $f$ is \emph{self-concordant}, then Newton's method has the appealing property of being invariant to rescaling and coordinate transformations, is independent of problem-dependent parameters, and thus needs little or no tuning of algorithmic hyperparameters. More precisely, we recall that, given a current iterate $x$, the standard Newton's method computes the Hessian matrix $H(x)$ and the descent direction $\vne$ defined as 
\begin{align}
    & H(x) \defn \nabla^2 f_0(x) + \nabla^2 g(x)\,,\\
    & \vne \defn -H(x)^{-1} \nabla f(x)\,.
\end{align}
Given a step size $s > 0$, it then uses the update 
\begin{align}
    \xne \defn x + s \, \vne\,. 
\end{align}
Despite these advantages, Newton's method requires, at each iteration, forming and solving the high-dimensional linear system $H(x) \vne = -\gx$, which has complexity scaling as $\mathcal{O}(nd^2)$, and this becomes prohibitive in large-scale settings. To address this numerical challenge, a multitude of different approximations to Newton's method have been proposed in the literature. Quasi-Newton methods (e.g, DFP, BFGS and their limited memory versions~\cite{nocedal2006numerical}) are computationally cheaper, but their convergence guarantees require stronger assumptions and are typically much weaker than those of Newtons's method. On the other hand, random projections are an effective way of performing dimensionality reduction~\cite{vempala2005random, mahoney2011randomized, drineas2016randnla}, and many random projection (or \emph{sketching}) based algorithms were designed to reduce the cost of solving the linear Newton system. For instance, the respective methods in~\cite{gower2019rsn} and~\cite{lacotte2019high} embed the optimization variable into a lower dimensional subspace, so that solving the Newton system becomes cheaper;~\cite{qu2016sdna} propose to solve an approximate Newton system based on random principal sub-matrices of a global upper bound on the Hessian; \cite{doikov2018randomized} address a common setting, that of block-separable convex optimization problems, and propose a method combining the ideas of randomized coordinate descent with cubic regularization~\cite{nesterov2012efficiency, nesterov2006cubic}.

Our work builds specifically on a generic method, that is, the Newton sketch~\cite{pilanci2017newton}, which is based on a structured random embedding of the Hessian matrix $H(x)$. Formally, given a sketch size $m$ such that $m \ll n$ and an embedding matrix $S \in \real^{m \times n}$ to be precised, the Newton sketch computes the approximate Hessian $H_S(x)$ and the approximate descent direction $\vnsk$ defined as
\begin{align}
    &H_S(x) \defn (\nabla^2 f_0(x)^\frac{1}{2})^\top S^\top S \nabla^2 f_0(x)^\frac{1}{2} + \nabla^2 g(x)\,,\\
    &\vnsk \defn -H_S(x)^{-1} \gx\,.
\end{align}
Given a step size $s > 0$, it then uses the update 
\begin{align}
    \xnsk \defn x + s \, \vnsk\,. 
\end{align}
For classical embeddings (e.g., sub-Gaussian, randomized orthogonal systems), it has been shown by~\cite{pilanci2017newton} that, in general, a sketch size $m \asymp d$ is sufficient for the Newton sketch to achieve a linear-quadratic convergence rate with high probability (w.h.p.).

\textbf{Contributions.} Our first key contribution is to show that, under the assumption that $g$ is $\mu$-strongly convex, the scaling $m \asymp \overline{d}_\mu \log(\overline{d}_\mu)/\delta$ is sufficient for the Newton sketch to achieve a $\delta$-accurate solution at a \emph{quadratic} convergence rate with high probability. More generally, we show that convergence is geometric provided that $m$ scales appropriately in terms of $\overline{d}_\mu$. Here, the critical quantity $\overline{d}_\mu$ is the effective (Hessian) dimension, defined as
\begin{align}
\label{eqneffectivedimension}
    \overline{d}_\mu \defn \sup_{x \in \mathcal{S}(x_0)} d_\mu(x)\,,
\end{align}
where $x_0$ is the initial point of our algorithm, $\mathcal{S}(x_0)$ is the sublevel set of $f$ at $x_0$, and 
\begin{align}
    d_\mu(x) \defn \mbox{trace}(\nabla^2 f_0(x) (\nabla^2 f_0(x) + \mu I_d)^{-1})
\end{align}
is the \emph{local} effective dimension. Importantly, it always holds that $d_\mu(x) \less \overline{d}_\mu \less \min\{n,d\} = d$. In many applications, the effective dimension is substantially smaller than the ambient dimension $d$~\cite{bach2013sharp, alaoui2015fast, yang2017randomized}. However, in order to pick $m$ in terms of $\overline{d}_\mu$ which is usually unknown and then achieve computational and memory space savings, it is necessary to estimate $\overline{d}_\mu$. There exist randomized techniques for precise estimation of $d_\mu(x)$, but they provably work under stringent assumptions, e.g., $d_\mu(x)$ very small (e.g., see Theorem 60 in~\cite{avron2017sharper}). In the context of ridge regression,~\cite{effective2020lacotte} proposed a sketching-based method with adaptive (time-varying) sketch size scaling as the effective dimension, and without prior knowledge or estimation of it. Starting with a small sketch size, it checks at each iteration whether enough progress is achieved by the update. If not, it doubles the sketch size. The time and memory complexities of this method to return a \emph{certified} $\delta$-accurate solution w.h.p.~scale in terms of the effective dimension, i.e., it takes time $\mathcal{O}(nd \log^2(\overline{d}_\mu)\log(d/\delta))$ with a sketch size $m \lesssim \overline{d}_\mu \log(\overline{d}_\mu)$ for large values of $n$. This significantly improves on usual standard randomized pre-conditioning methods~\cite{rokhlin2008fast, avron2010blendenpik, meng2014lsrn} which require $m \gtrsim d$.

In a vein similar to this adaptive ridge regression solver, our second key contribution is to propose an adaptive sketch size version of the effective dimension Newton sketch. Importantly, we prove that the adaptive sketch size scales in terms of $\overline{d}_\mu$. Furthermore, our adaptive method offers the possibility to the user to choose the convergence rate, from linear to quadratic.

\textbf{Other related works.} Recent studies in the literature on randomized second-order and Sub-sampled Newton methods \cite{byrd2011use,bollapragada2019exact,roosta2019sub,berahas2020investigation} show that picking an embedding dimension proportional to $d$ and possibly smaller than $d$ under certain conditions do work empirically in many settings~\cite{xu2016sub, xu2020second, wang2018giant}. The recent work by~\cite{li2020subsampled} provides a more precise understanding of these phenomena. In the context of empirical risk minimization with $\ell^2_2$-regularization, they show that the subsampled Newton method with $m \asymp \overline{d}_\mu$ data points is enough to guarantee convergence. However, differently from our work, their method needs to estimate the effective dimension at each iteration. Furthermore, their convergence guarantees severely depend on the condition number of the problem (e.g., see their Theorems~1 and~2), whereas our results are independent of condition numbers and only involve the relevant dimensions of the problem ($n,d,\overline{d}_\mu$) and the target accuracy. Besides effective dimension based sampling, sketching-based methods are used in the context of distributed optimization where due to stringent memory and/or communication constraints, reducing the number of iterations and/or the size of second-order information is critical~\cite{shamir2014communication, derezinski2020debiasing, bartan2020distributed}.

\subsection{Notations and background}
A closed convex function $\varphi : \real^d \to \overline{\real}$ is \emph{self-concordant} if $|\varphi^{\prime\prime\prime}(x)| \less 2 \, (\varphi^{\prime\prime}(x))^{3/2}$. This definition extends to a closed convex function $f : \real^d \to \overline{\real}$ by imposing this requirement on the univariate functions $\varphi_{x,y}(t) \defn f(x+ty)$ for all choices of $x,y$ in the domain of $f$. Self-concordance is a typical assumption for the analysis of the classical Newton's method, in order to obtain convergence results which are independent of unknown problem parameters (e.g., strong convexity, smoothness or Lipschitz constants; see the books by~\cite{nesterov2003introductory} or~\cite{boyd2004convex} for further background), and this encompasses many widely used functions in practice, e.g., linear, quadratic, negative logarithm. Hence, in this work, we assume that \emph{$f_0$ and $g$ are self-concordant functions.}

The choice of the sketching matrix $S \in \real^{m \times n}$ is critical for statistical and computational performances. The well-structured subsampled randomized Hadamard transform (SRHT)~\cite{ailon2006approximate} usually serves as a reference for comparing sketching algorithms thanks to its strong subspace embedding properties~\cite{mahoney2011randomized, drineas2016randnla, dobriban2019asymptotics, lacotte2020limiting} and fast sketching time $\mathcal{O}(nd \log m)$ compared to the classical sketching cost $\mathcal{O}(ndm)$ of sub-Gaussian embeddings. Another typical choice is the sparse Johnson-Lindenstrauss transform (SJLT)~\cite{nelson2013osnap, woodruff2014sketching} with, for instance, one non-zero entry per column. With $A \in \real^{n \times d}$, a sketch $SA$ is then much faster to compute (it takes time $\mathcal{O}(\mathbf{nnz}(A))$) at the expense of weaker subspace embedding properties.

\subsection{Organization of the paper}

In Section~\ref{sectionpreliminaries}, we introduce critical quantities and preliminary results for both the implementation of our algorithms and their analysis. We show that the approximate Newton direction $\vnsk$ is close to the exact one $\vne$, provided that the sketch size scales in terms of $\overline{d}_\mu$. In Section~\ref{sectioneffdimnewtonsketch}, we formally introduce our (non-adaptive) effective dimension Newton sketch algorithm (see Algorithm~\ref{algeffdimnewtonsketch}), and we present several relevant applications. Assuming knowledge of $\overline{d}_\mu$, we prove that its convergence rate is geometric. In Section~\ref{sectionadaptiveeffdimnewtonsketch}, we introduce an adaptive version of Algorithm~\ref{algeffdimnewtonsketch} (see Algorithm~\ref{algorithmadaptive}): importantly, it does not require knowledge of $\overline{d}_\mu$, but still guarantees geometric convergence as well as low memory complexity in terms of $\overline{d}_\mu$. We summarize our complexity guarantees in Table~\ref{tablecomparisoncomplexity} and compare to standard first- and second-order methods and to the original Newton sketch algorithm~\cite{pilanci2017newton} whose implementation and guarantees are agnostic to the effective dimension of the problem. Finally, we show in Section~\ref{sectionexperiments} the empirical benefits of our adaptive method, compared to several standard optimization baselines.

\section{Preliminaries}
\label{sectionpreliminaries}

Critical to our algorithms and their analysis are the Newton and approximate Newton decrements, defined as
\begin{align}
    &\newtondec \defn \left(\gx^\top H(x)^{-1} \gx \right)^\frac{1}{2}\,,\\
    &\appnewtondec \defn \left(\gx^\top H_S(x)^{-1} \gx \right)^\frac{1}{2}\,.
\end{align}
Importantly, for a self-concordant function $f$, the optimality gap at any point $x \in \mathbf{dom}\,f$ is bounded in terms of the Newton decrement as
\begin{align}
    f(x) - f(x^*) \less \lambda_f(x)^2\,.
\end{align}
Due to the expensive cost of computing the Newton decrement $\lambda_f(x)$ as opposed to $\wtilde \lambda_f(x)$, we will aim to characterize, w.h.p.~over the randomness of the sketching matrix, similar optimality bounds and properties with $\wtilde \lambda_f(x)$.

Given $x \in \mathbf{dom}\,f$, a sketch size $m \gre 1$, a random embedding $S \in \real^{m \times d}$ and a sampling precision parameter $\varepsilon > 0$, we consider the following probability event which is critical to our convergence guarantees,
\begin{align}
\label{eqneventcloseness}
    \boxed{\mathcal{E}_{x,m, \varepsilon} \defn \left\{ (1-\frac{\varepsilon}{2}) I_d \preceq C_S \preceq (1+\frac{\varepsilon}{2}) I_d \right\}}\,,
\end{align}
where $C_S \defn H^{-\frac{1}{2}} H_S H^{-\frac{1}{2}}$, $H \equiv H(x)$ and $H_S \equiv H_S(x)$. In words, when $\mathcal{E}_{x,m, \varepsilon}$ holds true, the matrix $H^{-1/2} H_S H^{-1/2}$ is a close approximation of the identity, i.e., $H^{-1/2} H_S H^{-1/2} \approx H^{-1/2} H H^{-1/2} = I_d$. The next result bounds the probability for this event to hold for different choices of the sketching matrix.
\blems 
\label{lemmaconcentration}
Let $\varepsilon \in (0,1/4)$ and $p \in (0,1/2)$. It holds that $\mathbb{P}(\mathcal{E}_{x,\varepsilon,m}) \gre 1-p$, provided that $m = \Omega(d_\mu(x)^2 / (\varepsilon^2 p))$ for the SJLT with single nonzero element in each column, and, $m = \Omega((d_\mu(x) + \log(1/\varepsilon p) \log(d_\mu(x)/p))/\varepsilon^2)$ for the SRHT.
\elems 
We show next that conditional on $\mathcal{E}_{x,m, \varepsilon}$, the approximate Newton decrement $\wtilde \lambda_f(x)$ is close to $\lambda_f(x)$, as well as the approximate Newton direction $\vnsk$ to the exact one $\vne$.
\btheos[Closeness of Newton decrements]
\label{theoremclosenessnewtondecrements}
Let $\varepsilon \in (0,1/4)$. Conditional on the event $\mathcal{E}_{x,m,\varepsilon}$, it holds that
\begin{align}
    &\|\vne - \vnsk\|_{H(x)} \less \varepsilon \, \|\vne\|_{H(x)}\,,\\
    &\sqrt{1-\varepsilon} \, \lambda_f(x) \less \wtilde \lambda_f(x) \less \sqrt{1+\varepsilon} \, \lambda_f(x)\,.
\end{align}
\etheos 

Given $\varepsilon \in (0,1/4)$, we introduce positive parameters $a,b$ such that $1-\frac{1}{2}\left(\frac{1+\varepsilon}{1-\varepsilon}\right)^2 \gre a$, which we use for backtracking line-search (see Algorithm~\ref{algeffdimnewtonsketch} for details). Furthermore, we define the parameters 
\begin{align*}
    &\eta \defn \frac{1}{8} \, \big(1-\frac{1}{2} \big(\frac{1+\varepsilon}{1-\varepsilon}\big)^2 - a\big)/\big(\frac{1+\varepsilon}{1-\varepsilon}\big)^3\,,\\
    &\nu \defn ab \, \frac{\eta^2}{1+\frac{1+\varepsilon}{1-\varepsilon}\,\eta}\,.
\end{align*}
The next results aim to describe the empirical behavior of our methods. As for the classical Newton's method, we distinguish two phases. The algorithm follows a first phase with constant additive decrease in objective value. In a second phase, it converges faster, i.e., the Newton decrement converges to zero at a geometric rate up to quadratic for an appropriate choice of the hyperparameters.
\blems[First phase decrement]
\label{lemmadecreasefirstphase}
Let $\varepsilon\in(0,1/4)$. Suppose that $\mathcal{E}_{x,m,\varepsilon}$ holds true and that $\wtilde \lambda_f(x) > \eta$. Then, we have that 
\begin{align}
    f(\xnsk) - f(x) \less -\nu\,.
\end{align}
\elems 
We introduce the following numerical function which will prove to be useful to characterize the rate of convergence of our algorithms,
\begin{align}
    \alpha(\tau) \defn 0.57 + \frac{16^\tau}{15}\,.
\end{align}
It is easy to verify that $\alpha(\tau)^{1/\tau} \less 2$ for $\tau \in (0,1]$ and $\alpha(0) \less \frac{16}{25}$.
\blems[Second phase decrement]
\label{lemmadecreasesecondphase}
Let $x \in \mathbf{dom}\,f$, $\tau \in [0,1]$ and $\varepsilon \in (0,1/4)$. Set $\varepsilon^\prime = \varepsilon \, \min\{1, \lambda_f(x)^\tau\}$. We assume that the event $\mathcal{E}_{x, m, \varepsilon^\prime}$ holds and that $\wtilde \lambda_f(x) \less \eta$. Then, we have 
\begin{align}
    \lambda_f(\xnsk) \less \alpha(\tau) \, \lambda_f(x)^{1+\tau}\,.
\end{align} 
Consequently, the progress is geometric for any $\tau \in (0,1]$, i.e., 
\begin{align}
    \alpha(\tau)^{1/\tau} \, \lambda_f(\xnsk) \less \left(\alpha(\tau)^{1/\tau} \, \lambda_f(x)\right)^{1+\tau}\,. 
\end{align}
On the other hand, the progress is linear for $\tau = 0$, i.e., 
\begin{align}
    \lambda_f(\xnsk) \less \frac{16}{25} \, \lambda_f(x)\,.
\end{align}
\elems 
We conclude this section with a simple technical lemma which characterizes a sufficient number of iterations before termination, under geometric convergence.
\blems[Geometric convergence and sufficient iteration number] 
\label{lemmageometricdecrease}
Let $\delta \in (0,1)$, $\alpha > 0$, $\tau \in (0,1]$, and $\{\beta_t\}_{t \gre 0}$ be a sequence of positive numbers such that $\beta_0 \less \eta$, $\eta \alpha^{1/\tau} < 1$, $\sqrt{\delta} \alpha^{1/\tau} < 1$ and $\alpha^{1/\tau} \beta_{t+1} \less (\alpha^{1/\tau} \beta_t)^{1+\tau}$ for all $t \gre 0$. Then, it holds that $\beta_t \less \sqrt{\delta}$ for any $t \gre T_{\tau,\alpha,\delta}$ where
\begin{align}
    T_{\tau, \alpha, \delta} \defn \ceil{\frac{1}{\log(1+\tau)} \,  \log\left( \frac{1 + \frac{\tau \log(1/\delta)}{2 \log(1/\alpha)}}{1 + \frac{\tau \log(1/\eta)}{\log(1/\alpha)}} \right)}\,.
\end{align}
\elems 
Throughout this work, we will use the shorthand
\begin{align}
    \boxed{T_{\tau,\delta} \equiv T_{\tau,\alpha(\tau),\delta}}\,.
\end{align}
Note in particular that $T_{\tau, \delta} = \mathcal{O}(\log(\tau \log(1/\delta)))$ for small $\delta$. Further, it holds that $\lim_{\tau \to 0} T_{\tau, \delta} \less \ceil{\frac{\log(1/\delta)}{\log(25/16)}}$, which corresponds to the classical complexity of linear convergence with rate $16/25$.

\section{Effective dimension Newton sketch}
\label{sectioneffdimnewtonsketch}

We formally introduce our effective dimension Newton sketch method in Algorithm~\ref{algeffdimnewtonsketch}.
\begin{algorithm}[!h]
\caption{Effective dimension Newton sketch}
\label{algeffdimnewtonsketch}
\begin{algorithmic}[1]
\REQUIRE Initial point $x_0 \in \mathbf{dom} f$, threshold sketch sizes $\overline m_1$ and $\overline m_2$, initial sketch size $m_0 = \overline m_1$, line-search parameters $(a,b)$, target accuracy $\delta > 0$, convergence rate parameter $\tau \in [0,1]$ and sampling precision parameter $\varepsilon=1/8$.
\FOR{$t=0,\dots$}
\STATE Sample an $m_t \times n$ embedding $S_t$ independent of $\{S_{j}\}_{j=0}^{t-1}$. Compute $\vnsk$ and $\wtilde \lambda_f(x_t)$ based on $S_t$.
\STATE \textbf{if} $\wtilde \lambda_f(x_t)^2 \less \frac{3}{4}\delta$ \textbf{then return} $x_t$.
\STATE Starting at $s = 1$: \textbf{while} $f(x_t + s \vnsk) > f(x_t) + a s \nabla f(x_t)^\top \vnsk, \quad s \leftarrow b s$. 
\STATE Update $x_{t+1} \leftarrow x_t + s \, \vnsk$.
\STATE If $\wtilde \lambda_f(x_t) > \eta$, set $m_{t+1} = \overline{m}_1$. Otherwise, set $m_{t+1} = \overline{m}_2$.
\ENDFOR
\end{algorithmic}
\end{algorithm}
Algorithm~\ref{algeffdimnewtonsketch} takes as inputs the phase 1 and phase 2 sketch sizes $\overline m_1$ and $\overline m_2$. As we will see in Theorem~\ref{theoremquadraticconvergenceeffdimnewtonsketch}, sufficient values for $\overline m_1$ and $\overline m_2$ to guarantee convergence both depend on the effective dimension $\overline{d}_\mu$. Here and only for Algorithm~\ref{algeffdimnewtonsketch}, we make the idealized assumption that the quantity $\overline{d}_\mu$ is known. In contrast, we introduce in Section~\ref{sectionadaptiveeffdimnewtonsketch} an adaptive method that does not require knowledge of $\overline{d}_\mu$.
\btheos[Geometric convergence guarantees of the Newton sketch]
\label{theoremquadraticconvergenceeffdimnewtonsketch}
Let $\tau \in [0,1]$, $\delta \in (0,1/2)$ and $p_0 \in (0,1/2)$. Set $\varepsilon = 1/8$. Then, the total number of iterations $T_f$ and the total time complexity $\mathcal{C}$ for obtaining a $\delta$-approximate solution $\wtilde x$ in function value (i.e., $f(\wtilde x) -f(x^*) \less \delta$) via Algorithm~\ref{algeffdimnewtonsketch} satisfy
\begin{align}
    &T_f \less \overline{T} \defn \frac{f(x_0) - f(x^*)}{\nu} + T_{\tau,\frac{3}{8}\delta} + 1\,,\\
    &\mathcal{C} = \mathcal{O}\!\left(\overline{m}_2^2 d + nd \log \overline{m}_2\right) \overline{T}\,,
\end{align}
with probability at least $1-p_0$, provided that $\overline{m}_1 \gtrsim \overline{d}_\mu + \log(\frac{\overline{T}}{p_0})\log(\frac{\overline{d}_\mu \overline{T}}{p_0})$ and $\overline{m}_2 \gtrsim \delta^{-\tau} \left(\overline{d}_\mu + \log(\frac{\overline{T}}{p_0 \delta^{\tau/2}})\log(\frac{\overline{d}_\mu \overline{T}}{p_0})\right)$ for the SRHT, whereas for the SJLT, it is sufficient to have $\overline m_1 \gtrsim \frac{\overline{d}_\mu^2 \overline T}{p_0}$ and $\overline{m}_2 \gtrsim \frac{\overline{d}_\mu^2 \overline T}{\delta^\tau p_0}$.
\etheos 
We draw some immediate consequences of Theorem~\ref{theoremquadraticconvergenceeffdimnewtonsketch}, which will be useful for further discussions and comparisons of our complexity guarantees in Section~\ref{sectioncomplexity}. With the SRHT, consider the quadratic convergence case, i.e., $\tau = 1$. We pick a failure probability $p_0 \asymp \frac{1}{\overline{d}_\mu}$, and sketch sizes $\overline{m}_1 \asymp \overline{d}_\mu$ and $\overline{m}_2 \asymp \frac{\overline{d}_\mu \log(\overline{d}_\mu/\delta) }{\delta}$. We observe quadratic convergence with $T_f = \mathcal{O}(\log\log(\frac{1}{\delta}\big))$ iterations. Further, assuming that the sample size $n$ is large enough for the sketching cost $\mathcal{O}(nd \log \overline{m})$ to dominate the cost $\mathcal{O}(\overline{m}^2 d)$ of solving the randomized Newton system, i.e., $n \gtrsim \frac{\overline{d}_\mu^2 \log(\overline{d}_\mu / \delta)}{\delta^2}$, then the total complexity results in 
\begin{align}
    \mathcal{C} = \mathcal{O}\big( nd \log\big(\frac{\overline{d}_\mu}{\delta}\big) \log \log(\frac{1}{\delta}) \big).
\end{align}
Similarly, we consider the linear convergence case, i.e., $\tau=0$. For simplicity, suppose that $\overline{d}_\mu \gtrsim \log \log(1/\delta)$. We pick $p_0 \asymp \frac{1}{\overline{d}_\mu}$, and sketch sizes $\overline{m}_1 \asymp \overline{m}_2 \asymp \overline{d}_\mu$. We observe linear convergence with $T_f = \mathcal{O}(\log \frac{1}{\delta})$ iterations. Assuming again that the sample size $n$ is large enough for the sketching cost to dominate the cost of solving the randomized Newton system, i.e., $n \gtrsim \overline{d}_\mu^2/\log(\overline{d}_\mu)$, we obtain the total time complexity
\begin{align}
    \mathcal{C} = \mathcal{O}\!\left(nd \log(\overline{d}_\mu) \log(\frac{1}{\delta})\right)\,.
\end{align}
We proceed with a similar discussion for the SJLT at the end of the proof of Theorem~\ref{theoremquadraticconvergenceeffdimnewtonsketch} deferred to the Appendix.

\subsection{Some applications of the effective dimension Newton sketch}
\label{sectionexamples}

We discuss various concrete instantiations of the optimization problem~\eqref{eqnmainoptimizationprobleminitial} where the function $g$ satisfies $\mu$-strong convexity and for which forming the partially sketched Hessian $H_S(x)$ is amenable to fast computation.

\bexs[Ridge regression]
We consider the optimization problem
\begin{align}
    \min_{x\in \real^d} \left\{f(x) \defn \frac{1}{2} \|Ax-b\|_2^2 + \frac{\mu}{2}\|x\|_2^2 \right\}
\end{align}
where $A \in \real^{n \times d}$ with $n \gre d$ and whose solution is given in closed-form by $x^* = (A^\top A + \mu I_d)^{-1} A^\top b$. Direct methods yield the exact solution in time $\mathcal{O}(nd^2)$, whereas first-order methods (e.g., conjugate gradient method) yield an $\delta$-approximate solution in time $\mathcal{O}(\sqrt{\kappa} nd \log(1/\delta))$, where $\kappa$ is the condition number of $A$. Randomized pre-conditioning and sketching methods can improve on this complexity (see Section~\ref{sectioncomplexity} for further details). Here, our setting for the Newton sketch applies with $f_0(x) = \frac{1}{2}\| Ax-b\|_2^2$ (whose square-root Hessian is $A$) and $g(x) = \frac{\mu}{2} \|x\|_2^2$ which is $\mu$-strongly convex.
\eexs 

\bexs[Portfolio optimization]
The optimization problem takes the form 
\begin{align}
    \min_{x \gre 0,\,\sum_{i=1}^d x_j \less 1} \left\{ f_0(x) \defn -r^\top x + \alpha \,\lra{x, \Sigma x} \right\}\,, 
\end{align}
where $\Sigma = A^\top A$ is an empirical covariance matrix based on the data $A \in \real^{n \times d}$, with $n \gre d$. Using the barrier method, we need to solve its penalized version $\min_{x \in \real^d} \{f_0(x) + g(x)\}$, where  $g(x) \defn -\mu \, \sum_{i=1}^d \log(x_i) - \mu \, \log (1-\lra{\mathbf{1},x}) $. We clearly have the Hessian square-root $\nabla^2 f_0(x)^{1/2} = \sqrt{\alpha} A$. Further, $g$ is $\mu$-strongly convex over its domain: indeed, note that for $0 < x_i < 1$, the Hessian of $g$ is $\mu \, \mathbf{diag}(x_i^2)^{-1} + \mu \mathbf{1} \mathbf{1}^\top$ and the first term satisfies $\mu\,\mathbf{diag}(x_i^2)^{-1} \succeq \mu I_d$.
\eexs

\bexs[Solving Lasso via its dual]
Given $A \in \real^{n \times d}$ with $d \gg n$, the dual Lasso problem takes the form 
\begin{align}
    \max_{\|A^T x\|_\infty \less \lambda} \left\{-\frac{1}{2}\|y-x\|_2^2\right\}\,.
\end{align}
Applying the logarithmic barrier method, one needs to solve a sequence of problems of the form $\min_{x \in \real^n} \{f_0(x) + g(x)\}$ where $g(x) \defn \frac{\mu}{2} \|y-x\|_2^2$, $f_0(x) \defn - \sum_{j=1}^d \log(\lambda - \lra{a_j, x}) - \sum_{j=1}^d \log(\lambda + \lra{a_j,x})$ and $a_j$ is the $j$-th column of $A$. This form is amenable to the Newton sketch: a square-root of $\nabla^2 f_0(x)$ is given by $\nabla^2 f_0(x)^{1/2} = \mathbf{diag}(|\lambda - \lra{a_j, x}|^{-1} + |\lambda + \lra{a_j, x}|^{-1}) \, A^T$, and the function $g(x)$ is $\mu$-strongly convex.
\eexs

\bexs[Regularized logistic regression with $n \gg d$]\label{exm:logistic}
We consider data points $\{(a_i, y_i)\}_{i=1}^n$ where each $a_i$ is a $d$-dimensional feature vector with binary response $y_i \in \{\pm1\}$. We aim to find a linear classifier through regularized logistic regression, that is, 
\begin{align}
\min_{x \in \real^d} \left\{ \sum_{i=1}^n \log\!\left(1 + e^{-y_i a_i^\top x}\right) + \frac{\mu}{2} \|x\|_2^2 \right\}\,.
\end{align}
Setting $f_0(x) = \sum_{i=1}^n \log\!\left(1 + e^{-y_i a_i^\top x}\right)$ and $g(x) = \frac{\mu}{2} \|x\|_2^2$, we have that $\nabla^2 f_0(x)^\frac{1}{2} = \mathbf{diag}(h) A$ where the $i$-th coefficient of $h \in \real^n$ is given by $h_i = \frac{e^{y_i a_i^\top x/2}}{1+e^{y_i a_i^\top x}}$. More generally, empirical risk minimization with generalized linear models yields a Hessian square-root of the form 'diagonal times data matrix $A$'.
\eexs 

\bexs[Projection onto polyhedra]
Given $v \in \real^d$, $A \in \real^{n \times d}$ with $n \gg d$ and $b \in \real^n$ such that there exists $x_0\in \real^d$ that satisfies $A x_0 < b$, we aim to solve the optimization problem 
\begin{align}
    \min_{x \in \real^d}\, \frac{1}{2} \|x - v\|_2^2\,,\quad \mbox{s.t.}\quad Ax \less b\,.
\end{align}
Applying a barrier method, one needs to solve a sequence of optimization problems of the form $\min_x f_0(x) + g(x)$, where $f_0(x) \defn - \sum_{i=1}^n \log(b_i - a_i^\top x)$ and $g(x) = \frac{\mu}{2} \|x-v\|_2^2$. Clearly, $g$ is $\mu$-strongly convex, and a square-root of $\nabla^2 f_0(x)$ is given by $\mathbf{diag}(|b_i - a_i^\top x|^{-1}) A$.
\eexs

\section{Adaptive Newton Sketch with effective dimensionality}
\label{sectionadaptiveeffdimnewtonsketch}

We turn to the adaptive version of Algorithm~\ref{algeffdimnewtonsketch}, which starts with small sketch size and does not require knowledge or estimation of the effective dimension $\overline{d}_\mu$. Importantly, our method is guaranteed to converge at a tunable geometric rate, and with a sketch size scaling in terms of $\overline{d}_\mu$.

For $\tau \in [0,1]$ and $\varepsilon \in (0,1/4)$, we set 
\begin{align}
    \alpha(\tau,\varepsilon) \defn \frac{(1+\varepsilon)^\frac{1}{2}}{(1-\varepsilon)^{\frac{1+\tau}{2}}}\, \alpha(\tau)\,,
\end{align}
and we will consider in this section the sufficient number of iterations $T_{\tau, \alpha(\tau,\varepsilon), \frac{\delta}{d}}$ as defined in Lemma~\ref{lemmageometricdecrease} for $\alpha = \alpha(\tau, \varepsilon)$. Our adaptive method is formally described in Algorithm~\ref{algorithmadaptive}. It starts each iteration by checking whether $\wtilde \lambda_f(x_t) > \eta$. If so, assuming the sketch size $m_t$ large enough, we have w.h.p.~by Lemma~\ref{lemmadecreasefirstphase} that $f(\xnsk) - f(x) \less -\nu$ and we set $x_{t+1} = \xnsk$. Otherwise, if $\wtilde \lambda_f(x_t) \less \eta$, we have w.h.p.~by Lemma~\ref{lemmadecreasesecondphase} a condition similar to $\wtilde \lambda_f(\xnsk) \less \alpha(\tau, \varepsilon) (\wtilde \lambda_f(x_t))^{1+\tau}$, in which case we set $x_{t+1} = \xnsk$. If none of the above events happen, we increase the sketch size by a factor $2$. On the other hand, if the sketch size is not large enough for the guarantees of Lemmas~\ref{lemmadecreasefirstphase} and~\ref{lemmadecreasesecondphase} to hold w.h.p., then either the algorithm terminates with a potentially small sketch size, or, the sketch size must at some point become large enough due to the doubling trick.

\begin{algorithm}[!h]
\caption{Adaptive effective dimension Newton sketch}
\label{algorithmadaptive}
\begin{algorithmic}[1]
\REQUIRE Initial point $x_0 \in \mathbf{dom} f$, initial sketch size $m_0 = \overline m_0$, line-search parameters $(a,b)$, target accuracy $\delta \in (0,1/2)$, convergence rate parameter $\tau \in [0,1]$ and sampling precision parameter $\varepsilon=1/8$.
\FOR{$t=0,\dots$}
\STATE Sample $S_t \in \real^{m_t \times n}$ independent of $S_{t-1}, \dots, S_0$.
\STATE Compute $\vnsk$ and $\wtilde \lambda_f(x_t)$ based on $S_t$.
\STATE \textbf{if} $\wtilde \lambda_f(x_t)^2 \less \frac{\delta}{d}$ then \textbf{return} $x_t$.
\STATE Find step size $s$ with backtracking line search, and set $\xnsk = x_t + s \vnsk$.
\IF{$\wtilde \lambda_f(x_t) > \eta$}
\IF{$f(\xnsk) - f(x) \less -\nu$}
\STATE Set $x_{t+1} = \xnsk$ and $m_{t+1} = m_t$.
\ELSE 
\STATE Set $x_{t+1} = x_t$ and $m_{t+1} = 2 m_t$.
\ENDIF 
\ELSE 
\STATE Sample $S^+ \in \real^{m_t \times n}$ independent of $S_t, \dots, S_0$.
\STATE Compute $v^+ = -H_{S^+}^{-1}\nabla f(\xnsk)$ and $\wtilde \lambda_f(\xnsk) = (-\langle \nabla f(\xnsk), v^+\rangle)^{1/2}$.
\IF{$\wtilde \lambda_f(\xnsk) \less \alpha(\tau, \varepsilon) (\wtilde \lambda_f(x_t))^{1+\tau}$}
\STATE Set $x_{t+1} = \xnsk$, $m_{t+1} = m_t$, $\vnsk = v^+$ and go to step 4.
\ELSE 
\STATE Set $x_{t+1} = x_t$ and $m_{t+1} = 2 m_t$.
\ENDIF 
\ENDIF 
\ENDFOR
\end{algorithmic}
\end{algorithm}
Note that if Algorithm~\ref{algorithmadaptive} terminates, then it returns an iterate $x$ such that $\wtilde \lambda_f(x)^2 \less \frac{\delta}{d}$. We prove next (see Lemma~\ref{lemmaterminationcondition}) that this termination condition implies the $\delta$-approximation guarantee, i.e., $f(x) - f(x^*) \less \delta$ w.h.p., provided that the \emph{initial} sketch size is large enough, and regardless of the final sketch size.
\blems[Termination condition]
\label{lemmaterminationcondition}
Let $\delta \in (0,1/2)$ and $p \in (0,1/2)$, and suppose that Algorithm~\ref{algorithmadaptive} returns $x$. Then, it holds that $f(x) - f(x^*) \less \delta$ with probability at least $1-p$ provided that $m_0 \gtrsim \log^2(1/p)$ for the SRHT, and, $m_0 \gtrsim 1/p$ for the SJLT. 
\elems

\subsection{Time and memory space complexity guarantees}
\label{sectioncomplexity}

For conciseness, we present a succinct version of our complexity guarantees for the adaptive Newton sketch (only for the linear rate $\tau=0$ and for the quadratic rate $\tau=1$). A more general statement for any $\tau \in [0,1]$ can be found in the proof of Theorem~\ref{theoremconvergenceadaptive}.
\btheos[Geometric convergence guarantees of the adaptive Newton sketch]
\label{theoremconvergenceadaptive}
Let $\tau \in [0,1]$, $p_0 \in (0,1/2)$ and $\delta \in (0,1/2)$. Let $\overline{m}_0$ be an initial sketch size. Then, it holds with probability at least $1-p_0$ that Algorithm~\ref{algorithmadaptive} returns a $\delta$-approximate solution $\wtilde x$ in function value (i.e., $f(\wtilde x) - f(x^*) \less \delta$) in less than $\overline T = \mathcal{O}\!\left(T_{\tau, \alpha(\tau,\varepsilon), \frac{\delta}{d}} \log (\overline{d}_\mu)\right)$ iterations, with final sketch size bounded by $2\,\overline m$ and with total time complexity $\overline{\mathcal{C}}$. The values of $\overline{m}_0$, $\overline m$ and $\overline{\mathcal{C}}$ depend on the choice $S$ as follows.

\textbf{(SRHT).} For $\tau = 1$ (quadratic rate), picking $p_0 \asymp \frac{\delta}{d}$ and assuming $n$ large enough such that $n\gtrsim \frac{d^2 \overline{d}_\mu^2}{\delta^2}$, we have $\overline{m}_0 \asymp \frac{d}{\delta}\log(\frac{d}{\delta})$, $\overline m \asymp \frac{d}{\delta}(\overline d_\mu + \log(\frac{d}{\delta})\log(\overline{d}_\mu))$, $\overline T = \mathcal{O}(\log(\overline{d}_\mu) \log \log(d/\delta))$ and
\begin{align}
    \overline{\mathcal{C}} = \mathcal{O}\Big(nd \log(\overline{d}_\mu) \log(d/\delta) \log\log(d/\delta)\Big)\,.
\end{align}
For $\tau = 0$ (linear rate), picking $p_0 \asymp 1/\overline{d}_\mu$ and assuming $n$ large enough such that $n \gtrsim \frac{\overline{d}^2_\mu}{\log(\overline{d}_\mu)}$, we have $\overline{m}_0 \asymp \log^2(d/\delta)$, $\overline m \asymp \overline{d}_\mu$, $\overline T = \mathcal{O}(\log(\overline{d}_\mu)\log(d/\delta))$ and 
\begin{align}
    \overline{\mathcal{C}} = \mathcal{O}\Big(nd \log^2(\overline{d}_\mu) \log(d/\delta)\Big)\,.
\end{align}
\textbf{(SJLT).} For $\tau = 1$ (quadratic rate), assuming $n$ large enough such that $n \gtrsim \frac{\overline{d}_\mu^4 d^{2} \log(\log(d/\delta))^2}{\delta^{2} p_0^2}$, we have $\overline{m}_0 \asymp \frac{d \log(\log(d/\delta))}{p_0 \delta}$, $\overline{m} \asymp \frac{d \overline{d}^2_\mu \log(\log(d/\delta))}{p_0 \delta}$, $\overline{T} = \mathcal{O}(\log(\overline{d}_\mu) \log \log(d/\delta))$ and 
\begin{align}
    \overline{\mathcal{C}} = \mathcal{O}\!\left( nd \, \log(\overline{d}_\mu) \, \log(\log(d/\delta))) \right)\,.
\end{align}
For $\tau = 0$ (linear rate), assuming $n$ large enough such that $n \gtrsim \frac{\overline{d}_\mu^4 \log(d/\delta)^2}{p_0^2}$, we have $\overline{m}_0 \asymp \frac{\log(d/\delta)}{p_0}$, $\overline{m} \asymp \frac{\overline{d}_\mu^2 \log(d/\delta)}{p_0}$, $\overline{T} = \mathcal{O}(\log(\overline{d}_\mu) \log(d/\delta))$ and
\begin{align}
    \overline{\mathcal{C}} = \mathcal{O}\!\left( nd \, \log(\overline{d}_\mu) \, \log(d/\delta)) \right)\,.
\end{align}    
\etheos 
Note that adaptivity with convergence rate parameter $\tau$ comes at the cost of an additional $d^\tau$ factor for the final sketch size, compared to Algorithm~\ref{algeffdimnewtonsketch}. This is essentially due to our exit condition threshold $\delta/d$ that we choose for the following reason. For small $m \gtrsim 1$, the approximate Newton decrement $\wtilde \lambda^2_f (x)$ may fluctuate around $\lambda^2_f(x)$ by a factor up to $\overline{d}_\mu$ (see Theorem~1 in~\cite{cohen2015optimal}). In this case, the exit condition $\wtilde \lambda_f(x_t)^2 \approx \delta$ would result in $f(x_t) - f(x^*) \approx \delta \, \overline{d}_\mu$. To guarantee $\delta$-accuracy, it is sufficient to use the termination condition $\delta/\overline{d}_\mu$ to account for these fluctuations. As $\overline d_\mu$ is unknown, we choose to divide instead by $d$.

We summarize our complexity guarantees in Table~\ref{tablecomparisoncomplexity}. In contrast to gradient descent (GD), Nesterov's accelerated gradient descent (NAG) and Newton's method (NE), our time complexity has no condition number dependency and scales linearly in $nd$ up to log-factors, and so does the original Newton sketch (NS). The NS log-factor is at least $\log(d) \log(1/\delta)$ whereas our SRHT-quadratic mode adaptive method has a log-factor $\log(\overline{d}_\mu) \log(d/\delta) \log(\log(d/\delta))$. The latter is much smaller for effective dimension $\overline{d}_\mu$ small compared to $d$ and $1/\delta$. Furthermore, in terms of memory, our algorithm starts with small $m$ whereas NS uses a constant sketch size $m \gtrsim d$. For $\tau = 0$, our memory savings are drastic when $\overline{d}_\mu$ is small. There are downsides to our method, in comparison to NS. For $\overline{d}_\mu$ close to $d$, our time complexity bounds become worse than NS, by an adaptivity-cost factor $\log \overline{d}_\mu$ for both $\tau=0$ and $\tau=1$. For $\tau=1$, our \emph{worst-case} sketch size is always greater than that of NS, by a factor $\overline{d}_\mu/\delta$, which comes from enforcing quadratic convergence. 

When $\log \overline d_\mu \gg \log \log (d/\delta)$, then our SRHT/quadratic mode adaptive method yields a better time complexity than its linear mode counterpart, but at the expense of worse memory complexity. 

In the context of ridge regression, we note that the time complexity of our SRHT-linear mode adaptive method scales similarly to the complexity of the adaptive method proposed by~\cite{effective2020lacotte} for returning a certified $\delta$-accurate solution. Importantly, our method applies to a much broader range of optimization problems, and can achieve better complexity by tuning the convergence rate parameter $\tau$. 

We emphasize again that our guarantees hold in a worst-case sense. In practice, the sketch size can start from a small value and may remain significantly smaller than the bounds in Table~\ref{tablecomparisoncomplexity}, which we illustrate in our numerical experiments.
\begin{table*}[!h]
	\caption{We compare the time complexity of different optimization methods in order to achieve a $\delta$-accurate solution in function value, for a function with condition number $\kappa$. 'NS-effdim' (resp.~'NS-ada') refers to our Algorithm~\ref{algeffdimnewtonsketch} (resp.~our Algorithm~\ref{algorithmadaptive}); 'linear' (resp.~'quadratic') signifies the choice $\tau=0$ (resp.~$\tau=1$). We refer to~\cite{nesterov2003introductory} for gradient descent (GD), Nesterov's accelerated gradient method (NAG) and Newton's method (NE); we refer to~\cite{pilanci2017newton} for the Newton sketch (NS), and we refer to Theorems~\ref{theoremquadraticconvergenceeffdimnewtonsketch} and~\ref{theoremconvergenceadaptive} for our algorithms. Katyusha \cite{allen2017katyusha} is an instance of an accelerated SVRG method whose condition number dependency improves on NAG. For each algorithm, we assume that the sample size $n$ is large enough for the time complexity to scale at least linearly in the term $nd$.}
    \label{tablecomparisoncomplexity}
	\centering
	\begin{tabular}{|c|c|c|c|c|}
		\toprule
		Algorithm & Time complexity & Sketch size & Proba. & Linear scaling regime \\
		\midrule
		GD & $\kappa \, nd \, \log(1/\delta)$ & - & $1$ & - \\
		\midrule
		NAG & $\sqrt{\kappa} \, nd \, \log(1/\delta)$ & - & $1$ & -\\
		\midrule
		Katyusha & $(nd + d \sqrt{\kappa n}) \log(1/\delta)$ & - & $1$ & -\\
		\midrule
		NE & $nd^2 \log(\log(1/\delta))$ & - & $1$ & - \\
		\midrule 
		NS & $nd \log(d) \log(1/\delta)$ & $d$ & $1-\frac{1}{d}$ & $n \gtrsim \frac{d^2}{\log d}$ \\
		\midrule
		NS-effdim & $nd \log(\overline{d}_\mu) \log(1/\delta)$ & $\overline{d}_\mu$ & $1-\frac{1}{\overline{d}_\mu}$ & $n \gtrsim \frac{d_\mu^2}{\log(d_\mu)}$\\
		(SHRT, linear)& & & & \\
		\midrule
		NS-effdim & $nd \log(\overline{d}_\mu/\delta)\log(\log(1/\delta))$ & $\frac{\overline{d}_\mu}{\delta} \log(\overline{d}_\mu/\delta)$ & $1-\frac{1}{\overline{d}_\mu}$ & $n \gtrsim \frac{\overline{d}_\mu^2 \log(\frac{\overline{d}_\mu}{\delta})}{\delta^2}$\\
		(SHRT, quadratic)& & & &\\
		\midrule
		NS-ada & $nd \log(\overline{d}_\mu)^2 \log(d/\delta)$ & $\overline{d}_\mu$ & $1-\frac{1}{\overline{d}_\mu}$ & $n \gtrsim \frac{\overline{d}_\mu^2}{\log(\overline{d}_\mu)}$\\
		(SRHT, linear)& & & &\\
		\midrule
		NS-ada & $nd \log(\overline{d}_\mu) \log(\frac{d}{\delta}) \log(\log(\frac{d}{\delta}))$ & $\frac{d}{\delta} \left(\overline{d}_\mu + \log(\frac{d}{\delta})\log(\overline{d}_\mu)\right)$ & $1-\frac{1}{\overline{d}_\mu}$ & $n \gtrsim \frac{d^2 \overline{d}_\mu^2}{\delta^2}$\\
		(SRHT, quadratic)& & & &\\
		\bottomrule
	\end{tabular}
\end{table*}

\section{Numerical experiments}
\label{sectionexperiments}
In this section, we compare adaptive Newton Sketch (NS-ada) with other optimization methods in regularized logistic regression problems as in Example \ref{exm:logistic}. The compared methods include Newton Sketch (NS) with fixed sketching dimension, Newton's method (NE), gradient descent method (GD), Nesterov's accelerated gradient descent method (NAG) \citep{nesterov}, the stochastic variance reduced gradient method (SVRG) \cite{johnson2013accelerating} and Katyusha \cite{allen2017katyusha} . For NS-ada and NS, we consider both SJLT sketching matrices and random row sampling (RSS) sketching matrices. All numerical experiments are executed on a Dell PowerEdge R840 workstation. Specifically, we use 4 cores with 192GB ram for all compared methods. 

The datasets used in the numerical experiments are collected from LIBSVM\footnote{https://www.csie.ntu.edu.tw/~cjlin/libsvm/} \citep{libsvm}. The datasets for multi-class classification are manually separated into two categories. For example, in MNIST dataset, we classify even and odd digits. For each dataset, we randomly split it into a training set and a test set with the ratio $1:1$. Several additional numerical results and experimental details are reported in Appendix~\ref{sectionadditionalexperiments}.

\subsection{Regularized logistic regression}
We demonstrate the performance of all compared methods on regularized logistic regression problems. The relative error is calculated by $\frac{f(x)-f_\mathrm{ref}+\epsilon}{1+f_\mathrm{ref}}$. Here $f_\mathrm{ref}$ is the minimal training loss function value among all compared methods and $\epsilon=5\times 10^{-7}$ is a small constant.

\begin{figure}[!ht]
\centering
\begin{minipage}[t]{0.45\textwidth}
\centering
\includegraphics[width=\linewidth]{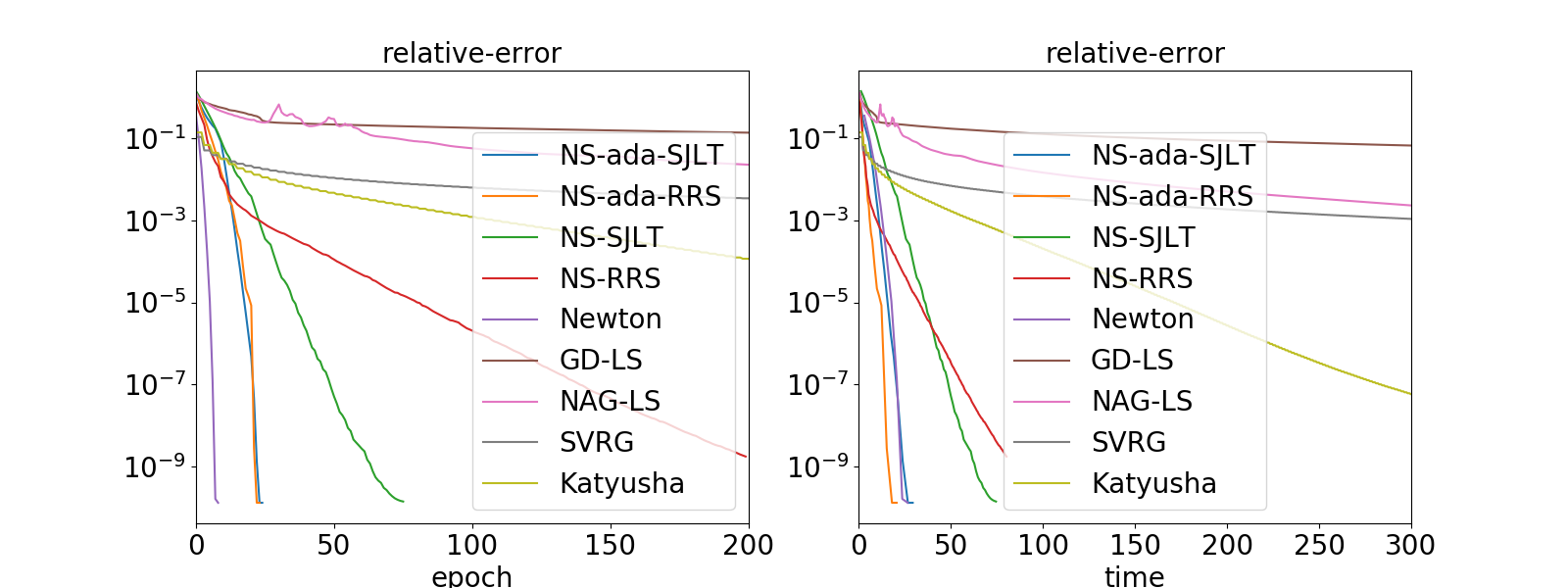}
\end{minipage}
\begin{minipage}[t]{0.45\textwidth}
\centering
\includegraphics[width=\linewidth]{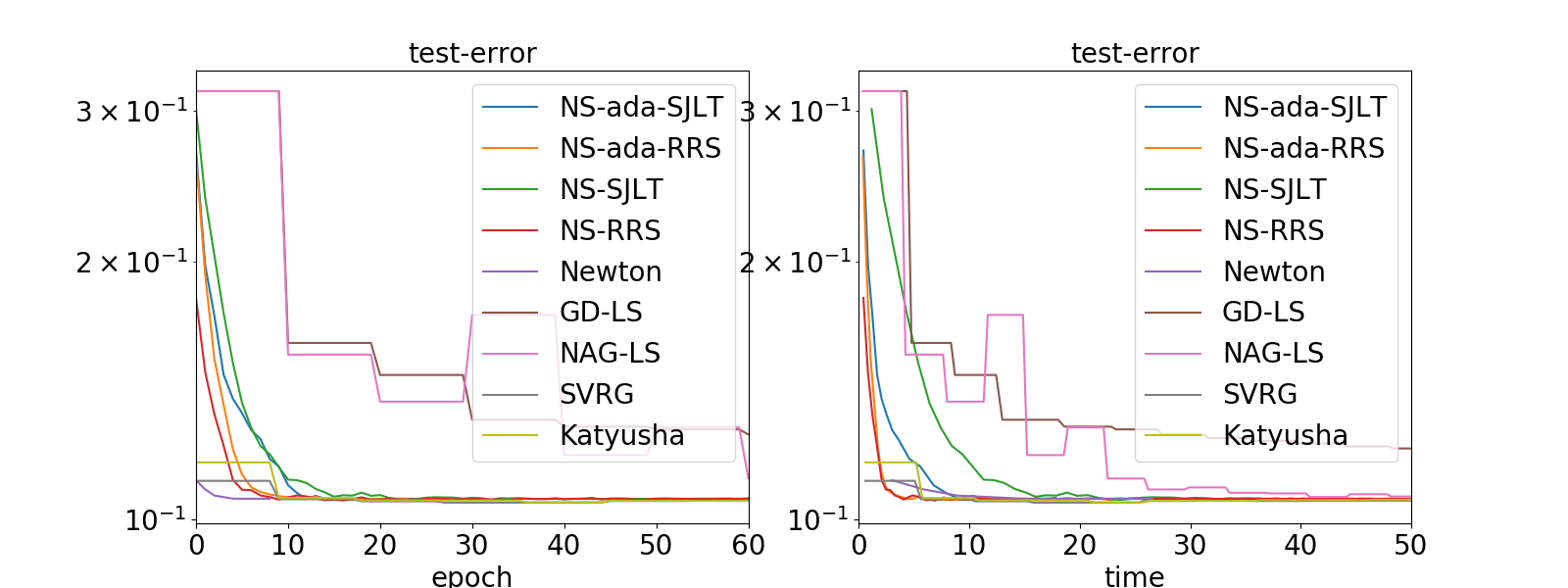}
\end{minipage}
\caption{MNIST. $n=30000, d=780, \mu=10^{-1}$. }
\label{fig:mnist}
\end{figure}

\begin{figure}[!ht]
\centering
\begin{minipage}[t]{0.45\textwidth}
\centering
\includegraphics[width=\linewidth]{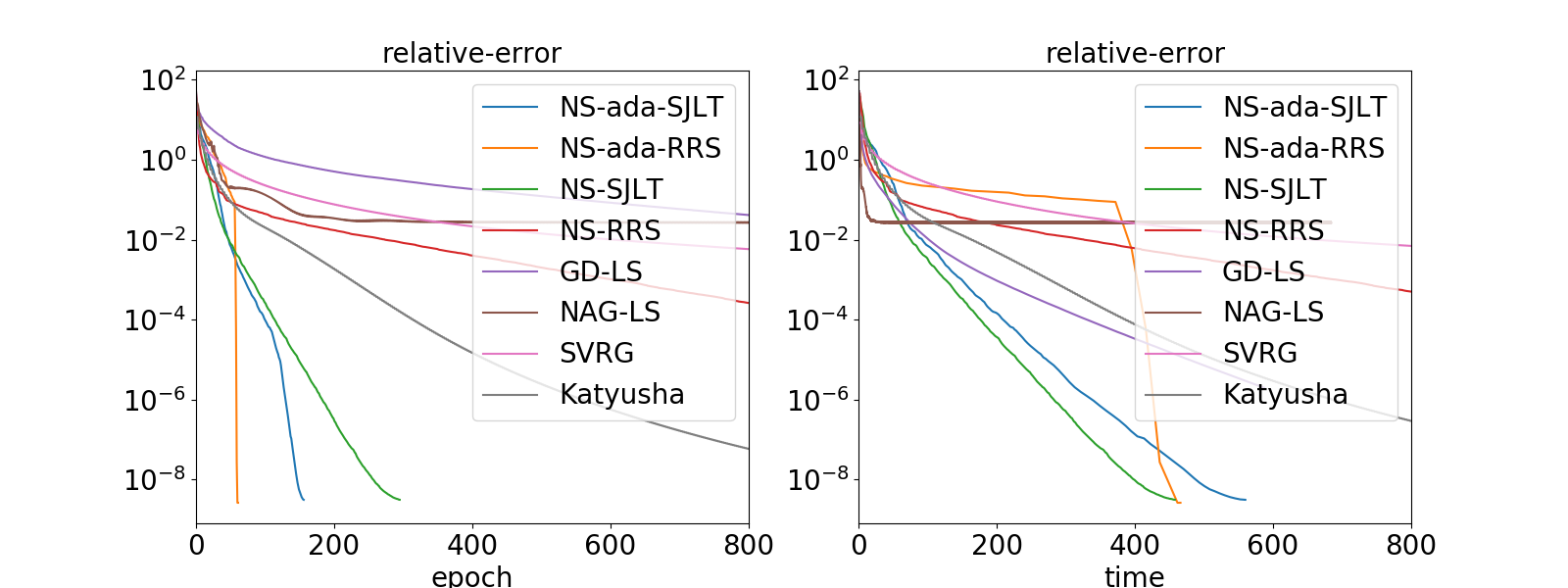}
\end{minipage}
\begin{minipage}[t]{0.45\textwidth}
\centering
\includegraphics[width=\linewidth]{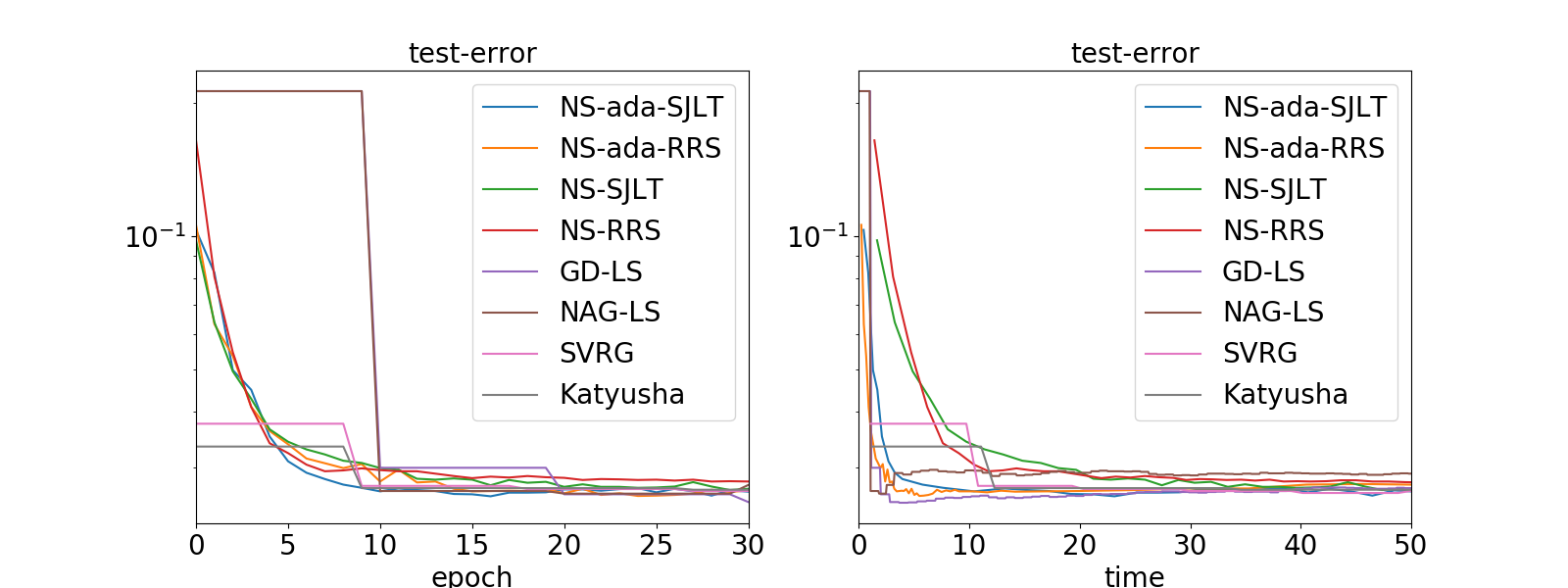}
\end{minipage}
\caption{realsim. $n=50000, d=20958, \mu=10^{-3}$. }
\label{fig:realsim}
\end{figure}

We report the relative error and the test error with respect to the iteration number and the CPU-time in Figures~\ref{fig:mnist} and~\ref{fig:realsim}. NS-ada-SJLT or NS-ada-RRS can achieve the best performance in the relative error with respect to the CPU-time. We can also obverse the super-linear asymptotic convergence rate of NS-ada as it gets closer to the optimum. Compared to NS, NS-ada requires less iterations and less time to converge to a solution with small relative error. Compared to methods utilizing second-order information, first-order methods are less competitive to find a high-precision solution.

\subsection{Regularized logistic regression with kernel matrix}
We also test on regularized logistic regression with the kernel matrix. The relative error and the test error with respect to the iteration number and the CPU-time are plotted in Figure \ref{fig:rcv1} to \ref{fig:epsilon}. NS-ada-SJLT and NS-ada-RRS demonstrate asymptotic super-linear convergence rate of the relative error as the Newton's method in terms of iteration numbers. They also achieve a rapid decrease in relative error in terms of CPU-time. First-order methods have worst performance in terms of relative error. This may come from that the kernel matrices are usually ill-conditioned, i.e., with large condition number.

\begin{figure}[!ht]
\centering
\begin{minipage}[t]{0.45\textwidth}
\centering
\includegraphics[width=\linewidth]{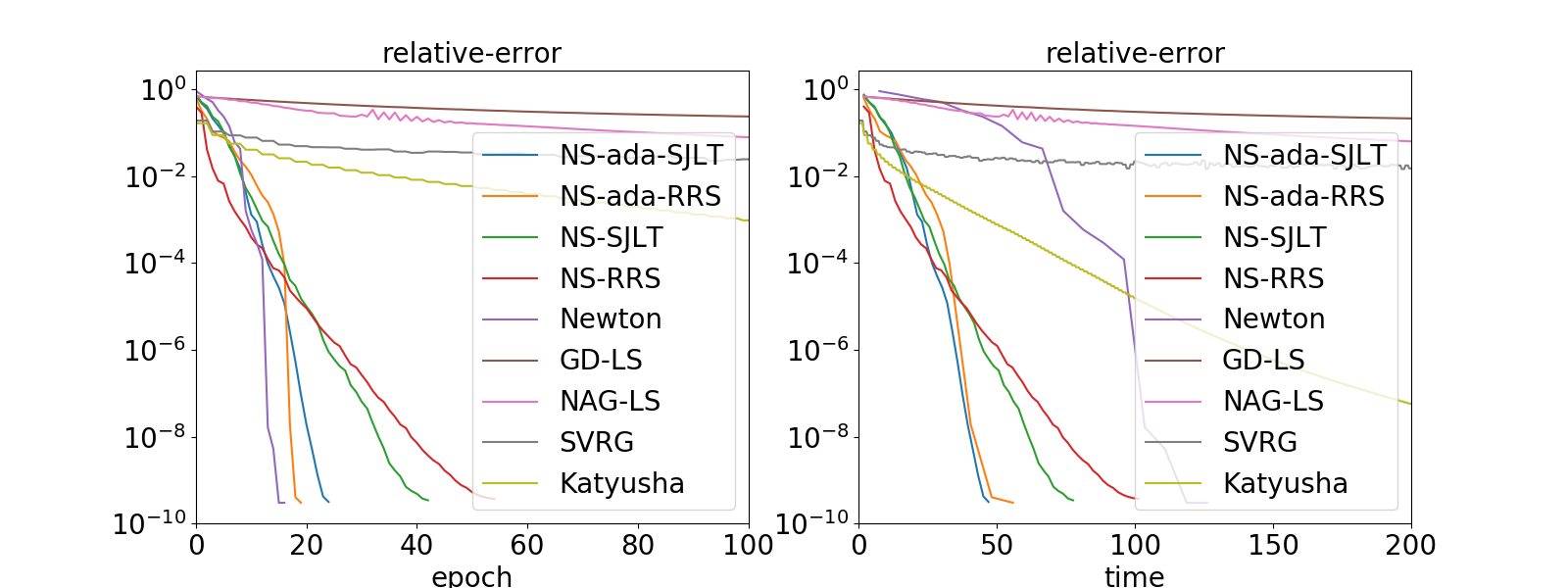}
\end{minipage}
\begin{minipage}[t]{0.45\textwidth}
\centering
\includegraphics[width=\linewidth]{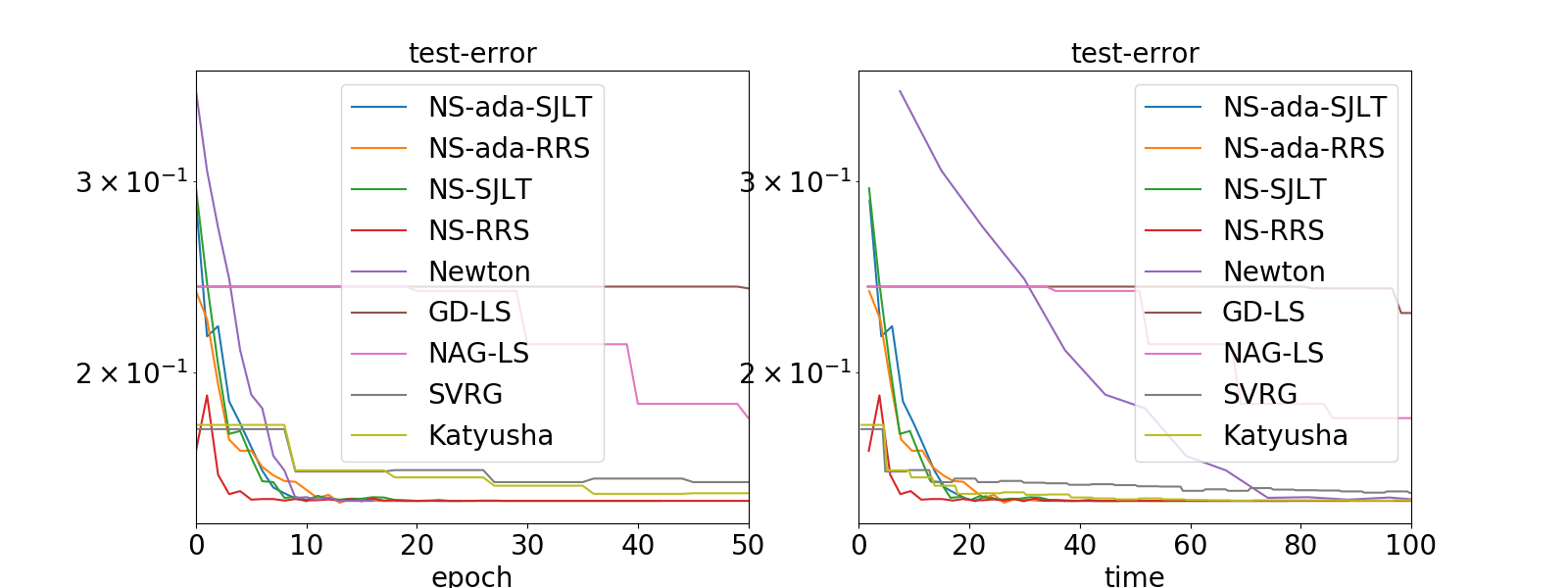}
\end{minipage}
\caption{a8a. kernel matrix. $n=10000,d=10000,\mu=10$.}  
\label{fig:a8a}
\end{figure}

\begin{figure}[!ht]
\centering
\begin{minipage}[t]{0.45\textwidth}
\centering
\includegraphics[width=\linewidth]{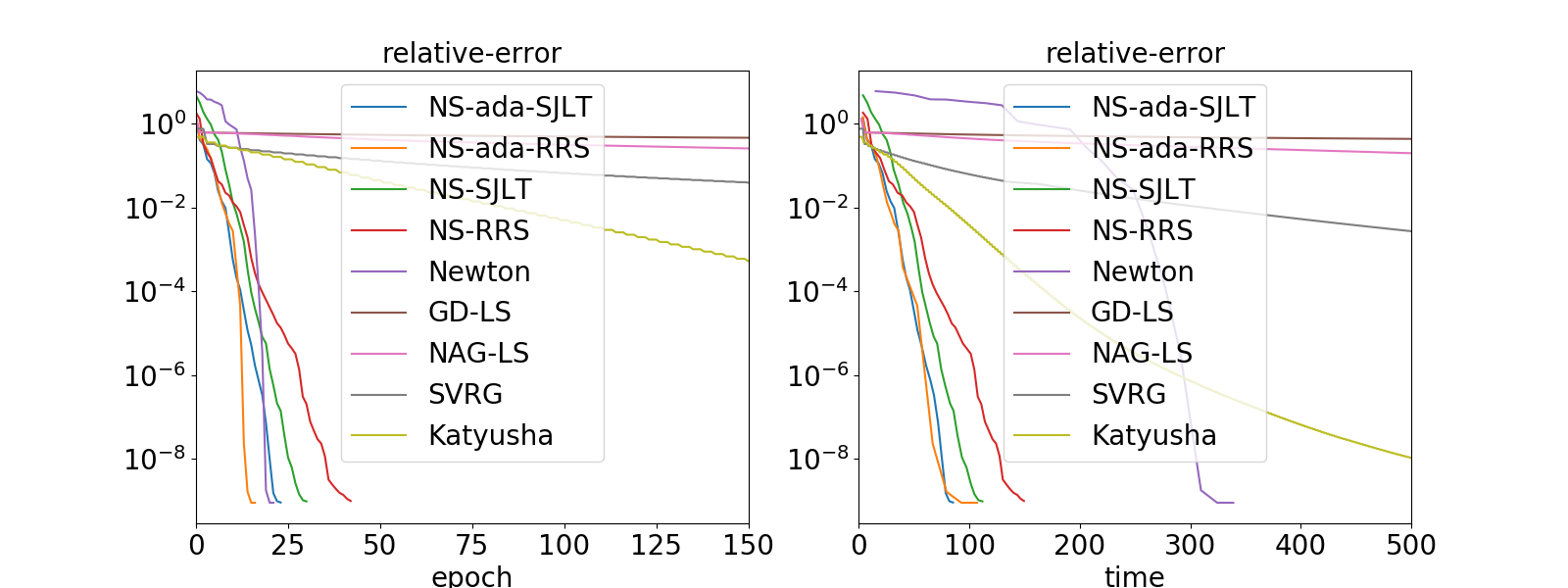}
\end{minipage}
\begin{minipage}[t]{0.45\textwidth}
\centering
\includegraphics[width=\linewidth]{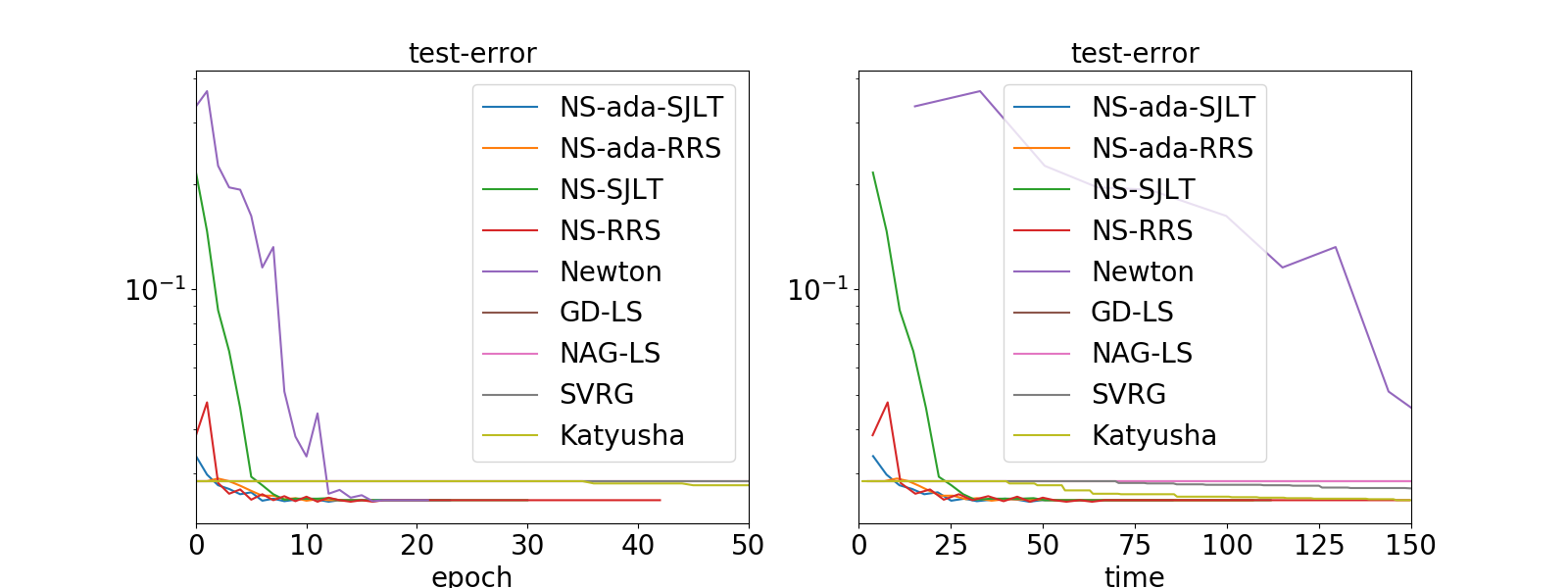}
\end{minipage}
\caption{w7a. kernel matrix. $n=12000,d=12000,\mu=10$. }\label{fig:w7a}
\end{figure}

\section*{Acknowledgements}

We would like to thank the reviewers for their comments. This work was partially supported by the National Science Foundation under grants IIS-1838179 and ECCS-2037304, Facebook Research, Adobe Research and Stanford SystemX Alliance.

\bibliographystyle{icml2021}

\begin{thebibliography}{43}
	\providecommand{\natexlab}[1]{#1}
	\providecommand{\url}[1]{\texttt{#1}}
	\expandafter\ifx\csname urlstyle\endcsname\relax
	\providecommand{\doi}[1]{doi: #1}\else
	\providecommand{\doi}{doi: \begingroup \urlstyle{rm}\Url}\fi
	
	\bibitem[Ailon \& Chazelle(2006)Ailon and Chazelle]{ailon2006approximate}
	Ailon, N. and Chazelle, B.
	\newblock Approximate nearest neighbors and the fast johnson-lindenstrauss
	transform.
	\newblock In \emph{Proceedings of the thirty-eighth annual ACM symposium on
		Theory of computing}, pp.\  557--563. ACM, 2006.
	
	\bibitem[Alaoui \& Mahoney(2015)Alaoui and Mahoney]{alaoui2015fast}
	Alaoui, A.~E. and Mahoney, M.~W.
	\newblock Fast randomized kernel ridge regression with statistical guarantees.
	\newblock In \emph{Proceedings of the 28th International Conference on Neural
		Information Processing Systems-Volume 1}, pp.\  775--783, 2015.
	
	\bibitem[Allen-Zhu(2017)]{allen2017katyusha}
	Allen-Zhu, Z.
	\newblock Katyusha: The first direct acceleration of stochastic gradient
	methods.
	\newblock \emph{The Journal of Machine Learning Research}, 18\penalty0
	(1):\penalty0 8194--8244, 2017.
	
	\bibitem[Asi \& Duchi(2019)Asi and Duchi]{asi2019importance}
	Asi, H. and Duchi, J.~C.
	\newblock The importance of better models in stochastic optimization.
	\newblock \emph{Proceedings of the National Academy of Sciences}, 116\penalty0
	(46):\penalty0 22924--22930, 2019.
	
	\bibitem[Avron et~al.(2010)Avron, Maymounkov, and Toledo]{avron2010blendenpik}
	Avron, H., Maymounkov, P., and Toledo, S.
	\newblock Blendenpik: Supercharging lapack's least-squares solver.
	\newblock \emph{SIAM Journal on Scientific Computing}, 32\penalty0
	(3):\penalty0 1217--1236, 2010.
	
	\bibitem[Avron et~al.(2017)Avron, Clarkson, and Woodruff]{avron2017sharper}
	Avron, H., Clarkson, K.~L., and Woodruff, D.~P.
	\newblock Sharper bounds for regularized data fitting.
	\newblock \emph{Approximation, Randomization, and Combinatorial Optimization.
		Algorithms and Techniques}, 2017.
	
	\bibitem[Bach(2013)]{bach2013sharp}
	Bach, F.
	\newblock Sharp analysis of low-rank kernel matrix approximations.
	\newblock In \emph{Conference on Learning Theory}, pp.\  185--209. PMLR, 2013.
	
	\bibitem[Bartan \& Pilanci(2020)Bartan and Pilanci]{bartan2020distributed}
	Bartan, B. and Pilanci, M.
	\newblock Distributed sketching methods for privacy preserving regression.
	\newblock \emph{arXiv preprint arXiv:2002.06538}, 2020.
	
	\bibitem[Berahas et~al.(2020)Berahas, Bollapragada, and
	Nocedal]{berahas2020investigation}
	Berahas, A.~S., Bollapragada, R., and Nocedal, J.
	\newblock An investigation of newton-sketch and subsampled newton methods.
	\newblock \emph{Optimization Methods and Software}, 35\penalty0 (4):\penalty0
	661--680, 2020.
	
	\bibitem[Bollapragada et~al.(2019)Bollapragada, Byrd, and
	Nocedal]{bollapragada2019exact}
	Bollapragada, R., Byrd, R.~H., and Nocedal, J.
	\newblock Exact and inexact subsampled newton methods for optimization.
	\newblock \emph{IMA Journal of Numerical Analysis}, 39\penalty0 (2):\penalty0
	545--578, 2019.
	
	\bibitem[Boyd \& Vandenberghe(2004)Boyd and Vandenberghe]{boyd2004convex}
	Boyd, S. and Vandenberghe, L.
	\newblock \emph{Convex optimization}.
	\newblock Cambridge university press, 2004.
	
	\bibitem[Byrd et~al.(2011)Byrd, Chin, Neveitt, and Nocedal]{byrd2011use}
	Byrd, R.~H., Chin, G.~M., Neveitt, W., and Nocedal, J.
	\newblock On the use of stochastic hessian information in optimization methods
	for machine learning.
	\newblock \emph{SIAM Journal on Optimization}, 21\penalty0 (3):\penalty0
	977--995, 2011.
	
	\bibitem[Chih-Chung \& Chih-Jen(2011)Chih-Chung and Chih-Jen]{libsvm}
	Chih-Chung, C. and Chih-Jen, L.
	\newblock {LIBSVM}: a library for support vector machines.
	\newblock \emph{ACM Transactions on Intelligent Systems and Technology}, 2011.
	
	\bibitem[Cohen et~al.(2015)Cohen, Nelson, and Woodruff]{cohen2015optimal}
	Cohen, M.~B., Nelson, J., and Woodruff, D.~P.
	\newblock Optimal approximate matrix product in terms of stable rank.
	\newblock \emph{arXiv preprint arXiv:1507.02268}, 2015.
	
	\bibitem[Derezinski et~al.(2020)Derezinski, Bartan, Pilanci, and
	Mahoney]{derezinski2020debiasing}
	Derezinski, M., Bartan, B., Pilanci, M., and Mahoney, M.~W.
	\newblock Debiasing distributed second order optimization with surrogate
	sketching and scaled regularization.
	\newblock In \emph{Conference on Neural Information Processing Systems}, 2020.
	
	\bibitem[Dobriban \& Liu(2019)Dobriban and Liu]{dobriban2019asymptotics}
	Dobriban, E. and Liu, S.
	\newblock Asymptotics for sketching in least squares regression.
	\newblock In \emph{Advances in Neural Information Processing Systems}, pp.\
	3670--3680, 2019.
	
	\bibitem[Doikov \& Richt{\'a}rik(2018)Doikov and
	Richt{\'a}rik]{doikov2018randomized}
	Doikov, N. and Richt{\'a}rik, P.
	\newblock {Randomized block cubic Newton method}.
	\newblock \emph{International Conference on Machine Learning}, 2018.
	
	\bibitem[Drineas \& Mahoney(2016)Drineas and Mahoney]{drineas2016randnla}
	Drineas, P. and Mahoney, M.~W.
	\newblock {RandNLA: randomized numerical linear algebra}.
	\newblock \emph{Communications of the ACM}, 59\penalty0 (6):\penalty0 80--90,
	2016.
	
	\bibitem[Gower et~al.(2019)Gower, Koralev, Lieder, and
	Richt{\'a}rik]{gower2019rsn}
	Gower, R., Koralev, D., Lieder, F., and Richt{\'a}rik, P.
	\newblock {RSN: Randomized subspace newton}.
	\newblock In \emph{Advances in Neural Information Processing Systems}, pp.\
	616--625, 2019.
	
	\bibitem[Johnson \& Zhang(2013)Johnson and Zhang]{johnson2013accelerating}
	Johnson, R. and Zhang, T.
	\newblock Accelerating stochastic gradient descent using predictive variance
	reduction.
	\newblock \emph{Advances in neural information processing systems},
	26:\penalty0 315--323, 2013.
	
	\bibitem[Lacotte \& Pilanci(2020)Lacotte and Pilanci]{effective2020lacotte}
	Lacotte, J. and Pilanci, M.
	\newblock Effective dimension adaptive sketching methods for faster regularized
	least-squares optimization.
	\newblock \emph{Advances in Neural Information Processing Systems}, 33, 2020.
	
	\bibitem[Lacotte et~al.(2019)Lacotte, Pilanci, and Pavone]{lacotte2019high}
	Lacotte, J., Pilanci, M., and Pavone, M.
	\newblock High-dimensional optimization in adaptive random subspaces.
	\newblock \emph{arXiv preprint arXiv:1906.11809}, 2019.
	
	\bibitem[Lacotte et~al.(2020)Lacotte, Liu, Dobriban, and
	Pilanci]{lacotte2020limiting}
	Lacotte, J., Liu, S., Dobriban, E., and Pilanci, M.
	\newblock Limiting spectrum of randomized hadamard transform and optimal
	iterative sketching methods.
	\newblock In \emph{Conference on Neural Information Processing Systems}, 2020.
	
	\bibitem[Li et~al.(2020)Li, Wang, and Zhang]{li2020subsampled}
	Li, X., Wang, S., and Zhang, Z.
	\newblock Do subsampled newton methods work for high-dimensional data?
	\newblock In \emph{Proceedings of the AAAI Conference on Artificial
		Intelligence}, volume~34, pp.\  4723--4730, 2020.
	
	\bibitem[Mahoney(2011)]{mahoney2011randomized}
	Mahoney, M.~W.
	\newblock {Randomized algorithms for matrices and data}.
	\newblock \emph{Foundations and Trends{\textregistered} in Machine Learning},
	3\penalty0 (2):\penalty0 123--224, 2011.
	
	\bibitem[Meng et~al.(2014)Meng, Saunders, and Mahoney]{meng2014lsrn}
	Meng, X., Saunders, M.~A., and Mahoney, M.~W.
	\newblock Lsrn: A parallel iterative solver for strongly over-or
	underdetermined systems.
	\newblock \emph{SIAM Journal on Scientific Computing}, 36\penalty0
	(2):\penalty0 C95--C118, 2014.
	
	\bibitem[Nelson \& Nguy{\^e}n(2013)Nelson and Nguy{\^e}n]{nelson2013osnap}
	Nelson, J. and Nguy{\^e}n, H.~L.
	\newblock Osnap: Faster numerical linear algebra algorithms via sparser
	subspace embeddings.
	\newblock In \emph{2013 ieee 54th annual symposium on foundations of computer
		science}, pp.\  117--126. IEEE, 2013.
	
	\bibitem[Nesterov(1983)]{nesterov}
	Nesterov, Y.
	\newblock A method of solving a convex programming problem with convergence
	rate {$O(1/k^2)$}.
	\newblock \emph{Soviet Mathematics Doklady}, 27\penalty0 (2):\penalty0
	372--376, 1983.
	
	\bibitem[Nesterov(2003)]{nesterov2003introductory}
	Nesterov, Y.
	\newblock \emph{Introductory lectures on convex optimization: A basic course},
	volume~87.
	\newblock Springer Science \& Business Media, 2003.
	
	\bibitem[Nesterov(2012)]{nesterov2012efficiency}
	Nesterov, Y.
	\newblock Efficiency of coordinate descent methods on huge-scale optimization
	problems.
	\newblock \emph{SIAM Journal on Optimization}, 22\penalty0 (2):\penalty0
	341--362, 2012.
	
	\bibitem[Nesterov \& Polyak(2006)Nesterov and Polyak]{nesterov2006cubic}
	Nesterov, Y. and Polyak, B.~T.
	\newblock {Cubic regularization of Newton method and its global performance}.
	\newblock \emph{Mathematical Programming}, 108\penalty0 (1):\penalty0 177--205,
	2006.
	
	\bibitem[Nocedal \& Wright(2006)Nocedal and Wright]{nocedal2006numerical}
	Nocedal, J. and Wright, S.
	\newblock \emph{Numerical optimization}.
	\newblock Springer Science \& Business Media, 2006.
	
	\bibitem[Pilanci \& Wainwright(2017)Pilanci and Wainwright]{pilanci2017newton}
	Pilanci, M. and Wainwright, M.~J.
	\newblock Newton sketch: A near linear-time optimization algorithm with
	linear-quadratic convergence.
	\newblock \emph{SIAM Journal on Optimization}, 27\penalty0 (1):\penalty0
	205--245, 2017.
	
	\bibitem[Qu et~al.(2016)Qu, Richt{\'a}rik, Tak{\'a}c, and Fercoq]{qu2016sdna}
	Qu, Z., Richt{\'a}rik, P., Tak{\'a}c, M., and Fercoq, O.
	\newblock {SDNA: Stochastic dual Newton ascent for empirical risk
		minimization}.
	\newblock In \emph{International Conference on Machine Learning}, pp.\
	1823--1832, 2016.
	
	\bibitem[Rokhlin \& Tygert(2008)Rokhlin and Tygert]{rokhlin2008fast}
	Rokhlin, V. and Tygert, M.
	\newblock A fast randomized algorithm for overdetermined linear least-squares
	regression.
	\newblock \emph{Proceedings of the National Academy of Sciences}, 105\penalty0
	(36):\penalty0 13212--13217, 2008.
	
	\bibitem[Roosta-Khorasani \& Mahoney(2019)Roosta-Khorasani and
	Mahoney]{roosta2019sub}
	Roosta-Khorasani, F. and Mahoney, M.~W.
	\newblock Sub-sampled newton methods.
	\newblock \emph{Mathematical Programming}, 174\penalty0 (1):\penalty0 293--326,
	2019.
	
	\bibitem[Shamir et~al.(2014)Shamir, Srebro, and Zhang]{shamir2014communication}
	Shamir, O., Srebro, N., and Zhang, T.
	\newblock Communication-efficient distributed optimization using an approximate
	newton-type method.
	\newblock In \emph{International conference on machine learning}, pp.\
	1000--1008. PMLR, 2014.
	
	\bibitem[Vempala(2005)]{vempala2005random}
	Vempala, S.~S.
	\newblock \emph{The random projection method}, volume~65.
	\newblock American Mathematical Society, 2005.
	
	\bibitem[Wang et~al.(2018)Wang, Roosta, Xu, and Mahoney]{wang2018giant}
	Wang, S., Roosta, F., Xu, P., and Mahoney, M.~W.
	\newblock Giant: Globally improved approximate newton method for distributed
	optimization.
	\newblock \emph{Advances in Neural Information Processing Systems},
	31:\penalty0 2332--2342, 2018.
	
	\bibitem[Woodruff et~al.(2014)]{woodruff2014sketching}
	Woodruff, D.~P. et~al.
	\newblock Sketching as a tool for numerical linear algebra.
	\newblock \emph{Foundations and Trends{\textregistered} in Theoretical Computer
		Science}, 10\penalty0 (1--2):\penalty0 1--157, 2014.
	
	\bibitem[Xu et~al.(2016)Xu, Yang, Roosta-Khorasani, R{\'e}, and
	Mahoney]{xu2016sub}
	Xu, P., Yang, J., Roosta-Khorasani, F., R{\'e}, C., and Mahoney, M.~W.
	\newblock Sub-sampled newton methods with non-uniform sampling.
	\newblock \emph{arXiv preprint arXiv:1607.00559}, 2016.
	
	\bibitem[Xu et~al.(2020)Xu, Roosta, and Mahoney]{xu2020second}
	Xu, P., Roosta, F., and Mahoney, M.~W.
	\newblock Second-order optimization for non-convex machine learning: An
	empirical study.
	\newblock In \emph{Proceedings of the 2020 SIAM International Conference on
		Data Mining}, pp.\  199--207. SIAM, 2020.
	
	\bibitem[Yang et~al.(2017)Yang, Pilanci, Wainwright,
	et~al.]{yang2017randomized}
	Yang, Y., Pilanci, M., Wainwright, M.~J., et~al.
	\newblock Randomized sketches for kernels: Fast and optimal nonparametric
	regression.
	\newblock \emph{The Annals of Statistics}, 45\penalty0 (3):\penalty0 991--1023,
	2017.
	
\end{thebibliography}

\newpage
\onecolumn 
\appendix

\section{Additional experimental details}
\label{sectionadditionalexperiments}
For NS-ada, we double the sketching dimension when $\tilde \lambda_f(x^{t+1})>c_1\tilde \lambda_f(x^t)\min(1,c_2\tilde \lambda_f(x^t)^{\tau})$. Here $c_1,c_2>0$ and $\tau\in[0,1]$. For all compared methods, we use the backtracking line search method to find a step size satisfying the Armijo condition. For NS-ada, NS and NE, we stop the algorithm when $\tilde \lambda_f(x)<10^{-6}$ or $\lambda_f(x)<10^{-6}$. For first-order methods, we first compute a referenced solution $\tilde x^*$ based on NS-ada. Then, we stop the algorithm when $\frac{f(x)-f(\tilde x^*)}{1+f(\tilde x^*)}<10^{-6}$. 

The parameters for NS-ada and NS for each dataset are summarized in Tables~\ref{tab:NS-ada-SJLT} to \ref{tab:NS}.

\begin{table}[htbp]
 \centering
    \begin{tabular}{|c|c|c|c|c|}
    \hline
         Dataset & $m_0$ & $c_1$ & $\tau$ & $c_2$ \\\hline
         RCV1 & 100&2&0&1\\\hline
         MNIST & 100&0.5&1&6\\ \hline
         gisette & 100&2&0&1\\ \hline
         realsim & 100&2&0&1\\ \hline
         epsilon & 100&1&0&1\\ \hline
    \end{tabular}
    \caption{Parameters for adaptive Newton sketch with SJLT sketching.}
    \label{tab:NS-ada-SJLT}
\end{table}

\begin{table}[htbp]
    \centering
    \begin{tabular}{|c|c|c|c|c|}
    \hline
         Dataset & $m_0$ & $c_1$ & $\tau$ & $c_2$ \\\hline
         RCV1 & 100&1&0&1\\\hline
         MNIST & 100&0.5&1&6\\ \hline
         gisette & 100&2&0&1\\ \hline
         realsim & 100&2&0&1\\ \hline
         epsilon & 100&1&0&1\\ \hline
    \end{tabular}
    \caption{Parameters for adaptive Newton sketch with RRS sketching.}
    \label{tab:NS-ada-RRS}
\end{table}

\begin{table}[htbp]
    \centering
    \begin{tabular}{|c|c|c|}
    \hline
         Dataset & $m$ (SJLT) & $m$ (RRS) \\\hline
         RCV1 & 800&800\\\hline
         MNIST & 800&1600\\ \hline
         gisette & 400&400\\ \hline
         realsim & 800&3200\\ \hline
         epsilon & 800&3200\\ \hline
    \end{tabular}
    \caption{Sketching dimensions of Newton Sketch.}
    \label{tab:NS}
\end{table}

\begin{figure}[htbp]
\centering
\begin{minipage}[t]{0.45\textwidth}
\centering
\includegraphics[width=\linewidth]{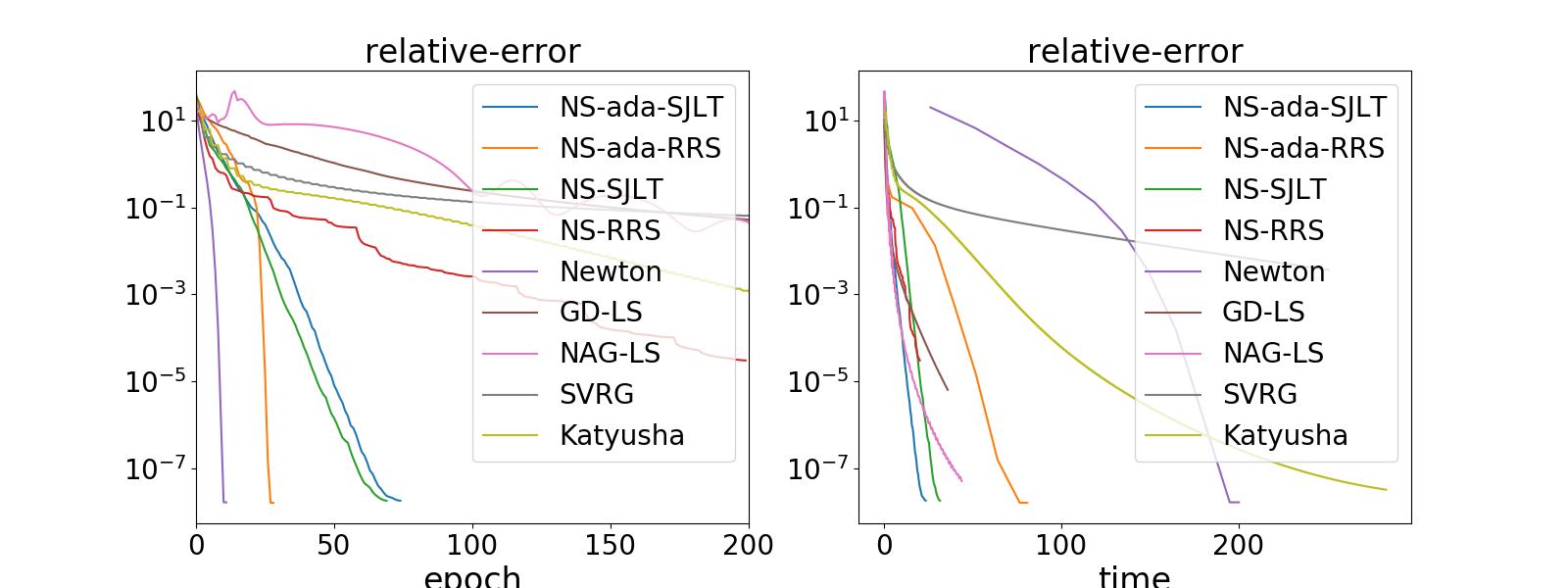}
\end{minipage}
\begin{minipage}[t]{0.45\textwidth}
\centering
\includegraphics[width=\linewidth]{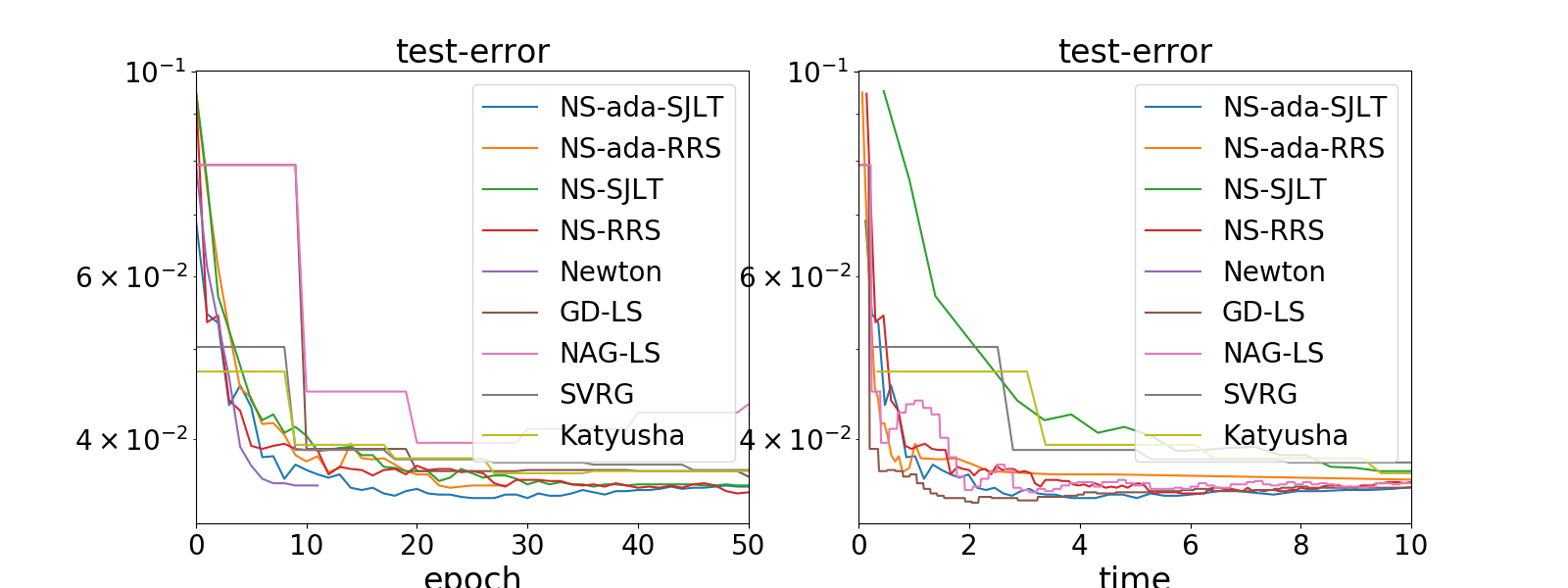}
\end{minipage}
\begin{minipage}[t]{0.45\textwidth}
\centering
\includegraphics[width=\linewidth]{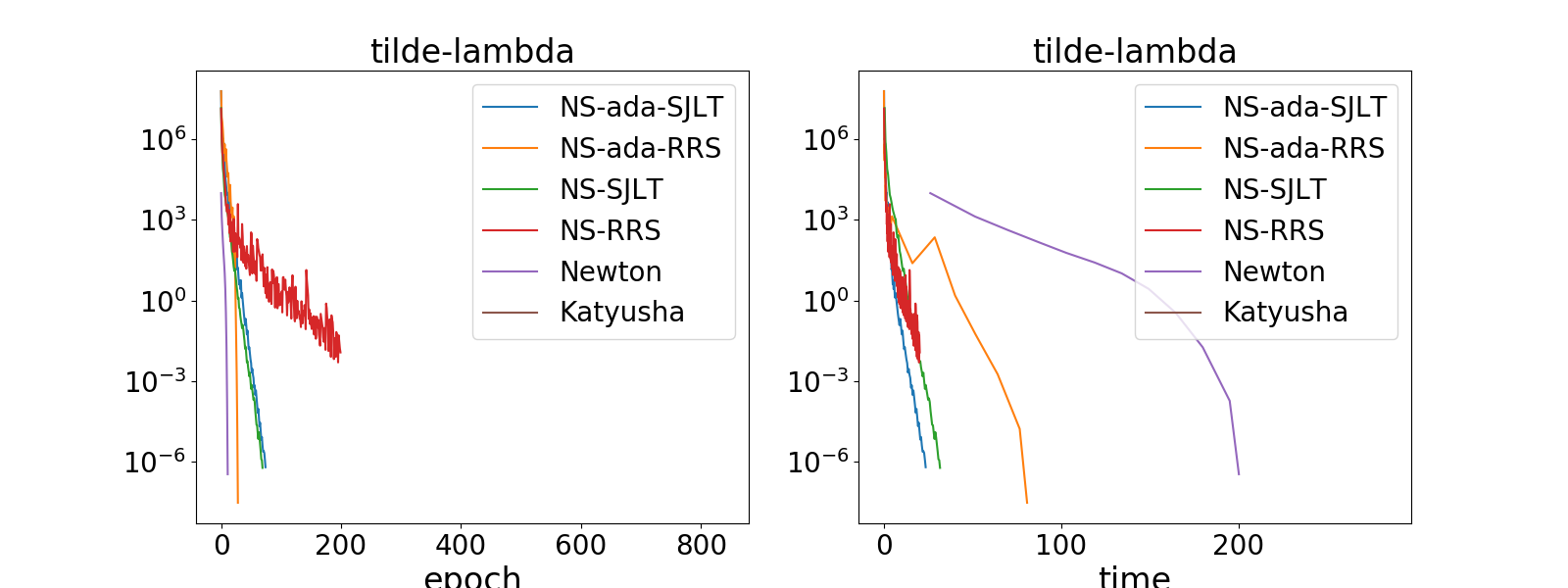}
\end{minipage}
\begin{minipage}[t]{0.45\textwidth}
\centering
\includegraphics[width=\linewidth]{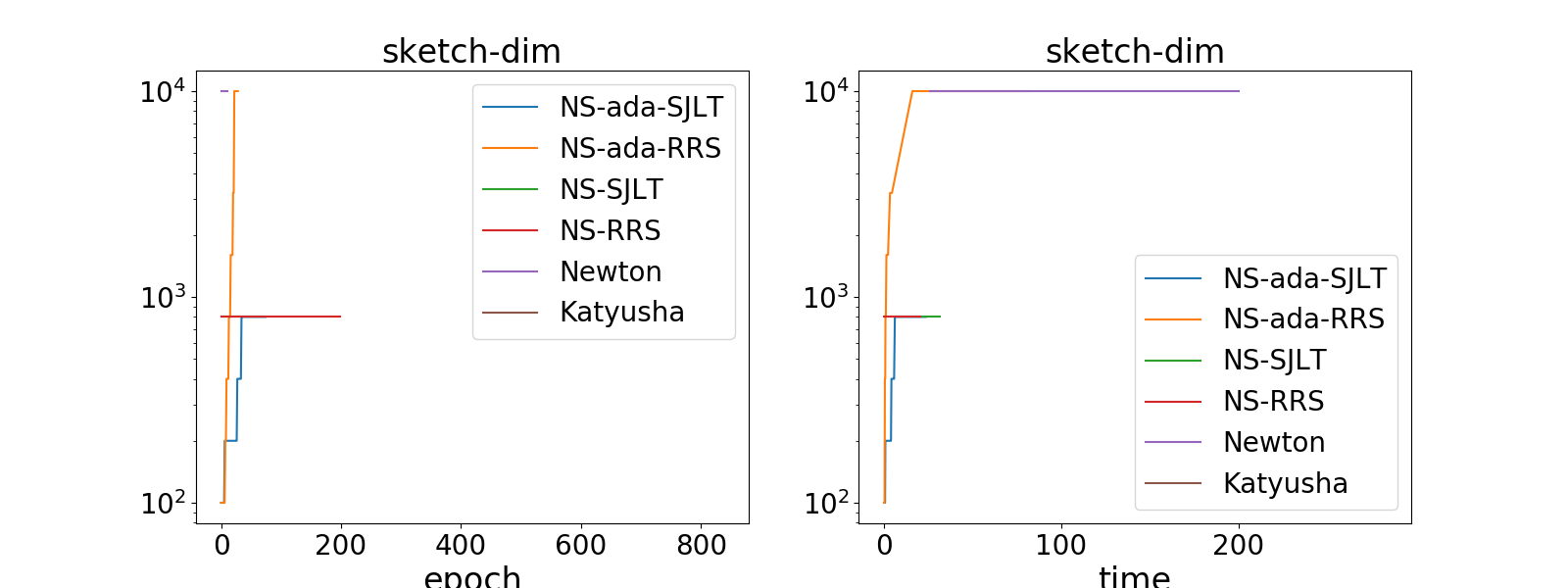}
\end{minipage}
\caption{RCV1. $n=10000, d=47236, \mu=10^{-3}$. }
\label{fig:rcv1}
\end{figure}

\begin{figure}[htbp]
\centering
\begin{minipage}[t]{0.45\textwidth}
\centering
\includegraphics[width=\linewidth]{figures/MNIST-n30000-d780-mu1.0e-01/relative-error.png}
\end{minipage}
\begin{minipage}[t]{0.45\textwidth}
\centering
\includegraphics[width=\linewidth]{figures/MNIST-n30000-d780-mu1.0e-01/test-acc.png}
\end{minipage}
\begin{minipage}[t]{0.45\textwidth}
\centering
\includegraphics[width=\linewidth]{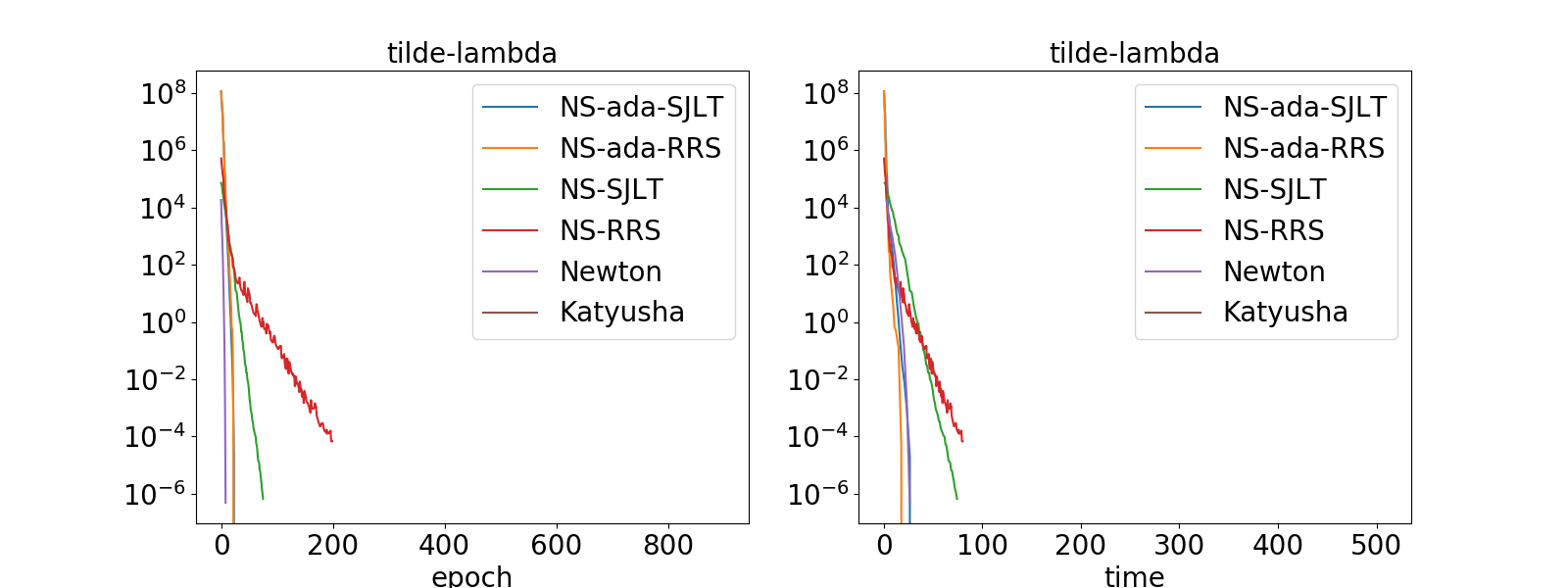}
\end{minipage}
\begin{minipage}[t]{0.45\textwidth}
\centering
\includegraphics[width=\linewidth]{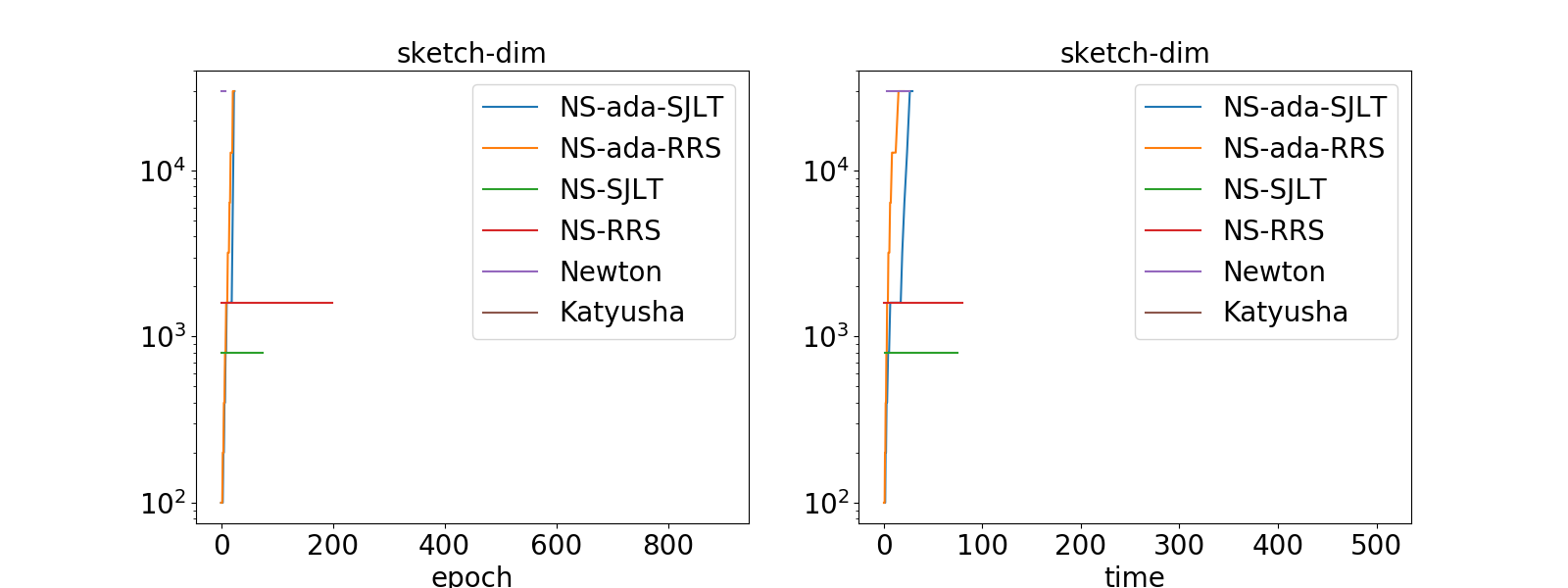}
\end{minipage}
\caption{MNIST. $n=30000, d=780, \mu=10^{-1}$. }
\label{fig:mnist-add}
\end{figure}

\begin{figure}[htbp]
\centering
\begin{minipage}[t]{0.45\textwidth}
\centering
\includegraphics[width=\linewidth]{figures/realsim-n50000-d20958-mu1.0e-03/relative-error.png}
\end{minipage}
\begin{minipage}[t]{0.45\textwidth}
\centering
\includegraphics[width=\linewidth]{figures/realsim-n50000-d20958-mu1.0e-03/test-acc.png}
\end{minipage}
\begin{minipage}[t]{0.45\textwidth}
\centering
\includegraphics[width=\linewidth]{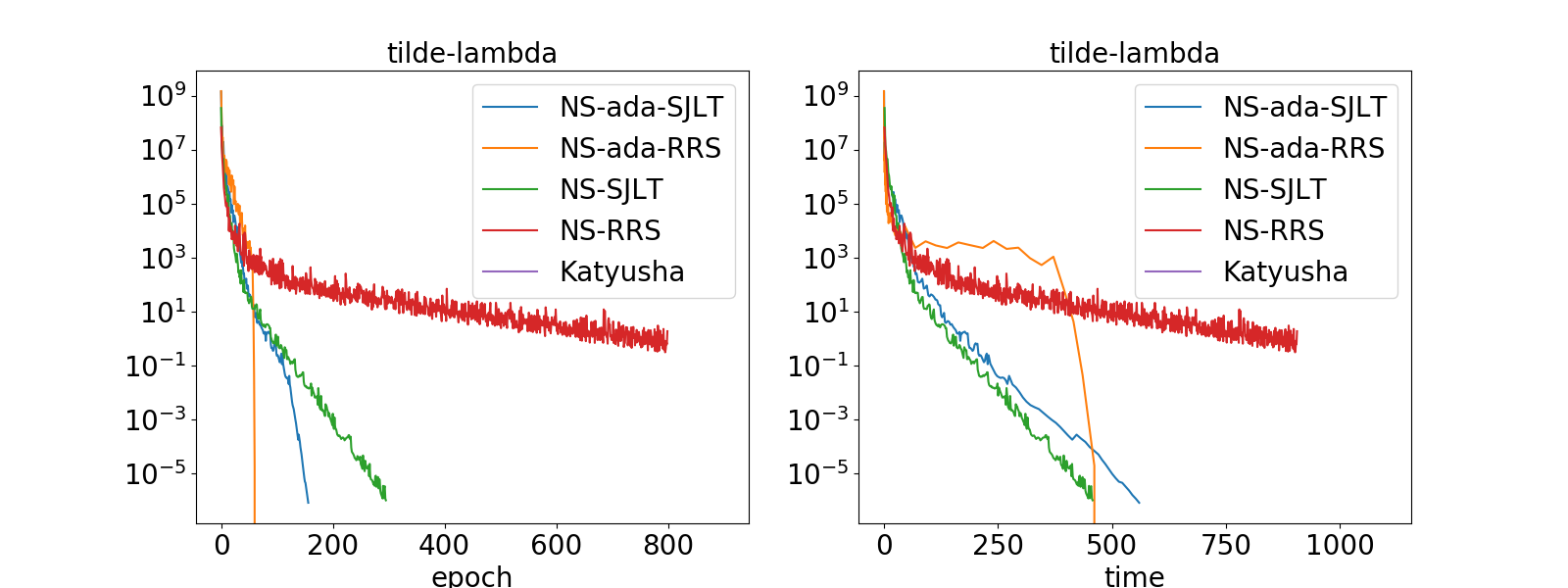}
\end{minipage}
\begin{minipage}[t]{0.45\textwidth}
\centering
\includegraphics[width=\linewidth]{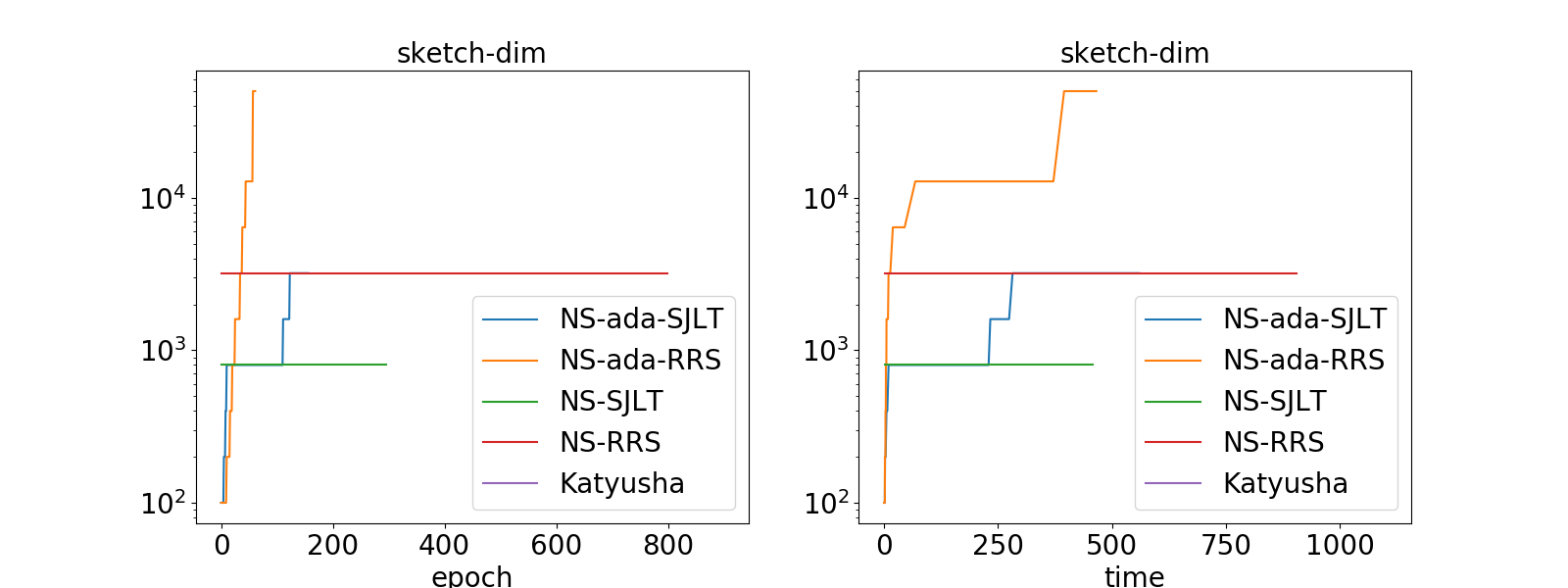}
\end{minipage}
\caption{realsim. $n=50000, d=20958, \mu=10^{-3}$. }
\label{fig:realsim-add}
\end{figure}

\begin{figure}[htbp]
\centering
\begin{minipage}[t]{0.45\textwidth}
\centering
\includegraphics[width=\linewidth]{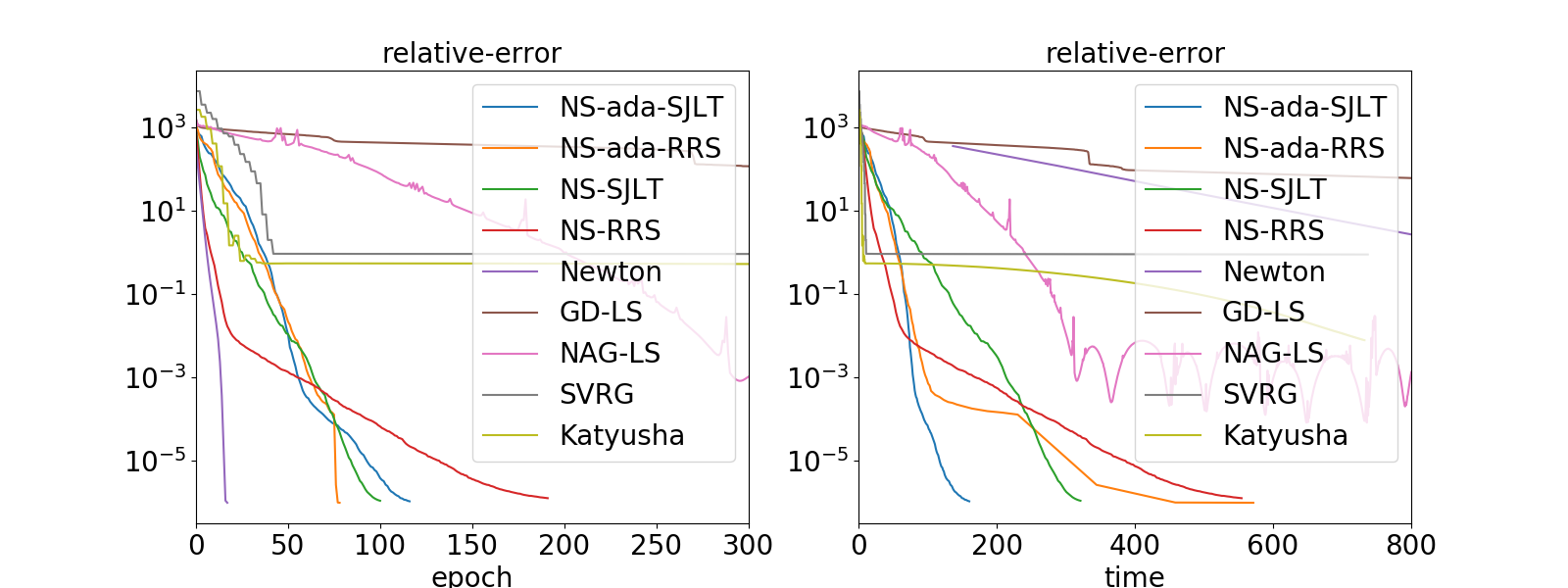}
\end{minipage}
\begin{minipage}[t]{0.45\textwidth}
\centering
\includegraphics[width=\linewidth]{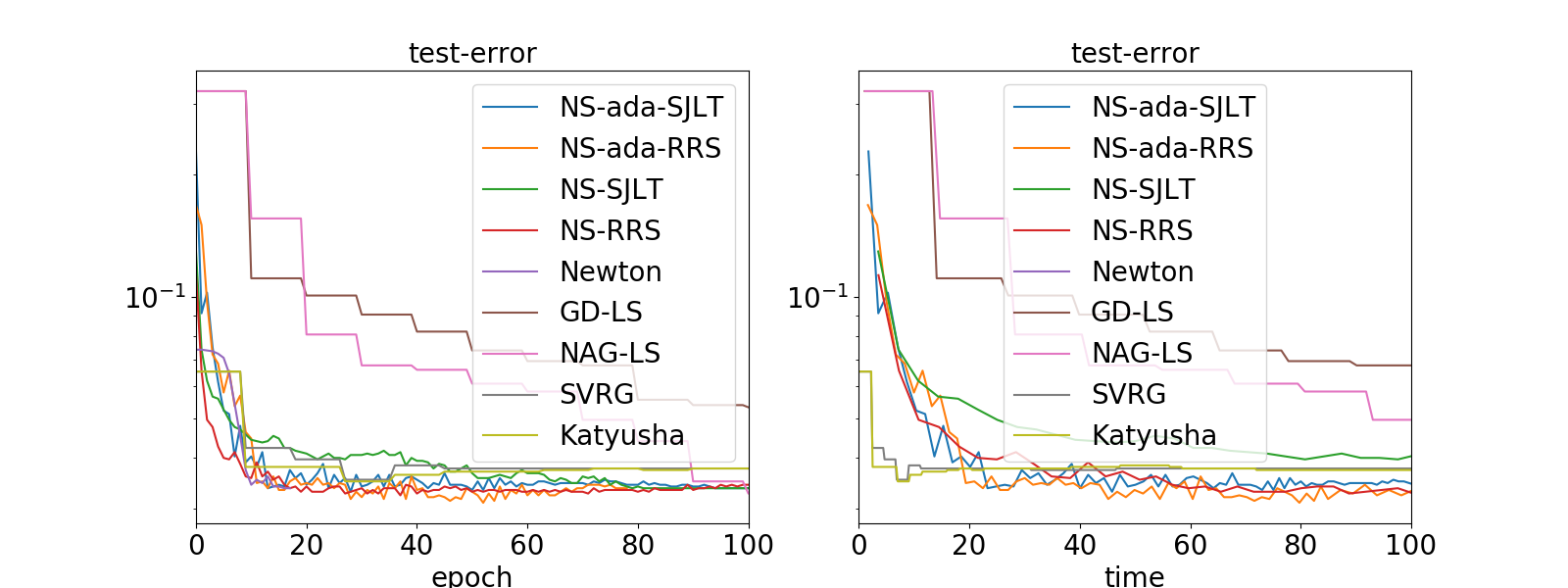}
\end{minipage}
\begin{minipage}[t]{0.45\textwidth}
\centering
\includegraphics[width=\linewidth]{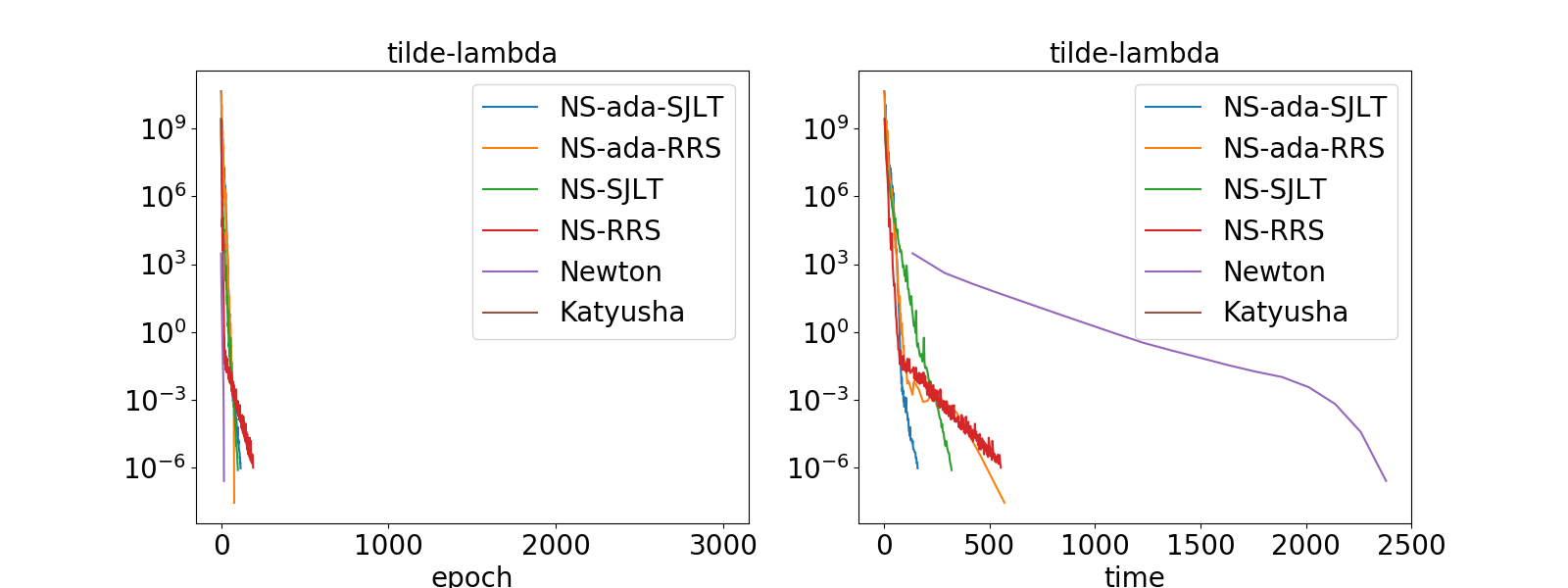}
\end{minipage}
\begin{minipage}[t]{0.45\textwidth}
\centering
\includegraphics[width=\linewidth]{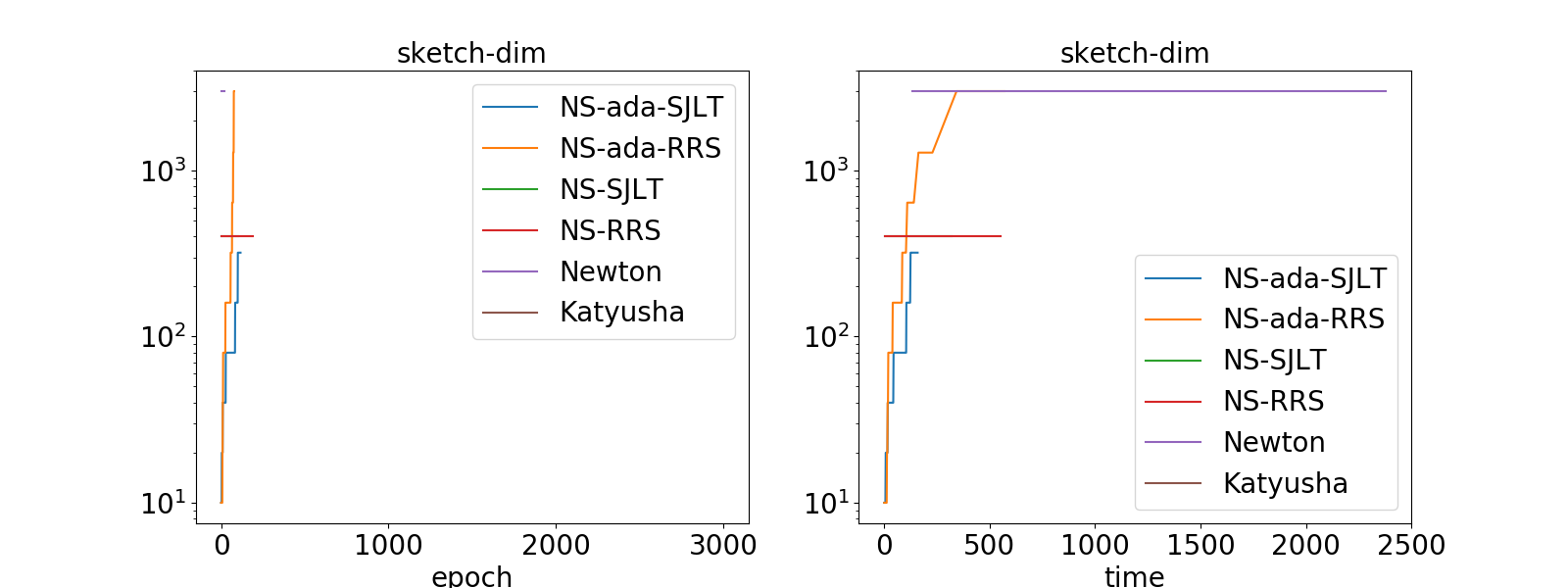}
\end{minipage}
\caption{gisette. $n=3000, d=5000, \mu=10^{-3}$. }
\label{fig:gisette}
\end{figure}

\begin{figure}[htbp]
\centering
\begin{minipage}[t]{0.45\textwidth}
\centering
\includegraphics[width=\linewidth]{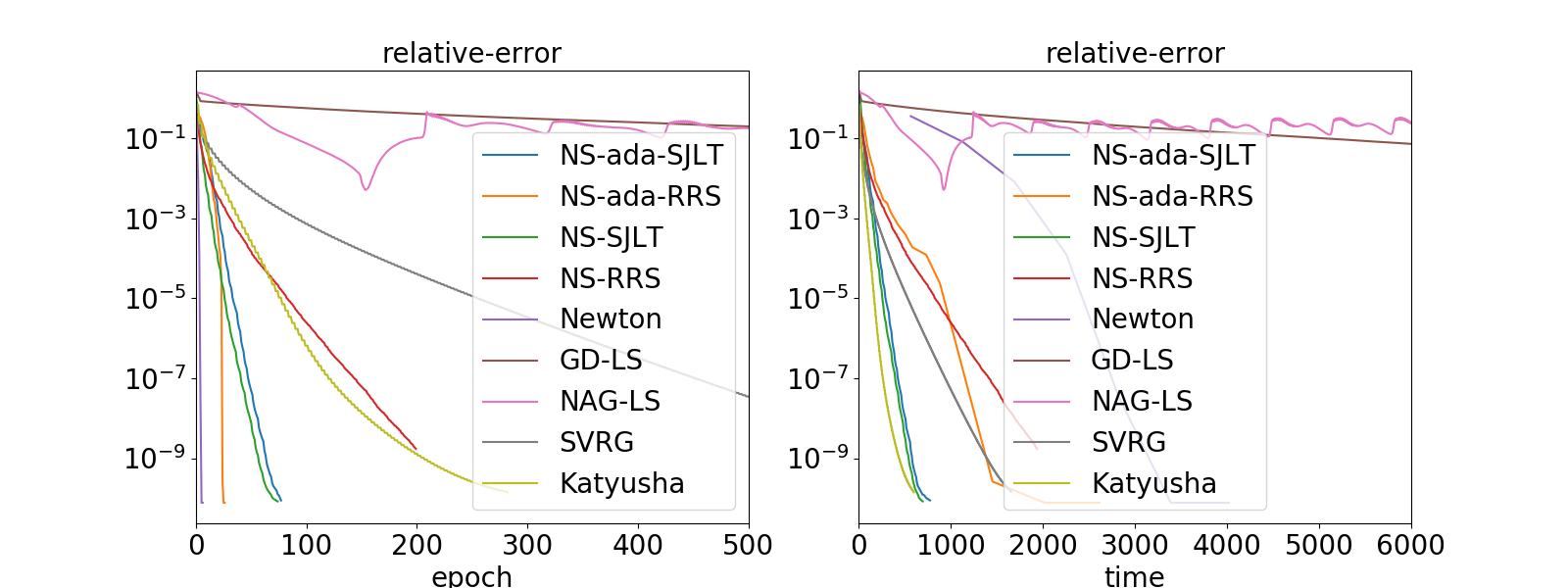}
\end{minipage}
\begin{minipage}[t]{0.45\textwidth}
\centering
\includegraphics[width=\linewidth]{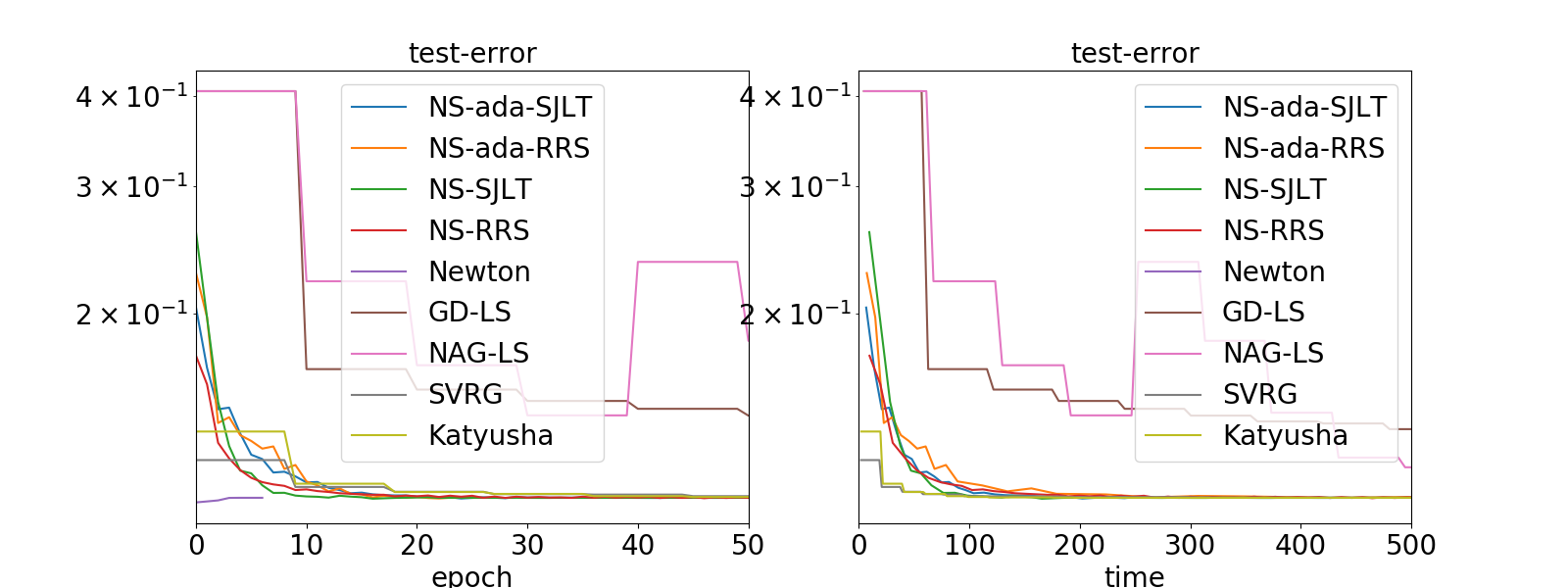}
\end{minipage}
\begin{minipage}[t]{0.45\textwidth}
\centering
\includegraphics[width=\linewidth]{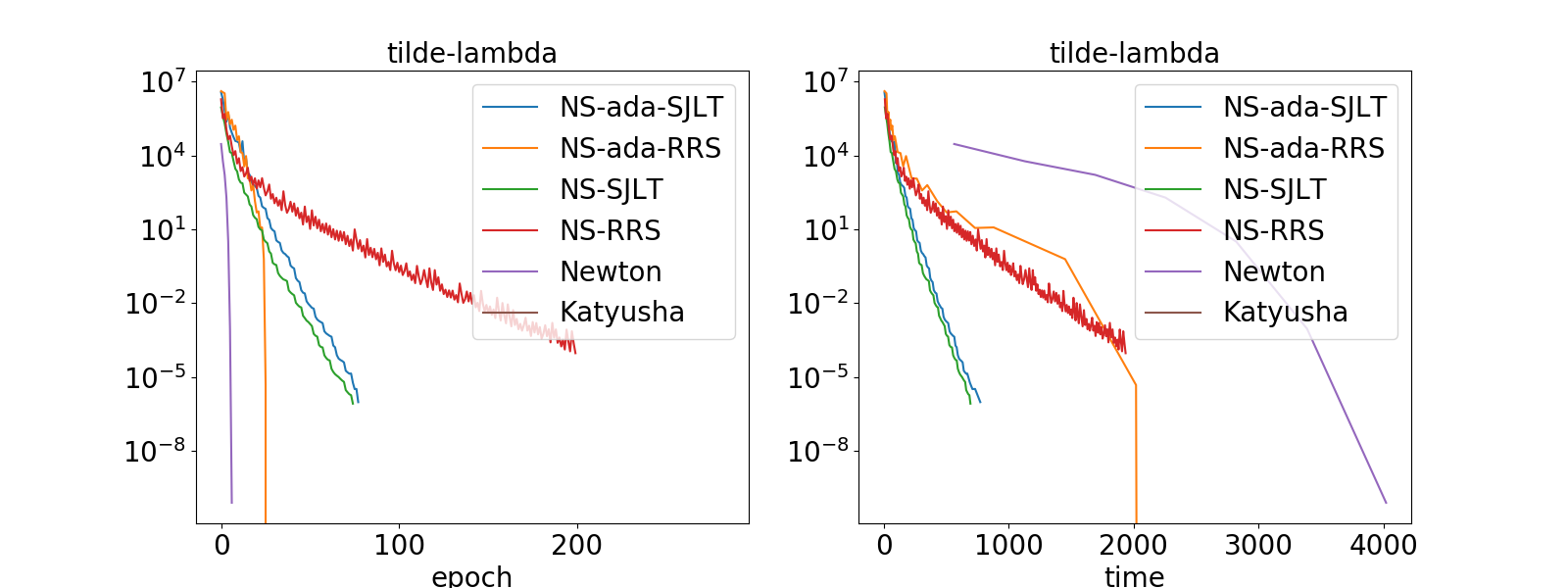}
\end{minipage}
\begin{minipage}[t]{0.45\textwidth}
\centering
\includegraphics[width=\linewidth]{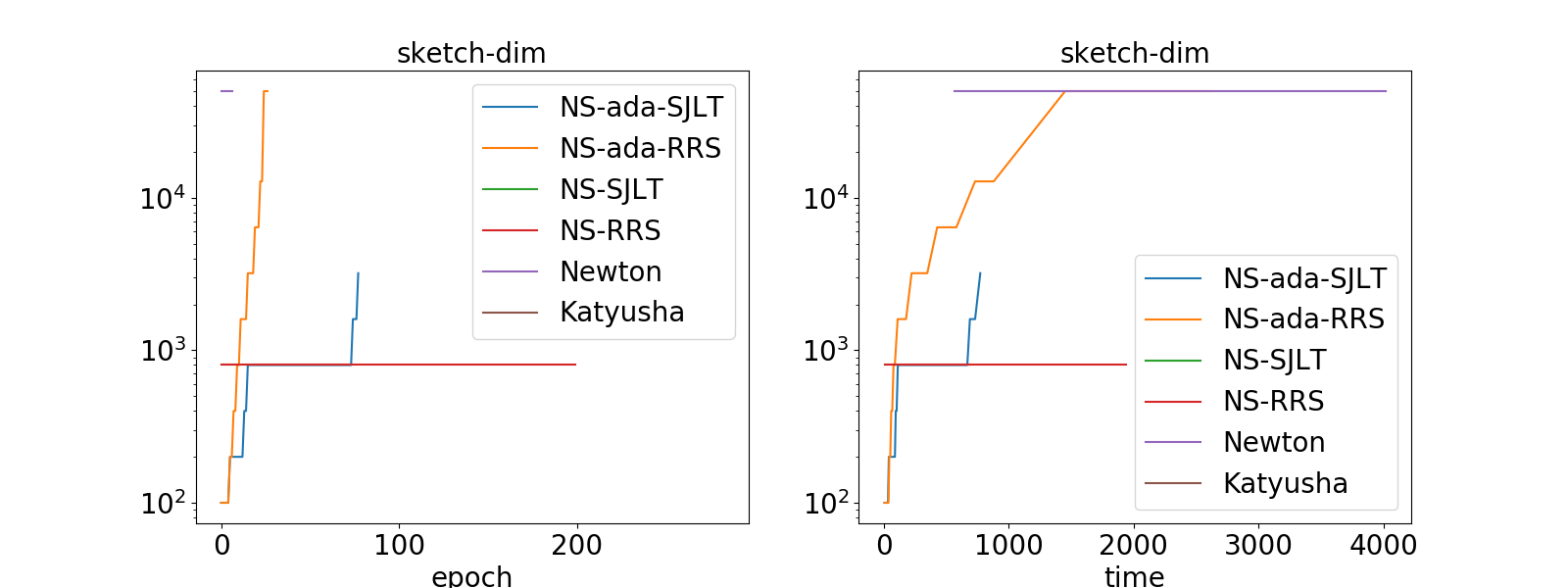}
\end{minipage}
\caption{epsilon. $n=50000,d=2000,\mu=10^{-1}$.}
\label{fig:epsilon}
\end{figure}

We present numerical performance of compared methods with additional details and additional numerical experiments in Figures~\ref{fig:rcv1} to~\ref{fig:epsilon}. Comparatively, NS-ada-RRS tends to have larger sketching dimension than NS-ada-SJLT. This may come from that NS-RRS has stronger oscillations than NS-SJLT in the plot of $\tilde \lambda_f(x^t)$. Thus, NSN-ada-RRS can be slower than NS-ada-SJLT in some test cases where $n$ is not significantly larger than $d$.  

For kernelized regularized logistic regression, the data matrices $A$ and $\tilde A$ are constructed as kernel matrices based on the original data features. Namely, it follows
$$
A_{i,j}=k(\tilde a_i, \tilde a_j), \quad A^\mathrm{test}_{i,j}=k(\tilde a_i^\mathrm{test}, \tilde a_j),
$$
where $\{\tilde a_i\}_{i=1}^n$ and $\{\tilde a_j^\mathrm{test}\}_{i=1}^{n_\mathrm{text}}$ are original data features from the training set and test set respectively. Here $k(x,x'):\mbR^d\times \mbR^d\to \mbR$ is a positive kernel function. We use the isotropic Gaussian kernel function:
$$
k(x,x')=(2\pi h)^{-d/2}\exp\left(-\frac{1}{2h}\|x-x'\|_2^2\right),
$$
where $h>0$ is the bandwidth. We set $h=10$ for a8a dataset and $h=20$ for w7a dataset. For NS-ada-SJLT and NS-ada-RRS, we let $c_1=0.5,\tau=0$ and $c_2=1$. For NS, the sketching dimensions are summarized in Table \ref{tab:NS-k}.

\begin{table}[htbp]
    \centering
    \begin{tabular}{|c|c|c|}
    \hline
         Dataset & $m$ (SJLT) & $m$ (RRS) \\\hline
         a8a-kernel & 100&800\\\hline
         w7a-kernel & 100&800\\ \hline
    \end{tabular}
    \caption{Sketching dimensions of Newton Sketch. kernel matrix.}
    \label{tab:NS-k}
\end{table}

\begin{figure}[htbp]
\centering
\begin{minipage}[t]{0.45\textwidth}
\centering
\includegraphics[width=\linewidth]{figures/a8a-kernel-n10000-d10000-mu1.0e+01/relative-error.png}
\end{minipage}
\begin{minipage}[t]{0.45\textwidth}
\centering
\includegraphics[width=\linewidth]{figures/a8a-kernel-n10000-d10000-mu1.0e+01/test-acc.png}
\end{minipage}
\begin{minipage}[t]{0.45\textwidth}
\centering
\includegraphics[width=\linewidth]{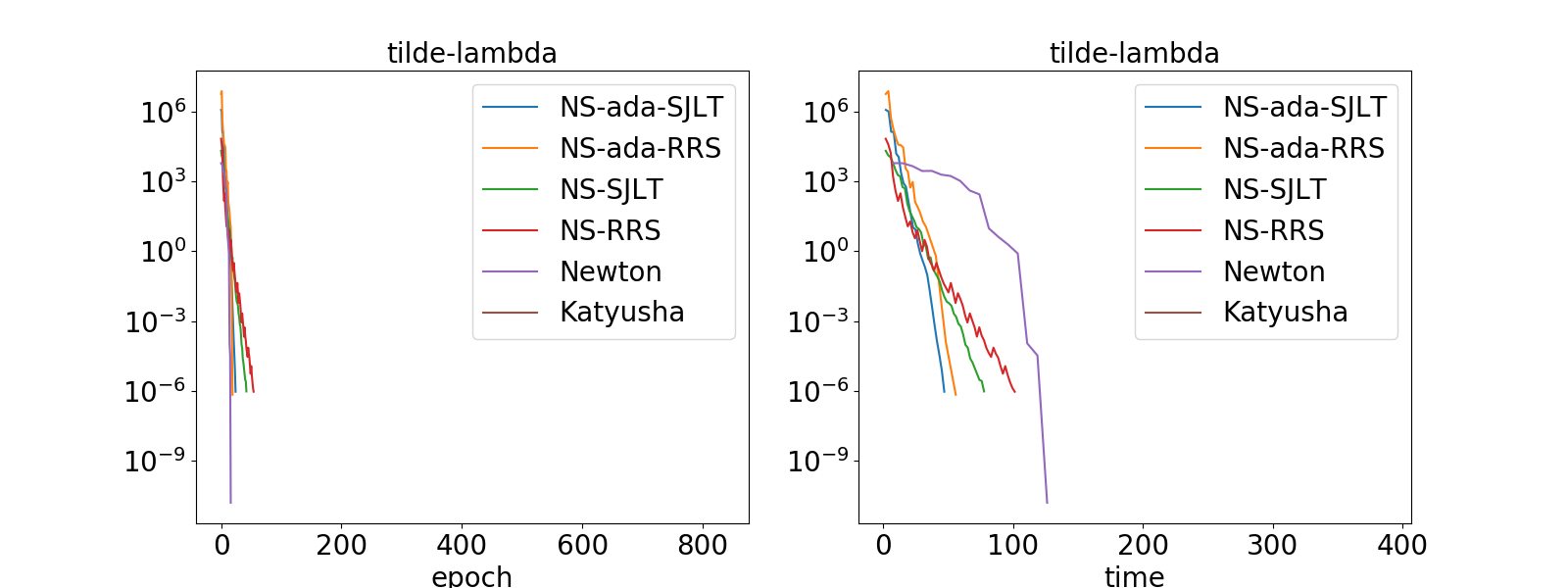}
\end{minipage}
\begin{minipage}[t]{0.45\textwidth}
\centering
\includegraphics[width=\linewidth]{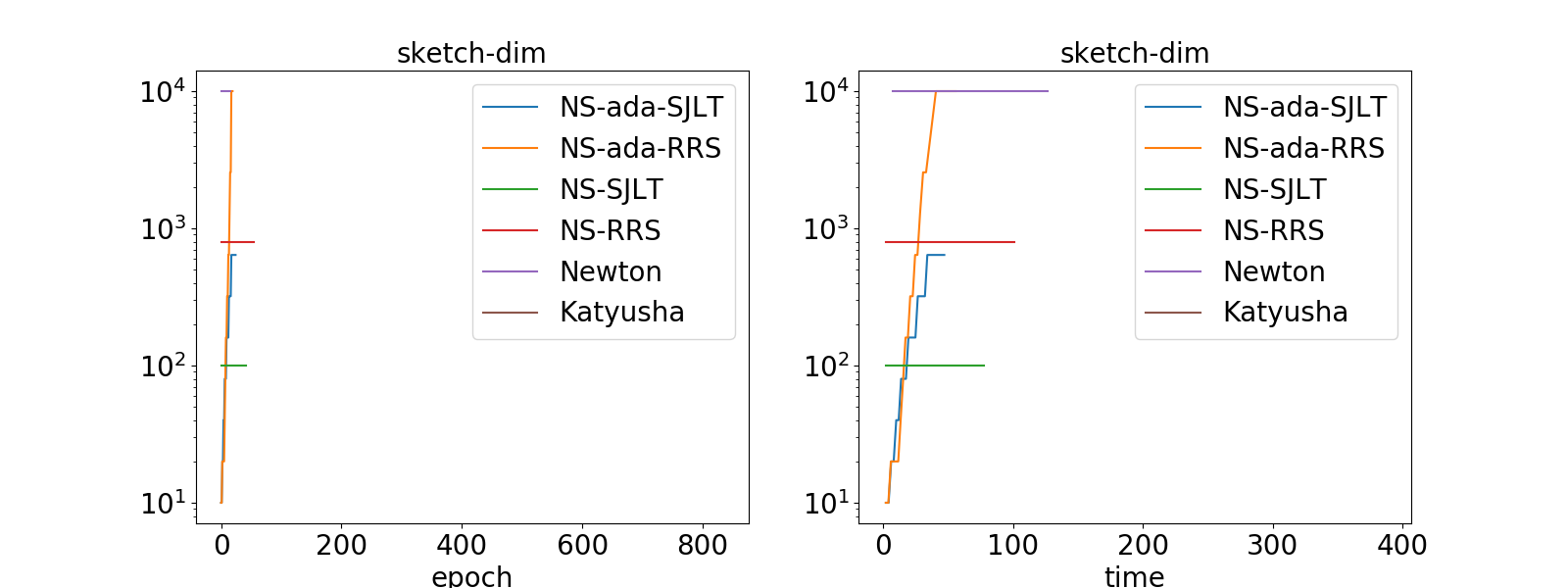}
\end{minipage}
\caption{a8a. kernel matrix. $n=10000,d=10000,\mu=10$.} 
\label{fig:a8a-add}
\end{figure}

\begin{figure}[htbp]
\centering
\begin{minipage}[t]{0.45\textwidth}
\centering
\includegraphics[width=\linewidth]{figures/w7a-kernel-n12000-d12000-mu1.0e+01/relative-error.png}
\end{minipage}
\begin{minipage}[t]{0.45\textwidth}
\centering
\includegraphics[width=\linewidth]{figures/w7a-kernel-n12000-d12000-mu1.0e+01/test-acc.png}
\end{minipage}
\begin{minipage}[t]{0.45\textwidth}
\centering
\includegraphics[width=\linewidth]{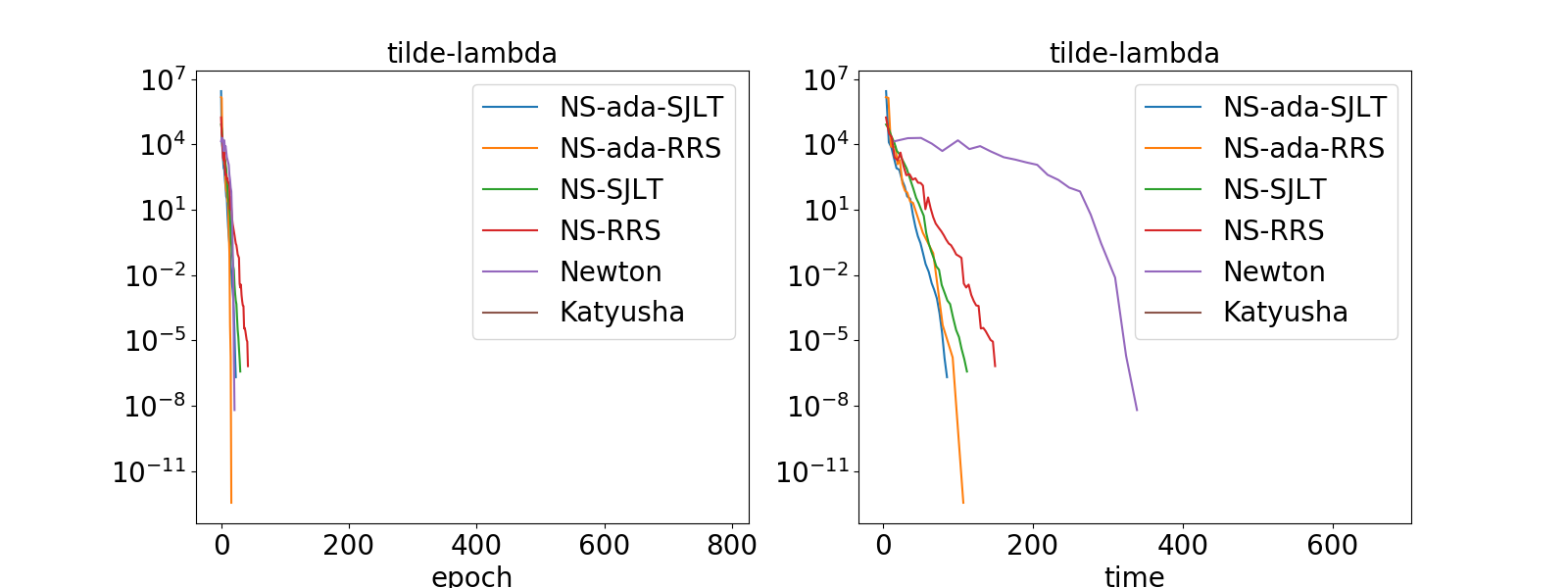}
\end{minipage}
\begin{minipage}[t]{0.45\textwidth}
\centering
\includegraphics[width=\linewidth]{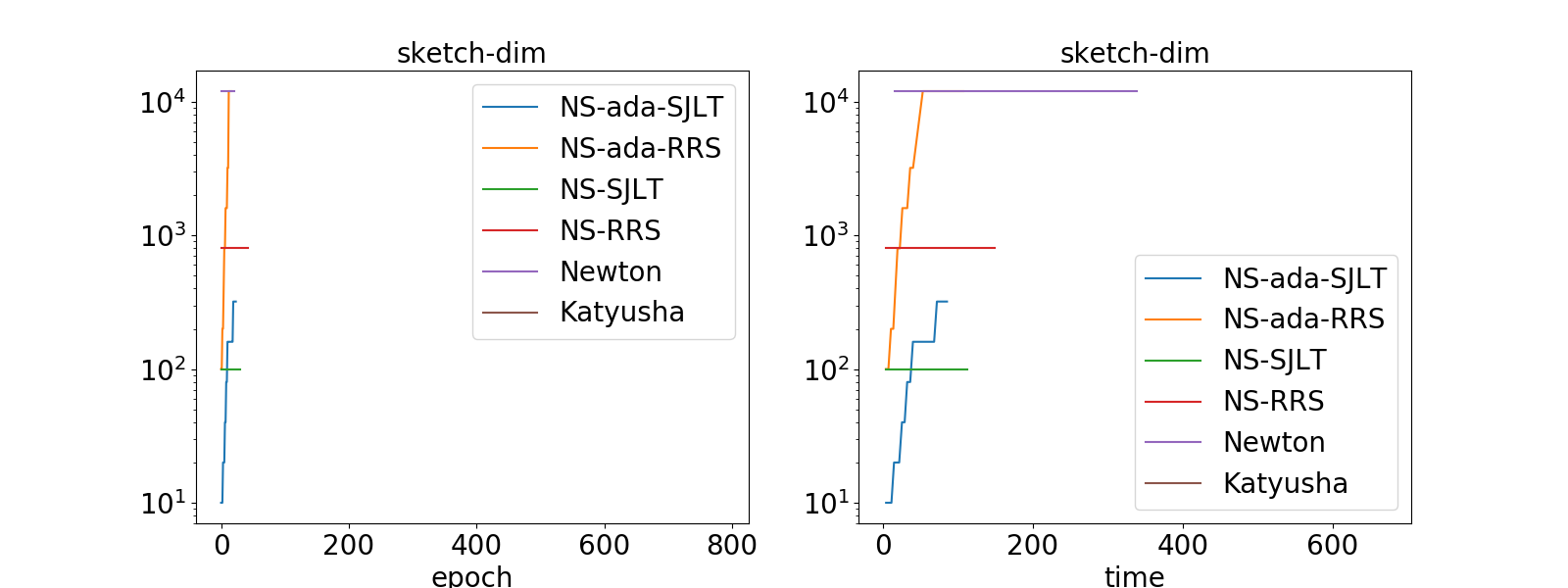}
\end{minipage}
\caption{w7a. kernel matrix. $n=12000,d=12000,\mu=10$.}\label{fig:w7a-add}
\end{figure}

We present numerical results with additional details in Figures~\ref{fig:a8a-add} and~\ref{fig:w7a-add}. We can also observe super linear convergence rate of NS-ada in the plot of $\tilde \lambda_f(x^t)$ when $x^t$ is close to the optimum of the optimization problem. Similarly, NS-ada-RRS tends to have larger sketching dimension than NS-ada-SJLT.

\section{Proof of main results}

\subsection{Proof of Lemma~\ref{lemmaconcentration}}

Let $x \in \textbf{dom}\,f$. We use the shorthand $A \defn \nabla^2 f_0(x)^{1/2}$, and we let $A = U \Sigma V^\top$ be a thin SVD of $A$. We denote by $H^{1/2}$ an invertible square-root matrix of the Hessian $H \equiv H(x) = A^\top A + \nabla^2 g(x)$. Recall that $H_S \equiv H_S(x) = A^\top S^\top S A + \nabla^2 g(x)$. Then, we have
\begin{align*}
    C_S = H^{-\frac{1}{2}} H_S H^{-\frac{1}{2}} &= H^{-\frac{1}{2}} (H + (H_S- H )) H^{-\frac{1}{2}}\\
    &= I_d + H^{-1/2}(H_S- H)H^{-1/2}\\
    &= I_d + H^{-1/2}V \Sigma (U^\top S^\top SU - I_d) \Sigma V^\top H^{-1/2}\,.
\end{align*}
We use the shorthand $M \defn \Sigma V^\top H^{-1/2}$. Using the fact that $\nabla^2 g(x) \succeq \mu \, I_d$, it follows that
\begin{align}
\label{eqnfrobnormeffdim}
    \|M\|_F^2 = \mbox{trace}(\Sigma V^\top H^{-1} V \Sigma) \less \mbox{trace}(\Sigma V^\top (A^\top A + \mu I_d)^{-1} V \Sigma) = d_\mu(x)\,.
\end{align}
It remains to control the spectral norm of $M^\top (U^\top S^\top SU - I_d) M$.

\textbf{(SJLT).} It was shown in~\cite{nelson2013osnap} that for $\varepsilon > 0$ and $p \in (0,1/2)$, it holds with probability at least $1-p$ that $\|M^\top (U^\top S^\top SU - I_d) M\|_2 \less \varepsilon$ provided that $m \gre c_0 \frac{\|M\|_F^4}{\varepsilon^2 p}$, where $c_0 > 0$ is a universal constant. Note that this lower bound on the sketch size is increasing as a function of $\|M\|_F^2$. From inequality~\eqref{eqnfrobnormeffdim}, it is then sufficient to have $m \gre c_0 \frac{d_\mu(x)^2}{\varepsilon^2 p}$ for the above inequality to hold with probability at least $1-p$.

\textbf{(SRHT).} According to Theorems~1 and~9 in~\cite{cohen2015optimal}, it holds with probability at least $1-p$ that $\|M^\top (U^\top S^\top SU - I_d) M\|_2 \less \varepsilon$ provided that $m \gre c_0 \, \varepsilon^{-2} \Big(\|M\|_F^2+\log(\frac{1}{\varepsilon p}) \log(\|M\|_F^2/p) \Big)$, where $c_0$ is a universal constant. Note that this lower bound on the sketch size is increasing as a function of $\|M\|^2_F$. From inequality~\eqref{eqnfrobnormeffdim}, it is then sufficient to have $m \gre c_0 \, \varepsilon^{-2} \Big(d_\mu(x)+\log(\frac{1}{\varepsilon p}) \log(d_\mu(x)/p) \Big)$ for the above inequality to hold with probability at least $1-p$.\qed

\subsection{Proof of Theorem~\ref{theoremclosenessnewtondecrements}}

Let $x \in \textbf{dom}\,f$. Plugging-in the definitions of $\vne$ and $\vnsk$, we have 
\begin{align*}
    \|\vne-\vnsk\|_{H(x)} = \|H^{1/2} (\vne - \vnsk)\|_2 &= \|H^{1/2} (H_S^{-1} \gx - H^{-1}\gx)\|_2\\
    &= \|(H^{1/2} H_S^{-1} H^{1/2} - I_d) H^{-1/2} \gx\|_2\\
    & \less \|C_S^{-1} - I_d\|_2 \, \|H^{-1/2} \gx\|_2\,.
\end{align*}
Using that $\|H^{-1/2} \gx\|_2 = \|\vne\|_{H(x)}$, we further obtain 
\begin{align*}
    \|\vne-\vnsk\|_{H(x)} \less \|C_S^{-1} - I_d\|_2 \, \|\vne\|_{H(x)}\,.
\end{align*}
Under the event $\mathcal{E}_{x,m,\varepsilon}$, it holds for $\varepsilon \in (0,1/4)$ that $(1+\varepsilon/2)^{-1} I_d \preceq C_S^{-1} \preceq (1-\varepsilon/2)^{-1} I_d$. Using the facts that $(1+\varepsilon/2)^{-1} \gre 1-\varepsilon$ and $(1-\varepsilon/2)^{-1} \less 1+\varepsilon$, we obtain the inequality $\|C_S^{-1} - I_d\|_2 \less \varepsilon$, whence
\begin{align*}
    \|\vne-\vnsk\|_{H(x)} \less \varepsilon \, \|\vne\|_{H(x)}\,,
\end{align*}
which proves the first inequality of Theorem~\ref{theoremclosenessnewtondecrements}. On the other hand, we have
\begin{align*}
    \wtilde \lambda_f(x)^2 = \lra{\gx, H_S^{-1} \gx} &= \lra{H^{-\frac{1}{2}} \gx, H^{\frac{1}{2}} H_S^{-1} H^\frac{1}{2} H^{-\frac{1}{2}} \gx}\\
    & = \|C_S^{-\frac{1}{2}} H^{-\frac{1}{2}} \gx\|_2\,.
\end{align*}
It follows that
\begin{align*}
    \frac{1}{\sigma_\text{max}(C_S)} \, \lambda_f(x)^2 \less \wtilde \lambda_f(x)^2 \less \frac{1}{\sigma_\text{min}(C_S)} \, \lambda_f(x)^2\,.
\end{align*}
Conditional on the event $\mathcal{E}_{x,m,\varepsilon}$ and using that $(1+\varepsilon/2)^{-1} \gre 1-\varepsilon$ and $(1-\varepsilon/2)^{-1} \less 1+\varepsilon$, we obtain the claimed result, i.e.,
\begin{align*}
    (1-\varepsilon) \, \lambda_f(x)^2 \less \wtilde \lambda_f(x)^2 \less (1+\varepsilon) \, \lambda_f(x)^2\,.
\end{align*}
\qed

\subsection{Proof of Lemma~\ref{lemmadecreasefirstphase}}

Our proof of this result closely follows the steps of the proof of Lemma~3(a) in~\cite{pilanci2017newton}: the core arguments are the same, but we adapt the proof to our technical framework, that is, conditional on the event $\mathcal{E}_{x,m,\varepsilon}$.

The strategy of the proof is to show that the backtracking line search leads to a step size $s > 0$ such that $f(\xnsk) - f(x) \less - \nu$. We define the univariate function $g(u) \defn f(x + u \vnsk)$ and we set $\varepsilon^\prime = \frac{2\varepsilon}{1-\varepsilon}$. We first show that $\hat u = \frac{1}{1+(1+\varepsilon^\prime)\wtilde \lambda_f(x)}$ satisfies the bound
\begin{align}
\label{eqndecreasefirstphaseintermediate1}
    g(\hat u) \less g(0) - a \hat u \wtilde \lambda_f^2(x)\,,
\end{align}
which implies that $\hat u$ satisfies the exit condition of backtracking line search. Therefore, the step size $s$ must be lower bounded as $s \gre b \hat u$, which further implies that the new iterate $\xnsk = x + s \vnsk$ satisfies the decrement bound
\begin{align*}
    f(\xnsk) - f(x) \less -ab \, \frac{\wtilde \lambda_f(x)^2}{1 + (1+\frac{2\varepsilon}{1-\varepsilon}) \wtilde \lambda_f(x)}\,.
\end{align*}
By assumption, $\wtilde \lambda_f(x) > \eta$. Using the fact that the function $u \mapsto \frac{u^2}{1+(1+\frac{2\varepsilon}{1-\varepsilon}u)}$ is monotone increasing, we get that
\begin{align*}
    f(\xnsk) - f(x) \less -ab \, \frac{\eta^2}{1+(1+\frac{2\varepsilon}{1-\varepsilon})\eta} = \nu\,,
\end{align*}
which is exactly the claimed result. It remains to prove the claims~\eqref{eqndecreasefirstphaseintermediate1}.

According to Lemma 4 in~\cite{pilanci2017newton}, we have for any $u \gre 0$ and $\gamma \gre 0$ that
\begin{align}
\label{eqnlemma4newtonsketch}
    g(u) \less g(0) - u \wtilde \lambda_f(x)^2 - \gamma - \log\!\left(1-\gamma\right)\,,
\end{align}
provided that $u \|\vnsk\|_{H(x)} \less \gamma < 1$. By assumption, the event $\mathcal{E}_{x,m,\varepsilon}$ holds true. As a consequence of Theorem~\ref{theoremclosenessnewtondecrements}, we have that 
\begin{align*}
    \|\vnsk\|_{H(x)} \less (1+\varepsilon) \lambda_f(x) \less \frac{1+\varepsilon}{1-\varepsilon}\, \wtilde \lambda_f(x) = (1+\varepsilon^\prime) \wtilde \lambda_f(x)\,.
\end{align*}
It follows that $\hat u \|\vnsk\|_{H(x)} \less \hat u (1+\varepsilon^\prime) \wtilde \lambda_f(x) < 1$. Plugging-in $u = \hat u$ and $\gamma = \hat u (1+\varepsilon^\prime \wtilde \lambda_f(x))$ into~\eqref{eqnlemma4newtonsketch}, we obtain that
\begin{align*}
    g(\hat u) &\less g(0) - \hat u \wtilde \lambda_f(x)^2 - \hat u (1+\varepsilon^\prime) \wtilde \lambda_f(x) - \log(1-\hat u(1+\varepsilon^\prime)\wtilde \lambda_f(x))\\
    & = g(0) - \left\{\hat u (1+\varepsilon^\prime)^2 \wtilde \lambda_f(x)^2 + \hat u(1+\varepsilon^\prime)\wtilde \lambda_f(x) + \log(1-\hat u(1+\varepsilon^\prime)\wtilde \lambda_f(x)) - \hat u ((1+\varepsilon^\prime)^2 -1) \wtilde \lambda_f(x)^2\right\}\,.
\end{align*}
Using that $\hat u (1+\varepsilon^\prime)^2 \wtilde \lambda_f(x)^2 + \hat u(1+\varepsilon^\prime)\wtilde \lambda_f(x) = (1+\varepsilon^\prime) \wtilde \lambda_f(x)$ and $\hat u ((1+\varepsilon^\prime)^2 -1) \wtilde \lambda_f(x)^2 = \frac{({\varepsilon^\prime}^2 + 2\varepsilon^\prime)\wtilde \lambda_f(x)^2}{1+(1+\varepsilon^\prime) \wtilde\lambda_f(x)}$, we find that 
\begin{align*}
    g(\hat u) \less g(0) - (1+\varepsilon^\prime) \wtilde \lambda_f(x) + \log(1 + (1+\varepsilon^\prime)\wtilde \lambda_f(x)) + \frac{({\varepsilon^\prime}^2 + 2\varepsilon^\prime)\wtilde \lambda_f(x)^2}{1+(1+\varepsilon^\prime) \wtilde\lambda_f(x)}\,.
\end{align*}
Applying the inequality $-z + \log(1+z) \less -\frac{1}{2}\frac{z^2}{(1+z)}$ with $z=(1+\varepsilon^\prime) \wtilde \lambda_f(x)$, we further obtain that 
\begin{align*}
    g(\hat u) &\less g(0) - \frac{\frac{1}{2}(1+\varepsilon^\prime)^2 \wtilde \lambda_f(x)^2}{1+(1+\varepsilon^\prime) \wtilde \lambda_f(x)} + \frac{({\varepsilon^\prime}^2 + 2\varepsilon^\prime) \wtilde \lambda_f(x)^2}{1+(1+\varepsilon^\prime)\wtilde \lambda_f(x)}\\
    & = g(0) - \left(\frac{1}{2} - \frac{{\varepsilon^\prime}^2}{2} - \varepsilon^\prime\right) \wtilde \lambda_f(x)^2 \hat u\\
    & \less g(0) - a \wtilde \lambda_f(x)^2 \hat u\,,
\end{align*}
where the final inequality follows by the assumption that $a \less 1 - \frac{1}{2}\left(\frac{1+\varepsilon}{1-\varepsilon}\right)^2$, that is, $a \less \frac{1}{2} - \frac{{\varepsilon^\prime}^2}{2} - \varepsilon^\prime$. This concludes the proof. \qed

\subsection{Proof of Lemma~\ref{lemmadecreasesecondphase}}

We recall Theorem 4.1.6~of~\cite{nesterov2003introductory} (see, also,~Exercise~9.17 in~\cite{boyd2004convex}): it guarantees that for a step size $s > 0$ such that $|1-s\|\vnsk\|_{H(x)}| < 1$, we have
\begin{align}
\label{eqnfundamentalnesterov}
    (1 - s \|\vnsk\|_{H(x)})^2 \, H(x) \preceq H(x + s \vnsk) \preceq \frac{1}{(1-s\|\vnsk\|_{H(x)})^2} \, H(x)\,.    
\end{align}
By assumption, the event $\mathcal{E}_{x, m, \varepsilon^\prime}$ holds. As a consequence of Theorem~\ref{theoremclosenessnewtondecrements}, we have $\|\vnsk\|_{H(x)} \less (1+\varepsilon^\prime) \|\vne\|_{H(x)}$. Plugging this bound into~\eqref{eqnfundamentalnesterov} and using that $\|\vne\|_{H(x)} = \lambda_f(x)$, we obtain
\begin{align}
\label{eqnsandwichrelation}
    (1 - s (1+\varepsilon^\prime) \lambda_f(x))^2 \, H(x) \preceq H(x + s \vnsk) \preceq \frac{1}{(1-s(1+\varepsilon^\prime) \lambda_f(x))^2} \, H(x)\,,      
\end{align}
for $s > 0$ such that $s (1+\varepsilon^\prime) \lambda_f(x) < 1$. Denote by $s_\mathrm{nsk}$ the step size obtained by backtracking line search. It satisfies $s_\mathrm{nsk} = 1$. Then, it holds that 
\begin{align*}
    s_\mathrm{nsk}(1+\varepsilon^\prime) \lambda_f(x) \less (1+\varepsilon^\prime) \lambda_f(x) & \underset{(i)}{\less} \frac{1+\varepsilon^\prime}{\sqrt{1-\varepsilon^\prime}} \, \wtilde \lambda_f(x) \\
    & \underset{(ii)}{\less} \frac{1+\varepsilon^\prime}{\sqrt{1-\varepsilon^\prime}} \, \eta\\
    & \underset{(iii)}{<} 1\,,
\end{align*}
where inequality (i) follows from the assumption that $\mathcal{E}_{x,m,\varepsilon^\prime}$ holds and from Theorem~\ref{theoremclosenessnewtondecrements}; inequality (ii) follows from the assumption that $\wtilde \lambda_f(x) \less \eta$. Furthermore, we have $\varepsilon^\prime \less \varepsilon < 1/4$, as well as $\eta < 1/16$ (see Lemma~\ref{lemmaboundeta}) and this yields inequality (iii).

Using~\eqref{eqnsandwichrelation}, we then obtain that
\begin{align*}
    \lambda_f(\xnsk) &= \|H(\xnsk)^{-1/2} \nabla f(\xnsk)\|_2\\
    &\less \frac{1}{(1-(1+\varepsilon^\prime)\lambda_f(x))} \, \|H(x)^{-1/2} \nabla f(\xnsk)\|_2\\
    &= \frac{1}{(1-(1+\varepsilon^\prime)\lambda_f(x))} \, \left\|H(x)^{-1/2} \left( \nabla f(x) + \int_0^1 H(x+s \vnsk) \vnsk \mathrm{d}s \right) \right\|_2\\
    &\less \frac{1}{(1-(1+\varepsilon^\prime)\lambda_f(x))} \, (M_1 + M_2)\,,
\end{align*}
where 
\begin{align*}
    &M_1 = \left\|H(x)^{-1/2} \left( \nabla f(x) + \int_0^1 H(x+s \vnsk) \vne \mathrm{d}s \right) \right\|_2\,,\\
    &M_2 = \left\|H(x)^{-1/2} \, \int_0^1 H(x+s \vnsk) (\vnsk - \vne) \mathrm{d}s\right\|_2\,.
\end{align*}
It remains to bound the terms $M_1$ and $M_2$. Regarding $M_1$, we have after re-arranging and using inequality~\eqref{eqnsandwichrelation} that
\begin{align*}
    M_1 &= \left\| \int_0^1 \left(H(x)^{-1/2} H(x + s\vnsk) H(x)^{-1/2} - I_d \right) \mathrm{d}s \, H(x)^{1/2} \vne \right\|_2\\
    & \less \left| \int_0^1 \frac{1}{(1-s(1+\varepsilon^\prime) \lambda_f(x))^2}\mathrm{d}s - 1\right| \, \left\| H(x)^{1/2} \vne \right\|_2\\
    & = \frac{(1+\varepsilon^\prime) \lambda^2_f(x)}{1-(1+\varepsilon^\prime)\lambda_f(x)}\,.
\end{align*}
Regarding $M_2$, we have
\begin{align*}
    M_2 &= \left\| \int_0^1 H(x)^{-1/2} H(x + s\vnsk) H(x)^{-1/2} \mathrm{d}s \, H(x)^{1/2} (\vnsk - \vne)\right\|_2\\
    & \less \left\| \int_0^1 \frac{1}{(1-s(1+\varepsilon^\prime) \lambda_f(x))^2}\mathrm{d}s \, H(x)^{1/2} (\vnsk - \vne) \right\|_2\\
    &= \frac{1}{1- (1+\varepsilon^\prime)\lambda_f(x)} \, \left\|H(x)^{1/2} (\vnsk - \vne)\right\|_2\\
    & \less \frac{\varepsilon^\prime \lambda_f(x)}{1-(1+\varepsilon^\prime)\lambda_f(x)}\,,
\end{align*}
where the last inequality follows from the assumption that the event $\mathcal{E}_{x,m,\varepsilon^\prime}$ holds and as a consequence of Theorem~\ref{theoremclosenessnewtondecrements}. Plugging these bounds on $M_1$ and $M_2$, we obtain that 
\begin{align}
\label{eqnintermediateboundnewtondecrement}
    \lambda_f(\xnsk) \less \frac{(1+\varepsilon^\prime) \lambda_f(x)^2 + \varepsilon^\prime \lambda_f(x)}{(1-(1+\varepsilon^\prime)\lambda_f(x))^2}\,.
\end{align}
Recall that $\varepsilon^\prime \less \varepsilon \, \lambda_f(x)^\tau$. Combining this inequality with~\eqref{eqnintermediateboundnewtondecrement}, we obtain
\begin{align*}
    \lambda_f(\xnsk) \less \frac{(1+\varepsilon \,\lambda_f(x)^\tau) \, \lambda_f(x)^2 + \varepsilon \, \lambda_f(x)^{1+\tau}}{(1-(1+\varepsilon \, \lambda_f(x)^\tau)\lambda_f(x))^2}
    & = \underbrace{\frac{\lambda_f(x)^{1-\tau}+\varepsilon \,\lambda_f(x) + \varepsilon}{(1-(1+\varepsilon \, \lambda_f(x)^\tau)\lambda_f(x))^2}}_{\defn \alpha(\tau,x)} \, \lambda_f(x)^{1+\tau}\,.
\end{align*}
On the event $\mathcal{E}_{x,m,\varepsilon^\prime}$, we have according to Theorem~\ref{theoremclosenessnewtondecrements} that $(1+\varepsilon)\lambda_f(x) \less \frac{(1+\varepsilon) \wtilde \lambda_f(x)}{\sqrt{1-\varepsilon}} \less \frac{(1+\varepsilon)\eta}{\sqrt{1-\varepsilon}} \less \frac{1}{16}$, where the last inequality follows from Lemma~\ref{lemmaboundeta}. Hence, the denominator of $\alpha(\tau,x)$ satisfies
\begin{align*}
    1-(1+\varepsilon \lambda_f(x)^\tau) \lambda_f(x) \gre 1- (1+\varepsilon) \lambda_f(x) \gre \frac{15}{16}\,,
\end{align*}
while the numerator of $\alpha(\tau,x)$ satisfies 
\begin{align*}
    \lambda_f(x)^{1-\tau}+\varepsilon \,\lambda_f(x) + \varepsilon \less \frac{1}{16^{1-\tau}} + \frac{1}{32} + \frac{1}{2}
\end{align*}
Combining these bounds together, we obtain that
\begin{align*}
    \alpha(\tau, x) \less \frac{8 + 1/2 + 16^\tau}{15} \less 0.57 + \frac{16^\tau}{15} = \alpha(\tau)\,.
\end{align*}
It is easy to verify that $\alpha(\tau)^{1/\tau} \less 2$ for any $\tau \in (0,1]$. Furthermore, for $\tau = 0$, we obtain that $\alpha(0) \approx 0.63333 \less 0.64 = \frac{16}{25}$, and this concludes the proof. Note that a similar linear convergence rate was obtained for the Newton sketch provided that $m \gtrsim d$ (see Lemma~3 in~\cite{pilanci2017newton}).\qed

\subsection{Proof of Lemma~\ref{lemmageometricdecrease}}

By induction, we obtain for any $t \gre 0$ that $\alpha^\frac{1}{\tau} \beta_t \less (\alpha^\frac{1}{\tau} \eta)^{(1+\tau)^t}$. To have $\beta_t \less \sqrt{\delta}$, it suffices that $(\alpha^\frac{1}{\tau} \eta)^{(1+\tau)^t} \less \alpha^\frac{1}{\tau} \sqrt{\delta}$. Taking the logarithm on both sides, this yields $(1+\tau)^t \log(\alpha^\frac{1}{\tau} \eta) \less \log(\alpha^\frac{1}{\tau} \sqrt{\delta})$, i.e., $(1+\tau)^t \log(1/\alpha^\frac{1}{\tau} \eta) \gre \log(1/\alpha^\frac{1}{\tau} \sqrt{\delta})$. By assumption, $\log(1/\alpha^\frac{1}{\tau} \eta) > 0$ and $\log(1/\alpha^\frac{1}{\tau} \sqrt{\delta}) > 0$. Therefore, after dividing both sides by $\log(1/\alpha^\frac{1}{\tau} \eta)$ and taking again the logarithm, we find that it is sufficient to have 
\begin{align*}
    t &\gre \ceil{\frac{1}{\log(1+\tau)} \, \log\!\left(\frac{\log(1/\alpha^\frac{1}{\tau} \sqrt{\delta})}{\log(1/\alpha^\frac{1}{\tau} \eta)}\right)}\\
    & = \ceil{\frac{1}{\log(1+\tau)} \, \log\!\left(\frac{1 + \frac{\tau \log(1/\delta)}{2\log(1/\alpha)}}{1 + \frac{\tau \log(1/\eta)}{\log(1/\alpha)}}\right)}\\
    & = T_{\tau, \alpha, \delta}\,.
\end{align*}
\qed

\subsection{Proof of Theorem~\ref{theoremquadraticconvergenceeffdimnewtonsketch}}

We denote $N_1 \defn \frac{f(x_0) - f(x^*)}{\nu}$ and $\wtilde p \defn \frac{p_0}{\overline{T}+2}$, where $\overline{T} \defn N_1 + 1 + T_{\tau,\frac{3}{8}\delta}$. Recall that we pick $\varepsilon = 1/8$.

Our proof strategy proceeds as follows. In a first phase, we show that $f(\xnsk) - f(x) \less -\nu$ until such a decrement cannot occur anymore, i.e., until $f(x_t) - f(x^*) < \nu$. Technical arguments for Phase 1 essentially follow from Lemma~\ref{lemmadecreasefirstphase}. Then, we enter a second phase where we observe a geometric decrease of the Newton decrement as described in Lemma~\ref{lemmadecreasesecondphase}.

We define
\begin{align*}
    t \defn \inf\left\{k \gre 0 \mid \wtilde \lambda_f(x_k)\less \eta\right\}\,,
\end{align*}
According to Lemma~\ref{lemmalengthofphase1}, we have $t \less N_1$ with probability at least $1-N_1 \wtilde p$. 

We turn to the analysis of Phase 2. We suppose that $T_f > t$ (i.e., the algorithm has not terminated during Phase 1), we define the additional number of iterations $J \defn \min\{T_{\tau,\frac{3}{8}\delta}, T_f-t-1\}$, and we introduce the event
\begin{align*}
    \mathcal{E}^{(2)} \defn \Big\{\mathcal{E}_{x_t, m_t, \varepsilon} \cap \bigcap_{j=0}^J \mathcal{E}_{x_{t+1+j}, m_{t+1+j}, \varepsilon \delta^{\frac{\tau}{2}} }\Big\}\,.
\end{align*}
Let us assume that $\mathcal{E}^{(2)}$ holds true, which happens with probability at least $1-(2+T_{\tau, \frac{3}{8}\delta}) \wtilde p$ according to Corollary~\ref{corollaryprobae2}. According to Lemma~\ref{lemmaphase2theorem2}, we have for any $j=0, \dots, J$ that $m_{t+1+j} = \overline{m}_2$ and,
\begin{align*}
    \alpha(\tau)^\frac{1}{\tau} \lambda_f(x_{t+1+j}) \less (\alpha(\tau)^{\frac{1}{\tau}} \lambda_f(x_{t+1}) )^{(1+\tau)^j}\,.
\end{align*}
Further, we have from Lemma~\ref{lemmadecreasesecondphase} and Theorem~\ref{theoremclosenessnewtondecrements} that $\lambda_f(x_{t+1}) \less \frac{16}{25} \, \lambda_f(x_t) \less \frac{\wtilde \lambda_f(x_t)}{\sqrt{1-\varepsilon}} \less \frac{\eta}{\sqrt{1-\varepsilon}} \less \frac{1}{16}$. Hence, $\alpha(\tau)^{\frac{1}{\tau}} \lambda_f(x_{t+1}) < 1/8$. As a consequence of Lemma~\ref{lemmageometricdecrease}, we must have that $\lambda_f(x_{t+1+j})^2 \less \frac{3}{8}\delta$ for some $j \less T_{\tau, \frac{3}{8}\delta}$, which further implies that
\begin{align*}
    \wtilde \lambda_f(x_{t+1+j})^2 \less (1+\varepsilon) \lambda_f(x_{t+1+j})^2 \less \frac{3(1+\varepsilon)}{8} \delta \less \frac{3}{4} \delta\,.
\end{align*}
The above inequality implies termination of the algorithm before the time $t+1+T_{\tau, \frac{3}{8}\delta}$. Using a union bound over $\{t \less N_1\}$ and $\mathcal{E}^{(2)}$, we find that the algorithm terminates within $N_1 + 1 + T_{\delta, \frac{3}{8}\delta}$ iterations with probability at least $1-(N_1 + 2 + T_{\tau,\frac{3}{8}\delta})\wtilde p$.

It remains to guarantee that the algorithm returns a point $\wtilde x$ such that $f(\wtilde x) - f(x^*) \less \delta$. Note that the exit criterion guarantees that $\wtilde \lambda_f(\wtilde x)^2 \less \frac{3}{4}\delta$. Furthermore, the final sketch size $\wtilde m$ necessarily satisfies $\wtilde m \gre \overline{m}_1$, so that, according to Theorem~\ref{theoremclosenessnewtondecrements}, we have with probability at least $1-\wtilde p$ that $\lambda_f(\wtilde x)^2 \less \frac{1}{1-\varepsilon}\wtilde \lambda_f(\wtilde x)^2 \less \delta$. Self-concordance of $f$ further implies that $f(\wtilde x) - f(x^*) \less \lambda_f(\wtilde x)^2 \less \delta$. 

In conclusion, we have shown that the algorithm returns a $\delta$-accurate solution within $N_1 + 1 + T_{\tau,\frac{3}{8}\delta}$ iterations with probability at least $1-(N_1 + 3 + T_{\tau,\frac{3}{8}\delta})\wtilde p = 1-p_0$. This concludes the proof.\qed

\subsubsection{Complexity guarantees for the SJLT}

With the SJLT, consider the quadratic convergence case, i.e., $\tau = 1$. Let $p_0 > 0$ be a failure probability, and consider the sketch sizes 
\begin{align*}
    \overline{m}_1 \asymp \frac{\overline{d}_\mu^2 \log \log 1/\delta}{p_0}\,,\qquad \overline{m}_2 \asymp \frac{1}{\delta} \, \frac{\overline{d}_\mu^2 \log \log 1/\delta}{p_0}\,.
\end{align*}
We observe quadratic convergence with $T_f = \mathcal{O}(\log\log(\frac{1}{\delta}\big))$ iterations. Further, assuming that the sketching cost $\mathcal{O}(nd)$ dominates the cost $\mathcal{O}(\overline{m}^2 d)$ of solving the randomized Newton system, i.e., $n \gtrsim \frac{\overline{d}_\mu^4 \log(\log(1/\delta))^2}{\delta^2 p_0^2}$, then the total complexity results in 
\begin{align*}
    \mathcal{C} = \mathcal{O}\big( nd \log \log 1/\delta \big)\,.
\end{align*}
Similarly, we consider the linear convergence case, i.e., $\tau=0$, and pick a failure probability $p_0 > 0$. Consider the sketch sizes 
\begin{align*}
    \overline{m}_1 \asymp \overline{m}_2 \asymp \frac{\overline{d}_\mu^2 \log 1/\delta}{p_0}\,.
\end{align*}
We observe linear convergence with $T_f = \mathcal{O}(\log \frac{1}{\delta})$ iterations. Assuming again that the sketching cost dominates the cost of solving the randomized Newton system, i.e., $n \gtrsim \frac{\overline{d}_\mu^4 \log^2(1/\delta)}{p_0^2}$, we obtain the total time complexity
\begin{align*}
    \mathcal{C} = \mathcal{O}\!\left(nd \log(1/\delta)\right)\,.
\end{align*}
\qed

\subsection{Proof of Lemma~\ref{lemmaterminationcondition}}

Let $S \in \real^{m \times n}$ be an embedding, and $C_S \defn H^{-1/2} H_S H^{-1/2}$. We use the notations $A \defn \nabla^2 f_0(x)^{1/2}$, and we let $A = U \Sigma V^\top$ be a thin SVD of $A$. Then, we have
\begin{align*}
    C_S = H^{-\frac{1}{2}} H_S H^{-\frac{1}{2}} &= H^{-\frac{1}{2}} (H + (H_S- H )) H^{-\frac{1}{2}}\\
    &= I_d + H^{-1/2}(H_S- H)H^{-1/2}\\
    &= I_d + M^\top (U^\top S^\top SU - I_d) M\,,
\end{align*}
where $M \defn \Sigma V^\top H^{-1/2}$. According to~\cite{cohen2015optimal}, it holds that $\|M^\top (U^\top S^\top SU - I_d) M\|_2 \less \frac{\overline{d}_\mu}{2}$ (i.e., $\|C_S\|_2 \less 1 + \frac{\overline{d}_\mu}{2}$) with probability at least $1-p$, provided that $m \gre \Omega(\log^2(1/p))$ for a SRHT $S$, and, $m \gre \Omega(1/p)$ for a SJLT $S$.

Then, we use the fact that 
\begin{align*}
    \wtilde \lambda_f(x)^2 = \langle H^{-1/2} \gx, H^{1/2} H_S^{-1} H^{1/2} H^{-1/2} \gx \rangle \gre \frac{1}{\|C_S\|_2} \, \lambda_f(x)^2\,.
\end{align*}
Conditional on $\|C_S\|_2 \less 1+\frac{\overline{d}_\mu}{2}$, it follows that
\begin{align*}
    \lambda_f(x)^2 \less \|C_S\|_2 \, \wtilde \lambda_f(x)^2 \less (1+\frac{\overline{d}_\mu}{2})\frac{\delta}{d} \less \delta\,.    
\end{align*}
Using the self-concordance of $f$, we obtain that $f(x) - f(x^*) \less \delta$. This concludes the proof.\qed

\subsection{Proof of Theorem~\ref{theoremconvergenceadaptive}}

We introduce the notations
\begin{align*}
    \overline{T} = T_{\tau, \alpha(\tau,\varepsilon), \frac{\delta}{d}} + N_1\,,\qquad \wtilde p = \frac{p_0}{\overline{T}} \qquad \mbox{and} \qquad \varepsilon^\prime = \varepsilon \, \left(\frac{\delta}{(1+\varepsilon)d}\right)^{\tau/2}\,.
\end{align*}
We consider $\overline m$ a sketch size such that $\mathcal{E}_{x, \overline m, \varepsilon^\prime}$ holds with probability at least $1-\wtilde p$, that is, 
\begin{align*}
    &\overline m = \Omega\Big(\frac{d^{\tau}\overline d_\mu^2 \overline T}{p_0 \delta^\tau}\Big) \qquad \mbox{for the SJLT}\,,\\ 
    &\overline m = \Omega\!\left( \frac{d^\tau}{\delta^\tau} \Big(\overline{d}_\mu + \log\big(\frac{\overline{T} d^{\tau/2}}{p_0 \delta^{\tau/2}}\big) \log\big(\frac{\overline{d}_\mu \overline{T}}{p_0}\big)\Big) \right) \qquad \mbox{for the SRHT}\,.
\end{align*}
\textbf{Phase 2}. Let $t \gre 0$ be the first iteration such that $m_t \gre \overline m$, if any. Let $x \equiv x_{t+j}$ be an iterate after time $t$, for some $j \gre 0$. The sketch size is non-decreasing, whence $m \equiv m_{t+j} \gre \overline m$. We assume that $\mathcal{E}_{x,m,\varepsilon^\prime}$ holds, and that the algorithm has not yet terminated, i.e., $\wtilde \lambda_f(x)^2 > \delta/d$. Note that $\varepsilon > \varepsilon^\prime$, whence $\mathcal{E}_{x,m,\varepsilon}$ also holds. By Theorem~\ref{theoremclosenessnewtondecrements}, this implies in particular that $\wtilde \lambda_f(x)^2 \less (1+\varepsilon) \lambda_f(x)^2$, and we further obtain that $\lambda_f(x)^2 > \frac{\delta}{(1+\varepsilon) d}$, i.e.,
\begin{align*}
\label{eqnsomeineqforlemma3}
    \varepsilon^\prime < \varepsilon \, \lambda_f(x)^\tau\,.
\end{align*}
There are two possible events.
\begin{itemize}
    \item $E_1$: Either $\wtilde \lambda_f(x) > \eta$. Using the fact that $\mathcal{E}_{x,m,\varepsilon}$ holds, it follows from Lemma~\ref{lemmadecreasefirstphase} that $f(\xnsk) - f(x) \less -\nu$. 
    \item $E_2$: Or $\wtilde \lambda_f(x) \less \eta$. Using the facts that $\mathcal{E}_{x,m,\varepsilon^\prime}$ holds and that $\varepsilon^\prime < \varepsilon \lambda_f(x)^\tau$, it follows from Lemma~\ref{lemmadecreasesecondphase} that $\lambda_f(\xnsk) \less \alpha(\tau)\, (\lambda_f(x))^{1+\tau}$. Assuming further that the event $\mathcal{E}_{\xnsk,m,\varepsilon^\prime}$ holds, we have according to Lemma~\ref{lemmadecreasingdecrementphase2} that $\wtilde \lambda_f(\xnsk) \less \wtilde \lambda_f(x) \less \eta$ and then 
    \begin{align*}
        \wtilde \lambda_f(\xnsk) \underset{(i)}{\less} \sqrt{1+\varepsilon} \, \lambda_f(\xnsk) &\less \sqrt{1+\varepsilon} \, \alpha(\tau) \, (\lambda_f(x))^{1+\tau}\\ & \underset{(ii)}{\less} \sqrt{1+\varepsilon} \,\alpha(\tau) \, (\wtilde \lambda_f(x)/\sqrt{1-\varepsilon})^{1+\tau}\\
        & = \alpha(\tau,\varepsilon) \, (\wtilde \lambda_f(x))^{1+\tau}\,,
    \end{align*}
    where inequalities (i) and (ii) are immediate consequences of Theorem~\ref{theoremclosenessnewtondecrements}.
\end{itemize}
Hence, conditional on $E_2$ occurs once, then the event $E_2$ occurs $K$ additional times in a row with probability at least $1-K \wtilde p$. According to Lemma~\ref{lemmageometricdecrease}, if $K \gre T_{\tau, \alpha(\tau, \varepsilon), \frac{\delta}{d}}$ then the algorithm terminates. On the other hand, the event $E_1$ can occur at most $N_1$ times. 

In summary, conditional on $m_t \gre \overline m$, the algorithm must terminate within $\overline T$ additional iterations with probability at least $1-\overline{T} \wtilde p = 1 - p_0$, and with final sketch size $m \less 2 \overline m$.

\textbf{Phase 1.} At each iteration, one of the following events must occur:
\begin{align*}
    & e_{1} \defn \{\wtilde \lambda_f(x) > \eta,\,f(\xnsk)-f(x) \less -\nu\}\\
    & e_{2} \defn \{\wtilde \lambda_f(x) \less \eta,\, \wtilde \lambda_f(\xnsk) \less \alpha(\tau, \varepsilon) (\wtilde \lambda_f(x))^{1+\tau}\}\\
    & e_{3} \defn \{m \leftarrow 2 m\}\,.
\end{align*}
Fix any iteration $t \gre 0$, and suppose that the algorithm has not yet terminated. Consider the sequence of events $c_0, \dots, c_{t} \in \{e_1, e_2, e_3\}$ up to time $t$. According to Lemma~\ref{lemmageometricdecrease}, any subsequence of $\{c_j\}_{j=0}^t$ which contains only the event $e_2$ would result in termination of Algorithm~\ref{algorithmadaptive} if its length is greater or equal to $T_{\tau, \alpha(\tau, \varepsilon),\delta/d}+1$. Consequently, any such subsequence must have length smaller or equal to $T_{\tau, \alpha(\tau, \varepsilon),\delta/d}$. Between two consecutive longest subsequences containing only $e_2$, either $e_1$ or $e_3$ occur. The event $e_1$ occurs at most $N_1$ times. By assumption on the choice of $m_0$, once $e_3$ has occurred at least $\mathcal{O}\!\left(\log(\overline{d}_\mu)\right)$ times then the sketch size is greater than $\overline m$. Consequently, there are at most $T_1 \defn \mathcal{O}\!\left( \big(N_1 + \log(\overline{d}_\mu)\big)T_{\tau, \alpha(\tau, \varepsilon),\delta/d} \right)$ iterations before reaching a sketch size $m$ such that $m \gre \overline m$ without termination. In the latter case, we enter Phase 2.

\textbf{Combining Phase 1 and Phase 2.} Combining the two above results, we obtain with probability at least $1-p_0$ that Algorithm~\ref{algorithmadaptive} terminates with a final sketch size $m$ smaller than $2 \overline m$ and within a number of iterations $T$ scaling as
\begin{align*}
    T = T_1 + T_2 = \mathcal{O}\!\left( \big(N_1 + \log(\overline{d}_\mu)\big)T_{\tau, \alpha(\tau, \varepsilon),\delta/d} \right) = \mathcal{O}\!\left(\log(\overline{d}_\mu)\, T_{\tau, \alpha(\tau, \varepsilon),\delta/d} \right)\,,
\end{align*}
where the last equality holds by treating $N_1$ as $\mathcal{O}(1)$. 

\textbf{Total complexity.} The worst-case complexity per iteration is given as follows.
\begin{enumerate}[label={(\arabic*)}]
    \item For a SJLT $S$, the sketching cost is at most $\mathcal{O}(nd)$ at each iteration, and forming and solving the linear system $H_S \vnsk = - \gx$ with a direct method using the Woodbury identity takes time $\mathcal{O}(\overline{m}^2 d)$. Multiplying by the number of iterations, we obtain the total time complexity
    \begin{align*}
        \overline{\mathcal{C}} = \mathcal{O}\!\left( \Big(nd + \frac{\overline{d}_\mu^4 d^{2\tau+1} T_{\tau, \alpha(\tau, \varepsilon),\delta/d}^2}{\delta^{2\tau} p_0^2}\Big) \, \log(\overline{d}_\mu)\, T_{\tau, \alpha(\tau, \varepsilon),\delta/d}\right)\,.
    \end{align*}
    For $\tau \approx 1$, we have that $T_{\tau, \alpha(\tau, \varepsilon),\delta/d} = \mathcal{O}(\log(\log(d/\delta)))$. For $n \gtrsim \frac{\overline{d}_\mu^4 d^{2} \log(\log(d/\delta))^2}{\delta^{2} p_0^2}$, the memory and time complexities simplify to
    \begin{align*}
        \overline{m} = \Omega\!\left( \frac{d \overline{d}^2_\mu \log(\log(d/\delta))}{p_0 \delta} \right)\,,\qquad \overline{\mathcal{C}} = \mathcal{O}\!\left( nd \, \log(\overline{d}_\mu) \, \log(\log(d/\delta))) \right)\,.
    \end{align*}
    For $\tau = 0$, we have $T_{\tau, \alpha(\tau, \varepsilon),\delta/d} = \mathcal{O}(\log(d/\delta))$. For $n \gtrsim \frac{\overline{d}_\mu^4 \log(d/\delta)^2}{p_0^2}$, the memory and time complexities simplify to
    \begin{align*}
        \overline{m} = \Omega\!\left( \frac{\overline{d}^2_\mu \log(d/\delta)}{p_0} \right)\,,\qquad \overline{\mathcal{C}} = \mathcal{O}\!\left( nd \, \log(\overline{d}_\mu) \, \log(d/\delta)) \right)\,.
    \end{align*}
    \item We assume for simplicity that $\overline{d}_\mu \gtrsim \log^2(\log(d/\delta))$. For the SRHT, the sketching cost is $\mathcal{O}(nd \, \log \overline m)$, whereas forming and solving the Newton linear system takes time $\mathcal{O}(\overline{m}^2 d)$. Thus, the total complexity is given by
    \begin{align*}
        \overline{\mathcal{C}} = \mathcal{O}\!\left( \left(nd \log \overline{m} + d \, \overline{m}^2\right) \log(\overline{d}_\mu) \, T_{\tau, \alpha(\tau, \varepsilon),\delta/d}\right)\,.
    \end{align*}
    For $\tau \approx 1$, we have $T_{\tau, \alpha(\tau,\varepsilon), \delta/d} = \mathcal{O}(\log(\log(d/\delta)))$. Picking $p_0 \asymp 1/\overline{d}_\mu$, we obtain the memory complexity
    \begin{align*}
        \overline{m} \asymp \frac{d}{\delta} \left(\overline{d}_\mu + \log(d/\delta)\log(\overline{d}_\mu)\right)\,.
    \end{align*}
    Consequently, $\log \overline{m} \lesssim \log(d/\delta)$ and $\overline{m}^2 \lesssim \frac{d^2}{\delta^2} (\overline{d}_\mu^2 + \log^2(d/\delta) \log^2(\overline{d}_\mu))$. Hence, provided that $n \gtrsim \frac{d^2 \overline{d}_\mu^2}{\delta^2}$, we obtain 
    \begin{align*}
        \overline{\mathcal{C}} = \mathcal{O}\!\left( nd \log(d/\delta) \log(\overline{d}_\mu) \log(\log(d/\delta))\right)\,.
    \end{align*}
    For $\tau = 0$, we have $T_{\tau, \alpha(\tau,\varepsilon), \delta/d} = \mathcal{O}(\log(d/\delta))$. Picking $p_0 \asymp 1/\overline{d}_\mu$, we obtain the memory complexity
    \begin{align*}
        \overline{m} \asymp \overline{d}_\mu\,.
    \end{align*}
    Consequently, $\log \overline{m} \lesssim \log(\overline{d}_\mu)$ and $\overline{m}^2 \lesssim \overline{d}_\mu^2$.
    Assuming that $n \gtrsim \overline{d}_\mu^2/\log(\overline{d}_\mu)$, the total time complexity is
    \begin{align*}
        \overline{\mathcal{C}} = \mathcal{O}\Big(nd \,\log(\overline{d}_\mu)^2 \, \log(d/\delta)\Big)\,.
    \end{align*}    
\end{enumerate}
This concludes the proof.\qed

\section{Auxiliary results}

\blems 
\label{lemmadecreasingdecrementphase2}
Let $x \in \mathbf{dom}\,f$ and $\varepsilon \in (0,1/4)$. Suppose that the event $\mathcal{E}_{x,m,\varepsilon} \cap \mathcal{E}_{\xnsk, m_\mathrm{nsk}, \varepsilon}$ holds, and that $\wtilde \lambda_f(x) \less \eta$. Then, we have that
\begin{align}
    \wtilde \lambda_f(\xnsk) \less \wtilde \lambda_f(x) \less \eta\,.
\end{align}
\elems 
\spro 
By assumption, the event $\mathcal{E}_{\xnsk, m_\mathrm{nsk}, \varepsilon}$ holds. It follows from Theorem~\ref{theoremclosenessnewtondecrements} that $\wtilde \lambda_f(\xnsk) \less \sqrt{1+\varepsilon} \, \lambda_f(\xnsk)$. We have by assumption that $\mathcal{E}_{x,m,\varepsilon}$ holds and that $\wtilde \lambda_f(x) \less \eta$.  As a consequence of Lemma~\ref{lemmadecreasesecondphase}, we have $\wtilde \lambda_f(x) \less \frac{16}{25} \, \lambda_f(x)$. As a consequence of Theorem~\ref{theoremclosenessnewtondecrements}, we have $\lambda_f(x) \less \frac{1}{\sqrt{1-\varepsilon}} \, \wtilde \lambda_f(x)$. Combining these bounds together, we obtain that
\begin{align*}
    \wtilde \lambda_f(\xnsk) \less \sqrt{\frac{1+\varepsilon}{1-\varepsilon}} \, \frac{16}{25} \, \wtilde \lambda_f(x)\,. 
\end{align*}
Finally, using that $\varepsilon \in (0,1/4)$, we get that $\sqrt{\frac{1+\varepsilon}{1-\varepsilon}} \, \frac{16}{25} \less 1$, whence,
\begin{align*}
    \wtilde \lambda_f(\xnsk) \less \wtilde \lambda_f(x) \less \eta\,.
\end{align*}
\fpro 
\blems 
\label{lemmaboundeta}
For $\varepsilon \in (0,1)$, it holds that
\begin{align}
    \eta \less \frac{1-\varepsilon}{1+\varepsilon} \, \frac{1}{16} \less \frac{1}{16}\,.    
\end{align}
\elems 
\spro 
Set $\gamma = \left(\frac{1+\varepsilon}{1-\varepsilon}\right)^2$. We aim to show that $\eta \, \sqrt{\gamma} \less 1/16$. Plugging-in the definition of $\eta$ and using that $a \gre 0$, we have $\eta \, \sqrt{\gamma} = \frac{1}{8} \, \frac{1-\frac{\gamma}{2} - a}{\gamma} \less \frac{1}{8} \, \frac{1 - \frac{\gamma}{2}}{\gamma}$. Since $\varphi(\gamma) \defn \frac{1}{8} \, \frac{1-\frac{\gamma}{2}}{\gamma}$ is monotone decreasing and since $\gamma \gre 1$, we obtain that $\eta \, \sqrt{\gamma} \less \varphi(1)$, i.e., $\eta \, \sqrt{\gamma} \less \frac{1}{16}$.
\fpro

\subsection{Technical lemmas for the proof of Theorem~\ref{theoremquadraticconvergenceeffdimnewtonsketch}}

\blems[Phase 1]  
\label{lemmalengthofphase1}
It holds that
\begin{align*}
    \boxed{t \less N_1\,,\quad \mbox{with probability at least } 1-N_1 \wtilde p\,.}\
\end{align*}
\elems 
\spro 
Let $j < t$ be any iteration before $t_1$. Note by construction of Algorithm~\ref{algeffdimnewtonsketch} that $m_j = \overline{m}_1$. Assuming that the event $\mathcal{E}_{x_j, m_j, \varepsilon}$ holds true, it follows from Lemma~\ref{lemmadecreasefirstphase} that we observe the decrement $f(\xnsk) - f(x_j) \less -\nu$. Consequently, under the event $\mathcal{E}^{(1)} \defn \bigcap_{j=0}^{t-1} \mathcal{E}_{x_j, m_j, \varepsilon}$, we obtain that
\begin{align*}
    f(x^*) - f(x_0) \less f(x_{t}) - f(x_0) = \sum_{j=0}^{t-1} f(x_{j+1}) - f(x_j) \less -{t} \, \nu\,.
\end{align*}
Hence, under $\mathcal{E}^{(1)}$, we must have $t \less \frac{f(x_0) - f(x^*)}{\nu}$, i.e., $t \less N_1$. According to Lemma~\ref{lemmaconcentration} and the choice of $\overline{m}_1$, each event $\mathcal{E}_{x_j, m_j, \varepsilon}$ holds with probability at least $1-\wtilde p$. Using a union bound, the event $\mathcal{E}^{(1)}$ holds with probability at least $1-N_1 \wtilde p$.
\fpro 
\blems[Phase 2]
\label{lemmaphase2theorem2}
Under the assumption that $\mathcal{E}^{(2)}$ holds, we have for any $j = 0, \dots, J$ that
\begin{align*}
\begin{cases}
    m_{t+1+j} = \overline{m}_2\,,\\
    \wtilde \lambda_f(x_{t+1+j}) \less \eta\,,\\
    \alpha(\tau)^{\frac{1}{\tau}}\lambda_f(x_{t+1+j+1}) \less (\alpha(\tau)^{\frac{1}{\tau}} \lambda_f(x_{t+1+j}))^{1+\tau}\,.
\end{cases}
\end{align*}
\elems 
\spro 
We prove this claim by induction. We start with $j=0$. By definition of the time $t$, we have $\wtilde \lambda_f(x_t) \less \eta$. Therefore, by construction of Algorithm~\ref{algeffdimnewtonsketch}, we have $m_{t+1} = \overline{m}_2$. From Lemma~\ref{lemmadecreasingdecrementphase2} and under $\mathcal{E}^{(2)}$, we get that $\wtilde \lambda_{f}(x_{t+1}) \less \wtilde \lambda_{f}(x_{t}) \less \eta$. Furthermore, before termination, we have that $\wtilde \lambda_f(x_{t+1})^2 > \frac{3}{4} \delta$. It follows from Theorem~\ref{theoremclosenessnewtondecrements} that 
\begin{align*}
    \lambda_f(x_{t+1})^2 \gre \frac{1}{1+\varepsilon} \wtilde \lambda_f(x_{t+1})^2 > \frac{3}{4(1+\varepsilon)} \delta = \frac{2}{3} \delta\,,
\end{align*}
and this implies in particular that $\varepsilon \delta^{\tau/2} \less \varepsilon (\frac{3}{2})^{\tau/2} \lambda_f(x_{t+1})^\tau \less 2\varepsilon \lambda_f(x_{t+1})^\tau$. Consequently, the hypotheses of Lemma~\ref{lemmadecreasesecondphase} are verified and we have $\alpha(\tau)^{\frac{1}{\tau}}\lambda_f(x_{t+2}) \less (\alpha(\tau)^{\frac{1}{\tau}} \lambda_f(x_{t+1}))^{1+\tau}$.

Now, we prove the induction hypothesis for any $j = 1, \dots, J$, assuming that it holds for $j-1$. Since $\wtilde \lambda_f(x_{t+1+j-1}) \less \eta$, it follows by construction of Algorithm~\ref{algeffdimnewtonsketch} that $m_{t+1+j} = \overline{m}_2$. From Lemma~\ref{lemmadecreasingdecrementphase2} and under $\mathcal{E}^{(2)}$, we get that $\wtilde \lambda_f(x_{t+1+j}) \less \wtilde \lambda_f(x_{t+1+j-1}) \less \eta$. Furthermore, before termination, we have $\wtilde \lambda_f(x_{t+1+j})^2 > \frac{3}{4} \delta$. It follows from Theorem~\ref{theoremclosenessnewtondecrements} that 
\begin{align*}
    \lambda_f(x_{t+1+j})^2 \gre \frac{1}{1+\varepsilon} \wtilde \lambda_f(x_{t+1+j})^2 > \frac{3}{4(1+\varepsilon)} \delta = \frac{2}{3} \delta\,,
\end{align*}
and this implies in particular that $\varepsilon \delta^{\tau/2} \less \varepsilon (\frac{3}{2})^{\tau/2} \lambda_f(x_{t+1+j})^\tau \less 2\varepsilon \lambda_f(x_{t+1+j})^\tau$. Consequently, the hypotheses of Lemma~\ref{lemmadecreasesecondphase} are verified and we have $\alpha(\tau)^{\frac{1}{\tau}}\lambda_f(x_{t+1+j+1}) \less (\alpha(\tau)^{\frac{1}{\tau}} \lambda_f(x_{t+1+j}))^{1+\tau}$.
\fpro 
\bcors 
\label{corollaryprobae2}
The event $\mathcal{E}^{(2)}$ holds true with probability at least $1-(2+T_{\tau,\frac{3}{8}\delta})\wtilde p$.
\ecors 
\spro 
Recall that $m_t = \overline{m}_1$ by definition of the time $t$. According to Lemma~\ref{lemmaphase2theorem2}, if $\mathcal{E}^{(2)}$ holds true, then $m_{t+1+j} = \overline{m}_2$ for $j=0,\dots,J$. From Lemma~\ref{lemmaconcentration}, we have that $\mathbb{P}(\mathcal{E}_{x_t, \overline{m}_1, \varepsilon}) \gre 1-\wtilde p$ and $\mathbb{P}(\mathcal{E}_{x_{t+1+j}, \overline{m}_2, \varepsilon \delta^{\tau/2}}) \gre 1-\wtilde p$. We obtain by a union bound that $\mathbb{P}(\mathcal{E}^{(2)}) \gre 1-(2+T_{\tau, \frac{3}{8}\delta})\wtilde p$.
\fpro

\end{document}